\pdfoutput=1
\PassOptionsToPackage{unicode}{hyperref}
\PassOptionsToPackage{hyphens}{url}
\PassOptionsToPackage{dvipsnames,svgnames,x11names}{xcolor}
\documentclass[
  12pt,
  letterpaper,
]{article}
\usepackage{amsmath,amssymb}
\usepackage{setspace}
\usepackage{iftex}
\ifPDFTeX
  \usepackage[T1]{fontenc}
  \usepackage[utf8]{inputenc}
  \usepackage{textcomp} 
\else 
  \usepackage{unicode-math} 
  \defaultfontfeatures{Scale=MatchLowercase}
  \defaultfontfeatures[\rmfamily]{Ligatures=TeX,Scale=1}
\fi
\usepackage{lmodern}
\ifPDFTeX\else
\fi
\IfFileExists{upquote.sty}{\usepackage{upquote}}{}
\IfFileExists{microtype.sty}{
  \usepackage[]{microtype}
  \UseMicrotypeSet[protrusion]{basicmath} 
}{}
\makeatletter
\@ifundefined{KOMAClassName}{
  \IfFileExists{parskip.sty}{%
    \usepackage{parskip}
  }{
    \setlength{\parindent}{0pt}
    \setlength{\parskip}{6pt plus 2pt minus 1pt}}
}{
  \KOMAoptions{parskip=half}}
\makeatother
\usepackage{xcolor}
\usepackage[margin=1in]{geometry}
\usepackage{longtable,booktabs,array}
\usepackage{calc} 
\usepackage{etoolbox}
\makeatletter
\patchcmd\longtable{\par}{\if@noskipsec\mbox{}\fi\par}{}{}
\makeatother
\IfFileExists{footnotehyper.sty}{\usepackage{footnotehyper}}{\usepackage{footnote}}
\makesavenoteenv{longtable}
\usepackage{graphicx}
\makeatletter
\def\maxwidth{\ifdim\Gin@nat@width>\linewidth\linewidth\else\Gin@nat@width\fi}
\def\maxheight{\ifdim\Gin@nat@height>\textheight\textheight\else\Gin@nat@height\fi}
\makeatother
\setkeys{Gin}{width=\maxwidth,height=\maxheight,keepaspectratio}
\makeatletter
\def\fps@figure{htbp}
\makeatother
\providecommand{\tightlist}{%
  \setlength{\itemsep}{0pt}\setlength{\parskip}{0pt}}
\usepackage{tcolorbox}
\tcbuselibrary{breakable,skins}

\usepackage{titlesec}
\titleformat{\section}{\Large\bfseries\sffamily}{}{0em}{}
\titleformat{\subsection}{\large\bfseries\sffamily}{}{0em}{}
\titleformat{\subsubsection}{\normalsize\bfseries\sffamily}{}{0em}{}
\titlespacing*{\section}{0pt}{18pt plus 4pt minus 2pt}{8pt plus 2pt}
\titlespacing*{\subsection}{0pt}{12pt plus 3pt minus 2pt}{6pt plus 2pt}
\usepackage{caption}
\usepackage{fancyhdr}
\allowdisplaybreaks
\usepackage{float}
\makeatletter
\renewcommand{\and}{\end{tabular}\par\vskip 0.6em\begin{tabular}[t]{c}}
\makeatother
\ifLuaTeX
  \usepackage{selnolig}  
\fi
\IfFileExists{bookmark.sty}{\usepackage{bookmark}}{\usepackage{hyperref}}
\IfFileExists{xurl.sty}{\usepackage{xurl}}{} 
\hypersetup{
  pdftitle={A Unified Theory of Scale},
  pdfauthor={Xiao Huang},
  colorlinks=true,
  linkcolor={blue!50!black},
  filecolor={Maroon},
  citecolor={Blue},
  urlcolor={blue!50!black},
  pdfcreator={LaTeX via pandoc}}

\title{A Unified Theory of Scale}
\author{Xiao Huang\\
Department of Environmental Sciences, Emory University\\
Atlanta, GA 30329\\
Email: xiao.huang2@emory.edu}
\date{November 2026}

\begin{document}
\maketitle
\begin{abstract}
Spatial data are often released by zone: a fine-scale field is averaged
within zones with respect to a support measure (area, population or
sampling weights). The modifiable areal unit problem, the ecological
fallacy, Simpson's paradox and the change-of-support problem all concern
one question: to what extent do zonal statistics reflect fine-scale
relationships? We formulate within-zone averaging as the conditional
expectation \(\mathbb{E}[\,\cdot\mid\mathcal{G}]\) and a zoning as a
sub-\(\sigma\)-algebra \(\mathcal{G}\), and prove three groups of
results in this framework. (A) Data-dependent zonings: if the support
measure is atomless and the two variables are not exactly linearly
related, then for any number of zones \(k\ge3\) and any prescribed zone
weights, the zonal correlation attains every value in \([-1,1]\)
exactly; requiring the zones of a square domain to be connected still
allows every value to be approached arbitrarily closely, so the
fine-scale correlation places no constraint on the zonal correlation.
(B) Compact zonings: the difference between the zonal and the fine-scale
correlation \(\rho\) is bounded explicitly in terms of \(\rho\) and the
two within-zone variance shares, and the bound is attained whenever it
does not exceed \(1-|\rho|\); for Lipschitz fields the bound is of
second order in the ratio of zone diameter to variation length, so
halving the largest zone diameter reduces the leading term of the bound
to one quarter. (C) Data-independent zoning designs: for jointly weakly
stationary fields, the design-level zonal correlation equals a weighted
average of the real coherence over spatial frequencies times a factor
not exceeding 1; it takes the same non-zero value for all connected zone
shapes and their mixtures if and only if the spectral measures of the
two variables are proportional and the real part of the cross-spectrum
is a non-zero multiple of this common spectrum (the intrinsic
correlation model); when proportionality fails and the fine-scale
correlation is non-zero, some design made of at most two connected zone
shapes moves the zonal correlation away from the fine-scale correlation,
while the correlations on any finite set of designs cannot certify
proportionality; as zones grow without bound, the zonal correlation
converges, for square-integrable kernels and spectral densities that are
bounded, continuous and positive at the origin, to the zero-frequency
coherence, whatever the zone shape. We also give a martingale
scale-spectrum decomposition along nested hierarchies together with its
renormalization fixed points, identification intervals for the
fine-scale covariance from zonal data alone that are sharp under an
atomless support measure, and a necessary and sufficient condition for
linear macro interventions to be well defined. All results are proved in
full, including the Sierpiński theorem, Lyapunov's convexity theorem,
the Bochner--Cramér spectral representation and the uniqueness theorem
for Fourier--Stieltjes transforms on which they rely, and the main
results are checked numerically.

\textbf{Keywords:} modifiable areal unit problem; ecological inference;
change of support; conditional expectation; convex order; Lyapunov
convexity theorem; spectral measure; intrinsic correlation model
\end{abstract}

\renewcommand*\contentsname{Contents}
{
\hypersetup{linkcolor=}
\setcounter{tocdepth}{2}
\tableofcontents
}
\setstretch{1.5}
\hypertarget{introduction}{%
\section{1 Introduction}\label{introduction}}

\hypertarget{the-problem}{%
\subsection{1.1 The problem}\label{the-problem}}

Many spatial data are available only as zonal aggregates. Population and
health statistics are released by census tract, county or postal code,
environmental monitoring and remote-sensing products by grid cell, and
survey data by the geographic units of a sampling frame. Behind such
data lies a field defined at a fine scale: a function defined pointwise
on a continuous domain, the pixel values of a raster, or individual
records. A zonal value is the average of this field over the zone with
respect to some support measure: area for a continuous field, population
counts for a demographic rate, sampling weights for a survey estimate.
Analysts then compute correlation coefficients, regression slopes and
other statistics from the zonal values and use them to draw conclusions
about fine-scale relationships.

The relation between zonal and fine-scale statistics has been studied in
four literatures. In geography it is the modifiable areal unit problem
(MAUP): Gehlke \& Biehl (1934) noticed that merging census tracts
changes correlation coefficients, and Openshaw \& Taylor (1979) and
Openshaw (1984) systematized the ``zoning effect'' of redrawing
boundaries with the number of zones held fixed and the ``scale effect''
of changing zone size. In sociology and epidemiology it is the
ecological fallacy and ecological inference: Robinson (1950) showed that
zonal and individual-level correlations can differ greatly, Goodman
(1953) and Duncan \& Davis (1953) proposed ecological regression and
deterministic bounds respectively, and King (1997) and Wakefield (2008)
discussed the problem systematically. Simpson's (1951) paradox is the
same phenomenon in contingency tables. Geostatistics treats it as the
change-of-support problem: the covariance of block averages is obtained
by integrating the point covariance, and block kriging gives predictable
and correctable changes of scale for Gaussian fields (Journel \&
Huijbregts 1978; Cressie 1993; Gelfand, Zhu \& Carlin 2001; Gotway \&
Young 2002; Chilès \& Delfiner 2012).

These strands leave inconsistent impressions. Zoning experiments show
that the correlation coefficient of one data set can take almost any
value under different zonings (Openshaw \& Taylor 1979), and in
multivariate regression the direction and strength of the MAUP are hard
to predict (Fotheringham \& Wong 1991); the geostatistical theory of
change of support, in contrast, says that the second moments of block
averages can be computed exactly once the covariance structure is known.
The question addressed in this paper is: what mathematical structure
underlies these phenomena, and under what conditions, and to what
extent, can zonal statistics be trusted?

\hypertarget{what-is-known-and-what-is-missing}{%
\subsection{1.2 What is known and what is
missing}\label{what-is-known-and-what-is-missing}}

Three facts are common knowledge in spatial statistics: redrawing
boundaries changes zonal correlations; smaller zones lose less
within-zone variation; changing the aggregation scale changes the
contributions of the different spatial scales to the zonal correlation.
What has been missing are the quantitative statements and the necessary
and sufficient conditions behind these qualitative facts.

\textbf{How large can the zoning effect be?} Redrawing boundaries
changes how the total variation is split into within-zone and
between-zone parts, but does the fine-scale correlation still constrain
the zonal correlation? If not, can fixing the number of zones and the
zone weights, or even requiring every zone to be connected, prevent an
arbitrary zonal correlation from being manufactured?

\textbf{How well can the aggregation error be controlled?} ``Smaller
zones, smaller error'' only gives a direction. What is needed is an
error bound that does not depend on a model and can be computed directly
from fine-scale data, the conditions under which it is attained, and,
under smoothness conditions, the rate at which the bound vanishes as
zone diameters shrink and whether that rate can be improved.

\textbf{Under what condition is the zonal correlation invariant to scale
and zone shape?} ``Change the scale and the correlation may change'' is
likewise only a direction. What is needed is a necessary and sufficient
condition for invariance, and conversely: when the condition fails, must
some connected zoning design move the zonal correlation away from the
fine-scale correlation, and can correlations observed at finitely many
scales certify invariance?

Two related questions also arise: to what extent are fine-scale
relationships identified from zonal data alone, and when are the effects
of macro interventions defined if zone means are treated as macro
variables?

Figure 1 illustrates the size of the zoning effect with a simulated
example. With the same fine-scale field and the same 400 basic units
grouped into 20 connected zones, the zonal correlation ranges from
\(-0.95\) to \(+1.00\), while the fine-scale correlation is only
\(+0.35\).

\begin{figure}
\centering
\includegraphics[width=1\textwidth,height=\textheight]{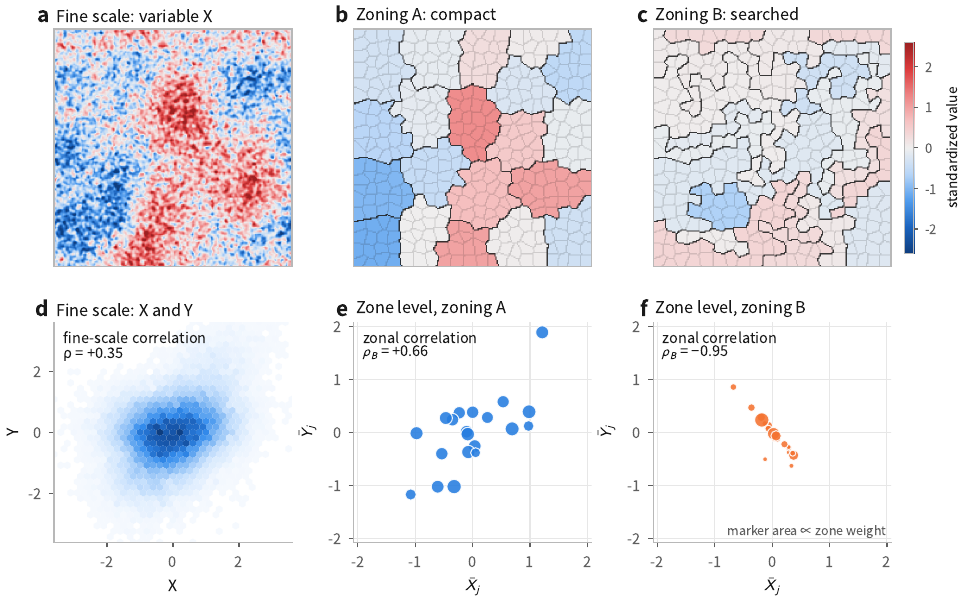}
\caption{The same fine-scale data, two zonings. (a) The fine-scale
variable \(X\) (\(200\times200\) pixels, standardized values); the
fine-scale correlation between \(Y\) and \(X\) is \(+0.35\). The study
region is tessellated into 400 irregular basic units (the thin grey
lines in (b) and (c)), and each zoning groups these units into 20 zones
(thick black lines). (b) Zoning A: compact zones colored by zone mean,
i.e.~by \(\mathbb{E}[X\mid\mathcal{G}]\). (c) Zoning B: obtained by a
local search over unit assignments that minimizes the zonal correlation
while keeping every zone connected. (d) Fine-scale \(X\) and \(Y\)
(hexagonal-bin density). (e) Zonal data under zoning A: zonal
correlation \(\rho_B=+0.66\). (f) Zonal data under zoning B:
\(\rho_B=-0.95\); marker areas are proportional to zone weights.
Grouping the same 400 units into 20 connected zones yields zonal
correlations anywhere between \(-0.95\) and \(+1.00\).}
\end{figure}

\hypertarget{approach}{%
\subsection{1.3 Approach}\label{approach}}

Our starting point is a simple observation: replacing the value of each
fine-scale unit by the weighted average of its zone produces exactly the
conditional expectation \(\mathbb{E}[Z\mid\mathcal{G}]\), where the
sub-\(\sigma\)-algebra \(\mathcal{G}\) is generated by the zoning, and
the weighted distribution of the zone means, in which each zone mean
carries the weight of its zone, is exactly the distribution of
\(\mathbb{E}[Z\mid\mathcal{G}]\) (Lemma 3.1). Every question about zonal
statistics computed from this weighted distribution (weighted means,
variances, correlations, regression slopes) thereby becomes a question
about conditional expectations and their distributions, and
well-developed tools of probability and harmonic analysis become
available:

\begin{itemize}
\tightlist
\item
  the tower property describes nested hierarchies of zonings, and the
  conditional Jensen inequality describes the contraction of
  distributions in the convex order caused by aggregation (Section 3);
\item
  Lyapunov's convexity theorem on atomless measures characterizes all
  zonal configurations that zonings can realize (Section 4), which
  yields the exact range of the zoning effect (Sections 5 and 6);
\item
  the within/between decomposition of variance gives exact error bounds
  for compact zonings (Section 7);
\item
  the spectral representation of stationary fields writes block
  averaging as low-pass filtering, so that the zonal correlation of a
  data-independent zoning design becomes a weighted average of the
  coherence (Section 8);
\item
  the martingale-difference decomposition distributes covariance
  orthogonally across the levels of a nested hierarchy (Section 9).
\end{itemize}

In this framework, how far zonal statistics can be trusted depends on
only three things: whether the zoning depends on the data; whether the
zones are compact; and how the coherence of the two variables is
distributed across spatial scales.

\hypertarget{main-results}{%
\subsection{1.4 Main results}\label{main-results}}

\textbf{Main Theorem (trichotomy).} In (A) and (B) let \(Z=(X,Y)\) be
square-integrable; in (A) let moreover the support measure \(\mu\) be
atomless and \(aX+bY\) not almost surely constant for any \((a,b)\neq0\)
((B) holds for every support measure). Denote by \(\rho\) the fine-scale
correlation and by \(\rho_B\) the zonal correlation under a zoning
\(\mathcal{G}\) (the correlation of zone means weighted by \(\mu\)).
Part (C) concerns jointly weakly stationary random fields and zoning
designs (Section 8).

\begin{enumerate}
\def\labelenumi{(\Alph{enumi})}
\item
  \textbf{Data-dependent zonings: complete indeterminacy.} For any
  \(k\ge3\), any prescribed zone weights \(p=(p_1,\dots,p_k)\)
  (\(p_j>0\), \(\sum_jp_j=1\)) and any \(\rho^*\in[-1,1]\) there is a
  \(k\)-zone zoning with zone weights exactly \(p\) and
  \(\rho_B=\rho^*\); for \(k=2\), \(\rho_B\) can only be \(\pm1\), and
  both values occur. On a planar square domain (\(\mu\) with a density),
  the conclusion still holds to arbitrary precision when every zone is
  required to be a path-connected polygon. The fine-scale correlation
  places no constraint on the zonal correlation (Theorem 6.2, Corollary
  6.3, Theorem 6.5, Corollary 6.6).
\item
  \textbf{Compact zonings: an exact bound and a second-order rate.} With
  \(a,b<1\) the within-zone variance shares of the two variables and
  \(s=\sqrt{(1-a)(1-b)}\),
  \(|\rho_B-\rho|\le\big(|\rho|(1-s)+\sqrt{ab}\big)/s\), and the bound
  is attained whenever it does not exceed \(1-|\rho|\). If \(X,Y\) are
  Lipschitz fields and every zone has diameter at most \(\delta\), put
  \(\tau=\delta L/(2\sigma)\) (\(\sigma/L\) is the variation length);
  then
  \(|\rho_B-\rho|\le\tau_X\tau_Y+\tfrac{|\rho|}{2}(\tau_X^2+\tau_Y^2)+O(\tau^4)\),
  and halving \(\delta\) reduces the leading term to one quarter. For
  stationary random fields and a zone shape scaled by \(T\), the
  design-level within-zone variance share likewise vanishes
  quadratically in \(T\) when the fields are mean-square differentiable
  (for kernels with finite second moments), but only at the order
  \(T^\alpha\) when \(C(0)-C(h)\sim c|h|^\alpha\) with \(\alpha<2\) (for
  kernels of bounded support), for example linearly for the exponential
  covariance (Theorem 7.1, Corollaries 7.2 and 7.3, Propositions 7.4 and
  8.16).
\item
  \textbf{Data-independent zoning designs: a necessary and sufficient
  condition.} For jointly weakly stationary random fields on
  \(\mathbb{R}^q\) and zoning designs fixed independently of them, the
  design-level zonal correlation equals a weighted average of the real
  coherence times a factor not exceeding 1, so it always lies in
  \([\min(0,\operatorname{ess\,inf}\operatorname{coh}),\)
  \(\max(0,\operatorname{ess\,sup}\operatorname{coh})]\). It takes the
  same non-zero value \(\rho_0\) for all connected zone shapes and their
  finite mixtures if and only if \(F_{YY}=\lambda F_{XX}\) and
  \(\operatorname{Re}F_{XY}=\rho_0\sqrt\lambda\,F_{XX}\) (the intrinsic
  correlation model), and it vanishes for all of them if and only if
  \(\operatorname{Re}F_{XY}=0\); if the intrinsic correlation model
  fails and the fine-scale correlation is non-zero, some design made of
  at most two connected zone shapes gives a zonal correlation different
  from the fine-scale correlation; and the correlations on finitely many
  designs cannot certify the intrinsic correlation model. As zones grow
  without bound, the zonal correlation converges to the zero-frequency
  coherence \(\operatorname{coh}(0)\), whatever the zone shape, provided
  the kernel is square-integrable and the spectral densities are bounded
  and continuous at the origin with \(f_{XX}(0)f_{YY}(0)>0\) (Theorems
  8.9 and 8.11, Corollary 8.12, Proposition 8.13, Theorem 8.15).
\end{enumerate}

The three parts are proved independently: (A) and the support-measure
statements of (B) rest on the conditional-expectation formulation of
Section 3, (C) and the stationary-field statement of (B) (Proposition
8.16) on the spectral representation of Section 8. The contrast between
(A) and (C) is the heart of the paper. For the same pair of variables,
once the zoning may depend on the data, any zonal correlation can be
manufactured; when the zoning does not depend on the data, the
design-level zonal correlation is held firmly within the limits set by
the coherence, and it takes the same non-zero value at all scales and
zone shapes only under the intrinsic correlation model.

Beyond the Main Theorem, the contributions include:

\begin{enumerate}
\def\labelenumi{\arabic{enumi}.}
\tightlist
\item
  \textbf{Only means are invariant to zoning (Theorem 5.5).} Under an
  atomless support measure, a zonal statistic that takes the same value
  for all \(k\)-zone zonings must take the value it has when every zone
  mean equals the overall mean, and a statistic invariant on every
  atomless support measure space is a function of the overall mean.
  Under the non-degeneracy condition of the Main Theorem, zonal
  variances, correlations and regression slopes all change with the
  zoning (Corollaries 6.3 and 6.4), so dispersion and association
  statistics ``immune to zoning'' do not exist.
\item
  \textbf{A necessary and sufficient condition for aggregation bias
  (Theorem 8.14).} For the linear model \(Y=\alpha+\beta X+e\), the
  design-level aggregated regression slope equals \(\beta\) for every
  connected design if and only if the symmetrized lagged
  cross-covariance of \(X\) and \(e\) vanishes everywhere; absence of
  correlation at the same location does not guarantee unbiased
  aggregation.
\item
  \textbf{Scale spectra along nested hierarchies (Theorems 9.2--9.4,
  Propositions 9.6--9.7).} Covariance decomposes orthogonally across the
  levels of a nested hierarchy; the zonal correlation at a level depends
  only on that level and the coarser ones; ``coarser aggregation, higher
  correlation'' holds when the scale spectra are proportional and the
  level coherence is non-decreasing from finer to coarser levels, and
  can fail when the scale spectra are not proportional, even if the
  level coherence increases; the fixed points of the coarse-graining map
  are exactly the geometric (power-law) scale spectra.
\item
  \textbf{Identification intervals from zonal data alone (Theorem 10.1,
  Corollary 10.2).} For variables taking values in \([0,1]\), the
  fine-scale covariance lies in an interval computable from zonal data,
  which is sharp when the support measure is atomless; its width can
  serve as an ``ecological fallacy index''.
\item
  \textbf{Well-defined macro interventions (Theorem 11.2).} For linear
  structural equations, the macro intervention that sets a zone mean to
  a given value is well defined if and only if the within-zone effects
  are homogeneous along the permitted intervention directions; for
  uniform within-zone shifts it is always well defined, and the macro
  effect equals the weighted mean of the within-zone effects.
\end{enumerate}

\textbf{Completeness.} Every result labeled theorem, proposition,
corollary or lemma is proved in full. Only standard textbook results of
measure theory, Fourier analysis, functional analysis, probability and
convex analysis are used without proof; they are listed at the end of
Section 3.1. Known results that are re-proved to keep the paper
self-contained include the results collected in Section 3 (the tower
property, the conditional Jensen inequality, the projection property and
their consequences), the Sierpiński theorem (Theorem 4.3), Lyapunov's
convexity theorem (Theorem 4.6), the Dvoretzky--Wald--Wolfowitz theorem
and its fixed-weight form (Theorem 4.8, Corollary 4.10), the law of
total covariance and Robinson's identity (Lemma 5.1, Theorem 5.2), the
uniqueness, Bochner and Cramér theorems (Theorems 8.1--8.3), the
block-covariance formula (Lemma 8.6), the sufficiency direction of
Theorem 8.11, the orthogonality of martingale increments (Lemma 9.1,
Theorem 9.2) and the binary case of Theorem 10.1 (the Duncan--Davis
bounds); Theorem 8.15(ii) echoes Tiao (1972) and Vargas-Guzmán, Warrick
\& Myers (1999). The other results are established in this paper; the
main ones are listed in Section 13.1.

\hypertarget{organization}{%
\subsection{1.5 Organization}\label{organization}}

Section 2 reviews related work. Section 3 sets up the framework: support
measures, zonings and conditional expectation, together with basic
properties of the convex order. Section 4 treats atomless measures and
proves the Sierpiński theorem, Lyapunov's convexity theorem and the
partition version of the Dvoretzky--Wald--Wolfowitz theorem. Section 5
gives the within/between decomposition of variance and the impossibility
theorem. Section 6 deals with data-dependent zonings, Section 7 with
compact zonings, and Section 8 develops the spectral theory of
data-independent zoning designs. Section 9 treats nested hierarchies,
Section 10 identification, and Section 11 macro interventions. Section
12 reports numerical checks, and Section 13 discusses the contributions,
implications for spatial analysis and limitations.

\hypertarget{related-work}{%
\section{2 Related work}\label{related-work}}

\hypertarget{the-modifiable-areal-unit-problem-and-zoning-experiments}{%
\subsection{2.1 The modifiable areal unit problem and zoning
experiments}\label{the-modifiable-areal-unit-problem-and-zoning-experiments}}

Using census-tract data, Gehlke \& Biehl (1934) found that successively
merging adjacent tracts systematically increases correlation
coefficients, one of the earliest reports of the effect of aggregation
on correlation. Openshaw (1977) treated the drawing of zones itself as
an optimization problem, searching over zonings of given basic units for
zonings that meet an objective, and proposed an automated zoning
procedure for this purpose; Openshaw \& Taylor (1979) enumerated a large
number of zonings of a set of county-level data, showed that with the
number of zones fixed the correlation coefficient can vary over a wide
range, and coined the name ``modifiable areal unit problem'',
distinguishing the scale effect from the zoning effect, a distinction
systematized by Openshaw (1984). Later simulation studies showed that
the strength and direction of the MAUP in multivariate regression are
hard to predict (Fotheringham \& Wong 1991), that means are insensitive
to aggregation while regression coefficients and correlation statistics
vary greatly (Amrhein 1995), and that within-zone homogeneity is the key
factor governing aggregation effects (Holt, Steel \& Tranmer 1996).
Arbia (1989) studied the effect of aggregation on statistical analysis
from the standpoint of spatial autocorrelation and data configuration.
Tobler (1989) argued that spatial analysis should use
``frame-independent'' methods, i.e., methods whose conclusions do not
change with the zoning.

These studies established the existence of zoning and scale effects and
their size in real data, but they did not answer two basic questions:
without model assumptions, what is the range of possible zonal
correlations, and are there association statistics that are invariant to
zoning? Section 6 answers the first exactly: under an atomless support
measure and for variables without an exact linear relation, with at
least three zones the attainable range is the whole interval \([-1,1]\),
and when every zone of a square domain must be connected, every value
can still be approached arbitrarily closely. Section 5 answers the
second: under an atomless support measure, a zonal statistic invariant
to zoning must take the value it has when all zone means equal the
overall mean, so for variables without an exact linear relation zonal
variances, correlations and slopes are not invariant.

\hypertarget{ecological-inference-and-the-ecological-fallacy}{%
\subsection{2.2 Ecological inference and the ecological
fallacy}\label{ecological-inference-and-the-ecological-fallacy}}

Using the 1930 US census, Robinson (1950) showed that the correlation
between race and illiteracy was 0.946 across nine geographic divisions
but only 0.203 across individuals, and gave the decomposition relating
zonal and individual-level correlations; in this paper the decomposition
is Theorem 5.2. Goodman (1953) proposed ecological regression: under
certain constancy assumptions on within-zone proportions and effects,
individual-level parameters can be recovered from a zonal regression.
Duncan \& Davis (1953) gave deterministic bounds that need no
assumptions: when only the two marginal proportions of each zone are
known, the individual-level joint proportion lies in a computable
interval. King's (1997) method of ecological inference adds
distributional assumptions within these deterministic bounds to obtain
point estimates. Greenland \& Robins (1994) discussed biases and
counterexamples in ecological studies systematically; Wakefield (2008)
pointed out that ecological bias arises because zonal data cannot
characterize the within-area variability of exposures and confounders.
Gelman et al.~(2007) used multilevel models to show that within-state
and between-state relationships between income and voting can point in
opposite directions.

This paper advances this line in three ways: Theorem 10.1 extends the
Duncan--Davis bounds to variables taking values in \([0,1]\) and proves
that they are sharp under an atomless support measure, the width of the
interval giving an ``ecological fallacy index'' computable from zonal
data alone; Section 7 gives an exact bound on the difference between
zonal and fine-scale correlations when fine-scale information is
available; and Section 11 states the conditions needed to interpret
zonal effects as macro-level causal effects.

\hypertarget{change-of-support-and-multivariate-geostatistics}{%
\subsection{2.3 Change of support and multivariate
geostatistics}\label{change-of-support-and-multivariate-geostatistics}}

Geostatistics calls the spatial extent to which an observation refers
its support; the change-of-support problem concerns inference between
point and block supports. The covariance of block averages equals the
double integral of the point covariance over the two blocks, which is
the basis of dispersion variances and block kriging (Journel \&
Huijbregts 1978; Cressie 1993; Chilès \& Delfiner 2012). Cressie (1996)
showed how a bivariate extension of the geostatistical change-of-support
approach can be used to tackle the MAUP; Gelfand, Zhu \& Carlin (2001)
handled changes of support between points and blocks in Bayesian
hierarchical models; Gotway \& Young (2002) gave a systematic review of
methods for combining incompatible spatial data; Banerjee, Carlin \&
Gelfand (2014) discussed spatial misalignment within hierarchical
models. In multivariate geostatistics, the intrinsic correlation model
(covariance and cross-covariance functions proportional to one another)
is the simplest coregionalization model, and under it correlations
between variables do not change with support (Wackernagel 2003).
Vargas-Guzmán, Warrick \& Myers (1999) studied how correlation
coefficients change with support and with the spatial scales of
variability; in time series, Tiao (1972) studied the limiting behavior
of temporal aggregates. The spectral theory of stationary random fields
can be found in Yaglom (1987) and Stein (1999).

These studies show that second moments of block averages can be computed
exactly once the covariance structure is given. Building on this,
Section 8 answers several open questions: within which limits the zonal
correlation of an arbitrary data-independent zoning design must lie (the
coherence bound of Theorem 8.9); what is the necessary and sufficient
condition for the zonal correlation to be invariant across all connected
zone shapes and their mixtures (Theorem 8.11: for a non-zero common
value the intrinsic correlation model is not only sufficient but also
necessary); when invariance fails, must the deviation be revealed by
some simple design (Corollary 8.12); can finitely many designs certify
invariance (Proposition 8.13); and what does the zonal correlation
converge to as zones grow without bound (Theorem 8.15)?

\hypertarget{aggregation-error-criteria-aggregation-bias-and-spatial-confounding}{%
\subsection{2.4 Aggregation error criteria, aggregation bias and spatial
confounding}\label{aggregation-error-criteria-aggregation-bias-and-spatial-confounding}}

Bradley, Wikle \& Holan (2017) defined a criterion for spatial
aggregation error (CAGE), a between-scale variance computed from
multiscale Karhunen--Loève expansions, and used it to choose
regionalizations; Daw et al.~(2025) call Proposition 2 of Bradley, Wikle
\& Holan (2017) a ``null-MAUP theorem'': when CAGE is zero, the MAUP and
the ecological fallacy do not occur. Daw et al.~extended the criterion
to the multivariate case. In econometrics, Theil (1954) studied the
aggregation bias that arises when micro relations are linearly
aggregated into macro relations. The spatial regression literature has
discussed how the spatial scales of exposure and unmeasured confounding
affect regression coefficients: adding spatially correlated random
effects does not remove unmeasured confounding (Hodges \& Reich 2010),
the bias depends on the spatial scales at which exposure and confounders
vary (Paciorek 2010), and it can be adjusted in the frequency domain
(Guan et al.~2023).

Theorem 8.14 writes aggregation bias as a frequency-domain condition:
the design-level aggregated regression slope is unbiased for every
connected design if and only if the real part of the cross-spectrum of
exposure and error vanishes, equivalently, the symmetrized lagged
cross-covariance vanishes everywhere. This agrees with the scale and
frequency-domain perspective of the literature above and gives an exact
criterion.

\hypertarget{convexity-of-vector-measures-fusions-of-distributions-and-zoning-algorithms}{%
\subsection{2.5 Convexity of vector measures, fusions of distributions
and zoning
algorithms}\label{convexity-of-vector-measures-fusions-of-distributions-and-zoning-algorithms}}

Lyapunov (1940) proved that the range of an atomless finite-dimensional
vector measure is compact and convex, and Lindenstrauss (1966) gave a
short proof using the Krein--Milman theorem. Dvoretzky, Wald \&
Wolfowitz (1951) extended it to partitions: for atomless vector
measures, every ``proportional allocation'' can be realized exactly by a
genuine partition. Fusions of distributions (collapsing part of the mass
to its barycenter) and their relation to the convex order are studied in
Elton \& Hill (1992); the general theory of stochastic orders is in
Shaked \& Shanthikumar (2007). On the side of zone design, automatic
zoning and spatially constrained clustering algorithms (Openshaw 1977;
Duque, Anselin \& Rey 2012) search over zonings of given basic units for
zonings satisfying homogeneity or size constraints; redistricting
research has turned to Markov-chain sampling of large ensembles of
compliant plans under contiguity and other constraints, using the
ensemble as a baseline for comparison (DeFord, Duchin \& Solomon 2021).

Section 4 uses the theorems of Lyapunov and Dvoretzky--Wald--Wolfowitz
as basic tools for zoning problems: under an atomless support measure,
the set of attainable zone weights and zone totals is convex, and so is
the set of attainable zone-mean vectors for prescribed zone weights, so
the combinatorial zoning problem becomes a convex one. The indeterminacy
theorems of Section 6 show that once the zoning may depend on the data,
a connectivity constraint by itself cannot prevent arbitrary zonal
correlations from being manufactured, which gives a mathematical reason
for the ensemble methods used in the evaluation of redistricting plans.

\hypertarget{causal-abstraction-and-urban-scaling-laws}{%
\subsection{2.6 Causal abstraction and urban scaling
laws}\label{causal-abstraction-and-urban-scaling-laws}}

Rubenstein et al.~(2017) and Beckers \& Halpern (2019) studied exact
transformations (abstractions) between micro and macro causal models:
under what conditions a macro model gives conclusions consistent with a
micro model for a family of interventions. Section 11 gives an explicit
form of this condition for linear spatial aggregation. Arcaute et
al.~(2015) showed that, for indicators that scale nonlinearly with city
size, the estimated scaling exponents can change considerably with the
definition of city boundaries, a manifestation of the MAUP in urban
science.

\hypertarget{position-of-this-paper}{%
\subsection{2.7 Position of this paper}\label{position-of-this-paper}}

The literature above established the size of zoning effects, the
deterministic bounds of ecological inference, the computation of second
moments of block averages, the scale mechanisms of aggregation bias and
spatial confounding, and ensemble methods for redistricting. The
contribution of this paper is to unify these problems through a single
mathematical object, the conditional expectation
\(\mathbb{E}[\,\cdot\mid\mathcal{G}]\) and its distribution, and to give
the exact boundaries between them: under data-dependent zonings the
attainable range is the whole interval (Section 6); under compact
zonings the error has an exact bound, attained whenever it does not
exceed one minus the absolute value of the fine-scale correlation, and,
for Lipschitz fields, a second-order rate in the zone diameter that
cannot be improved (Section 7); under data-independent designs the zonal
correlation is confined by the coherence and takes the same non-zero
value at all scales and zone shapes only under the intrinsic correlation
model (Section 8). All results are stated as theorems with complete
proofs.

\hypertarget{setting-and-preliminaries}{%
\section{3 Setting and preliminaries}\label{setting-and-preliminaries}}

\hypertarget{support-measures-fine-scale-variables-and-zonings}{%
\subsection{3.1 Support measures, fine-scale variables and
zonings}\label{support-measures-fine-scale-variables-and-zonings}}

\textbf{Definition (support measure space).} Let
\((\Omega,\mathcal{F},\mu)\) be a probability space. \(\Omega\) is the
set of fine-scale units: points of a domain \(D\subseteq\mathbb{R}^q\),
pixels of a raster, or individual records. \(\mu\) is the normalized
\textbf{support measure}: for a measurable set \(A\), \(\mu(A)\) is the
share of the total mass carried by \(A\). When \(\Omega=\{1,\dots,n\}\)
is finite, \(\mu(\{i\})=w_i/\sum_lw_l\) with unit weights \(w_i>0\).

\textbf{Remark (choice of support measure).} The support measure is
determined by the quantity under study: area or volume (normalized
Lebesgue measure) for continuous fields, the counts in the denominator
for rates (population counts, for example), and sampling weights for
survey estimates. If every zone is weighted equally in an analysis, this
amounts to replacing \(\mu\) by a measure \(\mu'\) that gives every zone
the same total mass; for that fixed zoning the results of this paper
hold for \(\mu'\) as well, with the ``fine-scale correlation'' computed
under \(\mu'\). Since \(\mu'\) depends on the zoning, results that
compare different zonings (for example Theorem 5.5, Sections 6 and 9 and
Corollary 11.4) do not carry over to equal zone weights.

\textbf{Fine-scale variables.} A fine-scale variable is an integrable
random vector \(Z:\Omega\to\mathbb{R}^d\),
i.e.~\(Z\in L^1(\mu;\mathbb{R}^d)\); when variances and correlations are
discussed we assume \(Z\in L^2\). Throughout, expectations, variances,
covariances and correlations are computed under \(\mu\).

\textbf{Definition (zoning).} A \(k\)-zone zoning is a measurable
partition \(\mathcal{Z}=\{A_1,\dots,A_k\}\) of \(\Omega\) with every
\textbf{zone weight} \(p_j=\mu(A_j)>0\). It generates the
sub-\(\sigma\)-algebra \(\mathcal{G}=\sigma(\mathcal{Z})\). The
\textbf{zone mean} of zone \(j\) is
\(\bar z_j=p_j^{-1}\int_{A_j}Z\,d\mu\), the \textbf{zonal data} are the
weighted point set \(\{(p_j,\bar z_j)\}_{j\le k}\), and a \textbf{zonal
statistic} is a statistic computed from the \(\bar z_j\) with weights
\(p_j\). Write
\(\Delta_k^\circ=\{p\in\mathbb{R}^k:p_j>0,\ \sum_jp_j=1\}\). For a
measurable set \(A\) with \(\mu(A)>0\) write
\(\mathbb{E}[Z\mid A]=\mu(A)^{-1}\int_AZ\,d\mu\), and for real
\(X,Y\in L^2\) write
\(\operatorname{Cov}(X,Y\mid A)=\mathbb{E}\big[(X-\mathbb{E}[X\mid A])(Y-\mathbb{E}[Y\mid A])\mid A\big]\)
and \(\operatorname{Var}(X\mid A)=\operatorname{Cov}(X,X\mid A)\); thus
\(\bar z_j=\mathbb{E}[Z\mid A_j]\).

The zonal value of an extensive quantity (a count or total) is
\(Np_j\bar z_j\) (\(N\) the total mass) and depends on both the zone
mean and the zone weight; this paper deals mainly with intensive
quantities, i.e.~zone means.

\textbf{Standard results used without proof.} The following results can
be found in graduate texts on measure theory, Fourier analysis,
functional analysis, probability and convex analysis (Folland 1999;
Rudin 1987; Rudin 1991; Durrett 2019; Rockafellar 1970):

\begin{itemize}
\tightlist
\item
  (S1) the Fubini--Tonelli theorem;
\item
  (S2) the monotone and dominated convergence theorems and Fatou's
  lemma, and the domination condition for differentiation under the
  integral sign;
\item
  (S3) the Radon--Nikodym theorem and the Hahn--Jordan decomposition of
  signed measures;
\item
  (S4) existence and almost-sure uniqueness of conditional expectations,
  together with their linearity and monotonicity;
\item
  (S5) regularity of finite Borel measures on \(\mathbb{R}^q\); density
  of \(C_c(U)\) in \(L^1(\nu)\) for every Borel measure \(\nu\) on an
  open set \(U\subseteq\mathbb{R}^q\) that is finite on compact sets
  (for example Lebesgue measure, or a finite Borel measure); density of
  \(C_c^\infty\) in \((C_0,\|\cdot\|_\infty)\); uniqueness of the
  representing measure in the Riesz representation theorem;
\item
  (S6) the Fourier transform is a bijection of the Schwartz space
  \(\mathcal{S}(\mathbb{R}^q)\) onto itself; if \(f\in L^1\) has an
  integrable Fourier transform, the inversion formula holds at every
  point of continuity of \(f\); Plancherel's theorem;
\item
  (S7) \(L^1(\lambda)^*=L^\infty(\lambda)\) for a finite measure
  \(\lambda\); the Banach--Alaoglu theorem; the Krein--Milman theorem;
\item
  (S8) Lévy's continuity theorem: if the characteristic functions of
  probability measures \(P_n\) on \(\mathbb{R}^q\) converge pointwise to
  a function \(\phi\) continuous at the origin, then \(P_n\) converges
  weakly to a probability measure \(P\) whose characteristic function is
  \(\phi\);
\item
  (S9) a convex function that is finite everywhere on \(\mathbb{R}^d\)
  or on an open interval is Lipschitz on every compact subset of its
  domain and has a subgradient at every point.
\end{itemize}

\textbf{Conventions.} The Jordan decomposition of a signed measure
\(\nu\) is \(\nu=\nu^+-\nu^-\), and its total variation is \(|\nu|\);
the total variation of a complex measure is also written \(|\nu|\). The
Fourier transform is \(\hat K(\omega)=\int e^{-is\cdot\omega}K(s)\,ds\),
and the Fourier--Stieltjes transform of a finite complex measure \(\nu\)
is \(\hat\nu(h)=\int e^{ih\cdot\omega}\,\nu(d\omega)\), with the same
sign as in the spectral representation
\(C(h)=\int e^{ih\cdot\omega}F(d\omega)\) of a covariance function.

\hypertarget{conditional-expectation-and-aggregation}{%
\subsection{3.2 Conditional expectation and
aggregation}\label{conditional-expectation-and-aggregation}}

For a sub-\(\sigma\)-algebra \(\mathcal{G}\subseteq\mathcal{F}\) and
\(Z\in L^1(\mu;\mathbb{R}^d)\), the conditional expectation
\(\mathbb{E}[Z\mid\mathcal{G}]\) is an integrable
\(\mathcal{G}\)-measurable function \(Y\) with
\(\int_GY\,d\mu=\int_GZ\,d\mu\) for all \(G\in\mathcal{G}\); it exists
and is almost surely unique (S4).

\textbf{Lemma 3.1 (aggregation is conditional expectation).} Let
\(\mathcal{Z}=\{A_1,\dots,A_k\}\) be a \(k\)-zone zoning and
\(\mathcal{G}=\sigma(\mathcal{Z})\). Then:

\begin{enumerate}
\def\labelenumi{(\alph{enumi})}
\item
  \(\mathcal{G}=\{\bigcup_{j\in J}A_j:\ J\subseteq\{1,\dots,k\}\}\);
\item
  a function is \(\mathcal{G}\)-measurable if and only if it is constant
  on each zone;
\item
  \(\mathbb{E}[Z\mid\mathcal{G}]=\sum_j\bar z_j\mathbf{1}_{A_j}\);
\item
  the law of \(\mathbb{E}[Z\mid\mathcal{G}]\) is
  \(\sum_jp_j\delta_{\bar z_j}\). Consequently, every zonal statistic
  that depends on the zonal data only through the weighted distribution
  \(\sum_jp_j\delta_{\bar z_j}\) (means, variances, covariances,
  correlation coefficients, regression slopes, quantiles and so on)
  equals the same statistic applied to the law of
  \(\mathbb{E}[Z\mid\mathcal{G}]\). In particular
  \(\sum_jp_j\bar z_j=\mathbb{E}Z\), and the weighted zonal covariance
  matrix equals \(\operatorname{Cov}(\mathbb{E}[Z\mid\mathcal{G}])\).
\end{enumerate}

\emph{Proof.} (a) Call the right-hand side \(\mathcal{U}\). Since the
\(A_j\) are pairwise disjoint with union \(\Omega\), the complement of
\(\bigcup_{j\in J}A_j\) is \(\bigcup_{j\notin J}A_j\), and a countable
union of such sets is the union indexed by the union of the index sets.
So \(\mathcal{U}\) is a \(\sigma\)-algebra containing every \(A_j\), and
\(\sigma(\mathcal{Z})\subseteq\mathcal{U}\); conversely, every
\(\sigma\)-algebra containing all \(A_j\) contains their finite unions,
so \(\mathcal{U}\subseteq\sigma(\mathcal{Z})\).

\begin{enumerate}
\def\labelenumi{(\alph{enumi})}
\setcounter{enumi}{1}
\item
  Let \(Y\) be \(\mathcal{G}\)-measurable and suppose it takes values
  \(y\neq y'\) at two points of \(A_j\). Then \(\{Y=y\}\in\mathcal{G}\)
  is, by (a), a union of zones; it contains a point of \(A_j\), hence
  all of \(A_j\), contradicting the value \(y'\) at a point of \(A_j\).
  Conversely, every preimage of a function constant on each zone is a
  union of zones and belongs to \(\mathcal{G}\).
\item
  \(Y=\sum_j\bar z_j\mathbf{1}_{A_j}\) is \(\mathcal{G}\)-measurable and
  integrable. For \(G=\bigcup_{j\in J}A_j\),
  \(\int_GY\,d\mu=\sum_{j\in J}\bar z_jp_j=\sum_{j\in J}\int_{A_j}Z\,d\mu=\int_GZ\,d\mu\).
  If \(Y'\) also satisfies the definition, then by (b)
  \(Y'=\sum_jc_j\mathbf{1}_{A_j}\), and \(G=A_j\) gives
  \(c_jp_j=\int_{A_j}Z\,d\mu\), i.e.~\(c_j=\bar z_j\).
\item
  For a Borel set \(B\),
  \(\mu(\mathbb{E}[Z\mid\mathcal{G}]\in B)=\sum_{j:\bar z_j\in B}p_j\),
  so the law is \(\sum_jp_j\delta_{\bar z_j}\) (zones with equal means
  merge into one atom). The remaining claims follow:
  \(\sum_jp_j\bar z_j=\sum_j\int_{A_j}Z\,d\mu=\mathbb{E}Z\); taking
  \(G=\Omega\) in the definition gives
  \(\mathbb{E}\,\mathbb{E}[Z\mid\mathcal{G}]=\mathbb{E}Z\), and hence \[
  \begin{aligned}\sum_jp_j(\bar z_j-\mathbb{E}Z)(\bar z_j-\mathbb{E}Z)^\top&=\mathbb{E}\Big[\big(\mathbb{E}[Z\mid\mathcal{G}]-\mathbb{E}Z\big)\big(\mathbb{E}[Z\mid\mathcal{G}]-\mathbb{E}Z\big)^\top\Big]\\&=\operatorname{Cov}(\mathbb{E}[Z\mid\mathcal{G}]).\qquad\blacksquare\end{aligned}
  \]
\end{enumerate}

Lemma 3.1 says that every zonal statistic computed from the weighted
distribution of the zone means (weighted moments, correlations,
regression slopes and so on) is a functional of the distribution of
\(\mathbb{E}[Z\mid\mathcal{G}]\); statistics that weight zones equally
or count zones, such as the unweighted mean of the zone means, are not.
The variable obtained by assigning each zone mean to every unit of the
zone, drawn as a map, is the choropleth map.

\textbf{Lemma 3.2 (tower property).} If
\(\mathcal{G}'\subseteq\mathcal{G}\) are sub-\(\sigma\)-algebras, then
\(\mathbb{E}\big[\mathbb{E}[Z\mid\mathcal{G}]\mid\mathcal{G}'\big]=\mathbb{E}[Z\mid\mathcal{G}']\)
a.e. In particular, if the zoning \(\mathcal{Z}'\) is coarser than
\(\mathcal{Z}\) (every zone of \(\mathcal{Z}'\) is a union of zones of
\(\mathcal{Z}\)), aggregating first over \(\mathcal{Z}\) and then over
\(\mathcal{Z}'\) is the same as aggregating directly over
\(\mathcal{Z}'\).

\emph{Proof.} For \(B\in\mathcal{G}'\) we have \(B\in\mathcal{G}\), and
two applications of the definition give \[
\int_B\mathbb{E}\big[\mathbb{E}[Z\mid\mathcal{G}]\mid\mathcal{G}'\big]\,d\mu=\int_B\mathbb{E}[Z\mid\mathcal{G}]\,d\mu=\int_BZ\,d\mu .
\] The integrand on the left is \(\mathcal{G}'\)-measurable, and
uniqueness of conditional expectations gives the first claim. If
\(\mathcal{Z}'\) is coarser than \(\mathcal{Z}\), every zone of
\(\mathcal{Z}'\) belongs to \(\sigma(\mathcal{Z})\), so
\(\sigma(\mathcal{Z}')\subseteq\sigma(\mathcal{Z})\) by Lemma 3.1(a).
\(\blacksquare\)

\hypertarget{convex-order-and-the-conditional-jensen-inequality}{%
\subsection{3.3 Convex order and the conditional Jensen
inequality}\label{convex-order-and-the-conditional-jensen-inequality}}

An integrable \(\mathbb{R}^d\)-valued random vector \(U\) is smaller
than \(V\) in the \textbf{convex order}, written \(U\le_{cx}V\), if
\(\mathbb{E}\varphi(U)\le\mathbb{E}\varphi(V)\) for every convex
function \(\varphi\) for which both expectations are defined. The convex
order depends only on the distributions; \(U\le_{cx}V\) implies
\(\mathbb{E}U=\mathbb{E}V\) (take \(\varphi(z)=\pm z_i\)).

\textbf{Lemma 3.3 (a convex function is a countable supremum of affine
functions).} Let \(O\) be \(\mathbb{R}^d\) or \((0,\infty)\) (then
\(d=1\)) and \(\varphi:O\to\mathbb{R}\) convex. For every rational point
\(r\in O\cap\mathbb{Q}^d\) choose a subgradient \(g_r\) (S9) and put
\(\ell_r(z)=\varphi(r)+\langle g_r,z-r\rangle\). Then for every
\(z\in O\), \[
\varphi(z)=\sup_{r\in O\cap\mathbb{Q}^d}\ell_r(z).
\]

\emph{Proof.} By definition of a subgradient, \(\ell_r\le\varphi\), so
the right-hand side is at most \(\varphi(z)\). Fix \(z\in O\) and choose
\(\varepsilon>0\) with the closed ball
\(\bar B(z,2\varepsilon)\subseteq O\). By (S9), \(\varphi\) is
\(L\)-Lipschitz on \(\bar B(z,2\varepsilon)\). For
\(r\in B(z,\varepsilon)\) and \(g_r\neq0\), taking \(u=g_r/|g_r|\) and
\(0<t\le\varepsilon\) in
\(\varphi(r+tu)\ge\varphi(r)+t\langle g_r,u\rangle\) gives
\(t|g_r|\le\varphi(r+tu)-\varphi(r)\le Lt\), i.e.~\(|g_r|\le L\). Taking
rational points \(r_n\in B(z,\varepsilon)\) with \(r_n\to z\), we get
\(\ell_{r_n}(z)\ge\varphi(r_n)-L|z-r_n|\to\varphi(z)\) (\(\varphi\) is
continuous). \(\blacksquare\)

\textbf{Proposition 3.4 (conditional Jensen inequality).} Let
\(O,\varphi\) be as in Lemma 3.3, let \(Z\in L^1(\mu;\mathbb{R}^d)\)
take values in \(O\) almost surely with \(\varphi(Z)\in L^1\), and
suppose \(\mathbb{E}[Z\mid\mathcal{G}]\) takes values in \(O\) almost
surely. Then \(\varphi(\mathbb{E}[Z\mid\mathcal{G}])\) is integrable and
\[
\varphi(\mathbb{E}[Z\mid\mathcal{G}])\le\mathbb{E}[\varphi(Z)\mid\mathcal{G}]\quad\text{a.e.},\qquad\mathbb{E}\,\varphi(\mathbb{E}[Z\mid\mathcal{G}])\le\mathbb{E}\,\varphi(Z).
\] For \(O=\mathbb{R}^d\) the last hypothesis holds automatically; it
also holds for \(O=(0,\infty)\) when \(Z>0\) a.e.

\emph{Proof.} For every rational point \(r\), \(\varphi(Z)\ge\ell_r(Z)\)
pointwise; by monotonicity and linearity of conditional expectation,
\(\mathbb{E}[\varphi(Z)\mid\mathcal{G}]\ge\ell_r(\mathbb{E}[Z\mid\mathcal{G}])\)
outside a null set \(N_r\). The set \(N=\bigcup_rN_r\) is null, outside
\(N\) the inequality holds for all rational points simultaneously, and
Lemma 3.3 gives the first claim. For one rational point \(r_0\),
\(\ell_{r_0}(\mathbb{E}[Z\mid\mathcal{G}])\le\varphi(\mathbb{E}[Z\mid\mathcal{G}])\le\mathbb{E}[\varphi(Z)\mid\mathcal{G}]\)
with integrable outer terms, so the middle term is integrable; taking
expectations and using
\(\mathbb{E}\,\mathbb{E}[\varphi(Z)\mid\mathcal{G}]=\mathbb{E}\varphi(Z)\)
gives the second claim. Finally, if \(O=(0,\infty)\) and \(Z>0\) a.e.,
let \(G=\{\mathbb{E}[Z\mid\mathcal{G}]\le0\}\in\mathcal{G}\); then
\(0\le\int_GZ\,d\mu=\int_G\mathbb{E}[Z\mid\mathcal{G}]\,d\mu\le0\), and
since \(Z>0\), \(\mu(G)=0\). \(\blacksquare\)

\textbf{Lemma 3.5 (projection property).} Let \(X\in L^2(\mu)\). Then
\(\mathbb{E}[X\mid\mathcal{G}]\in L^2(\mu)\), and for every
\(\mathcal{G}\)-measurable \(W\in L^2(\mu)\), \[
\mathbb{E}\big[(X-\mathbb{E}[X\mid\mathcal{G}])\,W\big]=0 .
\]

\emph{Proof.} By Proposition 3.4 with \(\varphi(x)=x^2\),
\(\mathbb{E}[X\mid\mathcal{G}]^2\le\mathbb{E}[X^2\mid\mathcal{G}]\)
a.e., and taking expectations gives
\(\mathbb{E}\big[\mathbb{E}[X\mid\mathcal{G}]^2\big]\le\mathbb{E}X^2<\infty\).
Write \(\tilde X=\mathbb{E}[X\mid\mathcal{G}]\). For
\(G\in\mathcal{G}\),
\(\mathbb{E}[X\mathbf{1}_G]=\mathbb{E}[\tilde X\mathbf{1}_G]\) by
definition; by linearity, \(\mathbb{E}[Xs]=\mathbb{E}[\tilde Xs]\) for
\(\mathcal{G}\)-measurable simple functions \(s\). For a bounded
\(\mathcal{G}\)-measurable \(W_0\), take \(\mathcal{G}\)-measurable
simple functions \(s_n\) converging uniformly to \(W_0\); from
\(|X(s_n-W_0)|\le\|s_n-W_0\|_\infty|X|\) and the same estimate for
\(\tilde X\) we get \(\mathbb{E}[XW_0]=\mathbb{E}[\tilde XW_0]\). For a
general \(W\in L^2\), put \(W_n=W\mathbf{1}\{|W|\le n\}\), which is
bounded and \(\mathcal{G}\)-measurable; \(|XW_n|\le|XW|\) and
\(|\tilde XW_n|\le|\tilde XW|\) with integrable right-hand sides by
Cauchy--Schwarz, and dominated convergence gives
\(\mathbb{E}[XW]=\mathbb{E}[\tilde XW]\). \(\blacksquare\)

\textbf{Corollary 3.6 (aggregation is a contraction in the convex
order).} For every convex function \(\varphi\) on \(\mathbb{R}^d\),
\(\mathbb{E}\varphi(Z)\in(-\infty,+\infty]\) is defined and
\(\mathbb{E}\,\varphi(\mathbb{E}[Z\mid\mathcal{G}])\le\mathbb{E}\,\varphi(Z)\);
that is, \(\mathbb{E}[Z\mid\mathcal{G}]\le_{cx}Z\). If \(Z\in L^2\),
then
\(\operatorname{Cov}(\mathbb{E}[Z\mid\mathcal{G}])\preceq\operatorname{Cov}(Z)\)
in the Loewner order.

\emph{Proof.} By Lemma 3.3, \(\varphi\ge\ell_{r_0}\), so the negative
parts of \(\varphi(Z)\) and \(\varphi(\mathbb{E}[Z\mid\mathcal{G}])\)
are integrable and the expectations are defined. If
\(\mathbb{E}\varphi(Z)=+\infty\) the inequality is trivial; otherwise it
follows from Proposition 3.4. For \(a\in\mathbb{R}^d\), taking
\(\varphi(z)=(a^\top(z-\mathbb{E}Z))^2\) and using
\(\mathbb{E}\,\mathbb{E}[Z\mid\mathcal{G}]=\mathbb{E}Z\) gives
\(a^\top\operatorname{Cov}(\mathbb{E}[Z\mid\mathcal{G}])\,a\le a^\top\operatorname{Cov}(Z)\,a\).
\(\blacksquare\)

\textbf{Proposition 3.7 (monotonicity of dispersion functionals).} Let
\(Z>0\) be an integrable real random variable, \(m=\mathbb{E}Z\) and
\(\tilde Z=\mathbb{E}[Z\mid\mathcal{G}]\).

\begin{enumerate}
\def\labelenumi{(\roman{enumi})}
\item
  If \(Z\in L^2\), then
  \(\operatorname{Var}\tilde Z\le\operatorname{Var}Z\).
\item
  If the Theil index
  \(\operatorname{Th}(Z)=\mathbb{E}\big[(Z/m)\log(Z/m)\big]\) is finite,
  then \(\operatorname{Th}(\tilde Z)\le \operatorname{Th}(Z)\).
\item
  With the Gini coefficient
  \(\operatorname{Gini}(U)=\mathbb{E}|U-U'|/(2\,\mathbb{E}U)\) (\(U'\)
  an independent copy of \(U\)),
  \(\operatorname{Gini}(\tilde Z)\le \operatorname{Gini}(Z)\).
\end{enumerate}

\emph{Proof.} (i) is the case \(d=1\) of Corollary 3.6. (ii)
\(\varphi(z)=(z/m)\log(z/m)\) is convex on \((0,\infty)\)
(\(\varphi''(z)=1/(mz)>0\)); by Proposition 3.4 with \(O=(0,\infty)\)
and \(\mathbb{E}\tilde Z=m\),
\(\operatorname{Th}(\tilde Z)=\mathbb{E}\varphi(\tilde Z)\le\mathbb{E}\varphi(Z)=\operatorname{Th}(Z)\).
(iii) We first show: if real integrable variables satisfy \(U\le_{cx}V\)
and \(U',V'\) are independent copies, then
\(\mathbb{E}|U-U'|\le\mathbb{E}|V-V'|\). For fixed \(u\),
\(v\mapsto|u-v|\) is convex, so \(\mathbb{E}|u-U'|\le\mathbb{E}|u-V'|\).
By independence and Fubini's theorem, \[
\mathbb{E}|U-U'|=\int\mathbb{E}|u-U'|\,P_U(du)\le\int\psi(u)\,P_U(du),\qquad \psi(u):=\mathbb{E}|u-V'| .
\] \(\psi\) is an average of convex functions, hence convex, and
\(\psi(u)\le|u|+\mathbb{E}|V|\) is finite; applying \(U\le_{cx}V\) once
more gives \(\int\psi\,dP_U\le\int\psi\,dP_V=\mathbb{E}|V-V'|\). Take
\(U=\tilde Z\) and \(V=Z\) (Corollary 3.6); both have mean \(m\), which
gives (iii). \(\blacksquare\)

The indices in Proposition 3.7 are computed under \(\mu\) after
assigning each zone mean to all units of its zone. Applied to
\(\mathbb{E}[Z\mid\mathcal{G}]\) in place of \(Z\) and a coarser
\(\mathcal{G}'\subseteq\mathcal{G}\) (Lemma 3.2), the proposition says
that as one aggregates up a nested hierarchy the variance, the Theil
index and the Gini coefficient, like every dispersion functional that is
monotone in the convex order, can only decrease or stay the same.

\hypertarget{notation}{%
\subsection{3.4 Notation}\label{notation}}

\begin{longtable}[]{@{}
  >{\raggedright\arraybackslash}p{(\columnwidth - 4\tabcolsep) * \real{0.3077}}
  >{\raggedright\arraybackslash}p{(\columnwidth - 4\tabcolsep) * \real{0.5769}}
  >{\raggedright\arraybackslash}p{(\columnwidth - 4\tabcolsep) * \real{0.1154}}@{}}
\caption{Main notation}\tabularnewline
\toprule\noalign{}
\begin{minipage}[b]{\linewidth}\raggedright
Symbol
\end{minipage} & \begin{minipage}[b]{\linewidth}\raggedright
Meaning
\end{minipage} & \begin{minipage}[b]{\linewidth}\raggedright
First used
\end{minipage} \\
\midrule\noalign{}
\endfirsthead
\toprule\noalign{}
\begin{minipage}[b]{\linewidth}\raggedright
Symbol
\end{minipage} & \begin{minipage}[b]{\linewidth}\raggedright
Meaning
\end{minipage} & \begin{minipage}[b]{\linewidth}\raggedright
First used
\end{minipage} \\
\midrule\noalign{}
\endhead
\bottomrule\noalign{}
\endlastfoot
\((\Omega,\mathcal{F},\mu)\) & fine-scale units and normalized support
measure & 3.1 \\
\(\mathcal{Z}=\{A_1,\dots,A_k\}\), \(p_j=\mu(A_j)\) & zoning and zone
weights & 3.1 \\
\(\mathcal{G}=\sigma(\mathcal{Z})\) & sub-\(\sigma\)-algebra generated
by the zoning & 3.1 \\
\(\bar z_j\), \(\mathbb{E}[Z\mid\mathcal{G}]\) & zone means; the
variable assigning zone means to the units of each zone & 3.1, 3.2 \\
\(\mathbb{E}[Z\mid A]\), \(\operatorname{Var}(X\mid A)\) & mean and
variance over a set \(A\) with \(\mu(A)>0\) & 3.1 \\
\(\rho\), \(\rho_B\), \(\rho_W\) & fine-scale, zonal and pooled
within-zone correlations & 5 \\
\(\eta_X^2\), \(a=1-\eta_X^2\) & correlation ratio (between-zone
variance share), within-zone variance share & 5, 7 \\
\(\ell=\sigma/L\), \(\tau=\delta/(2\ell)\) & variation length; half the
ratio of zone diameter to variation length & 7 \\
\(F\), \(\Gamma=\operatorname{Re}F_{XY}\), \(\operatorname{coh}\) &
spectral matrix measure, real part of the cross-spectrum, real coherence
& 8 \\
\(K\), \(\hat K\), \(\pi\), \(g_\pi\) & aggregation kernel (zone shape)
and its Fourier transform, zoning design, filter of a design & 8 \\
\(v_l\), \(c_l\), \(\kappa_l\) & variance component, covariance
component and coherence at level \(l\) & 9 \\
\(L\), \(U\) & lower and upper bounds on the fine-scale covariance from
zonal data alone & 10 \\
\end{longtable}

\hypertarget{atomless-support-measures-and-the-convexity-of-zonings}{%
\section{4 Atomless support measures and the convexity of
zonings}\label{atomless-support-measures-and-the-convexity-of-zonings}}

This section builds the tools needed in Sections 5 and 6. The central
fact is that, under an atomless support measure, every ``proportional
allocation'' given by a fuzzy partition can be realized exactly by a
genuine zoning; consequently the set of attainable zone weights and zone
totals is convex, and so is the set of attainable zone-mean vectors for
prescribed zone weights (Corollary 4.10). This turns the combinatorial
zoning problem into a convex one.

\hypertarget{atoms-and-the-sierpiux144ski-theorem}{%
\subsection{4.1 Atoms and the Sierpiński
theorem}\label{atoms-and-the-sierpiux144ski-theorem}}

\textbf{Definition 4.1 (atom).} Let \(\lambda\) be a finite measure on a
measurable space. A measurable set \(A\) is an \textbf{atom} of
\(\lambda\) if \(\lambda(A)>0\) and every measurable subset \(B\) of
\(A\) has \(\lambda(B)\in\{0,\lambda(A)\}\). A measure without atoms is
\textbf{atomless}. A signed measure is called atomless if its total
variation is atomless.

For a support measure on finitely many units (\(\Omega=\{1,\dots,n\}\)),
every unit is an atom; a measure with a density with respect to Lebesgue
measure (for instance area on a region of \(\mathbb{R}^q\)) is atomless:
if \(\mu(A)>0\), then \(t\mapsto\mu(A\cap\{s:s_1<t\})\) is continuous
(hyperplanes are Lebesgue-null) and increases from \(0\) to \(\mu(A)\),
so it takes the value \(\mu(A)/2\). Atomlessness can be seen as the
idealization of basic units that are small relative to the zones. The
atomless assumption is used only in Corollaries 4.9--4.10, Theorem 5.5,
Section 6, Theorem 10.1(ii) with the sharpness statements that follow
from it (the last sentence of Theorem 10.1 and the identification
statement of Corollary 10.2), and the last sentence of Corollary 11.4;
Theorems 4.6 and 4.8 are themselves statements about atomless vector
measures. All other results of Sections 3, 5, 7, 9, 10 and 11 hold for
every support measure, and Section 8 deals with random fields on
continuous space.

\textbf{Lemma 4.2.} Let \(\lambda\) be a finite atomless measure and
\(\lambda(A)>0\). Then for every \(\varepsilon>0\) there is a measurable
\(B\subseteq A\) with \(0<\lambda(B)<\varepsilon\).

\emph{Proof.} \(A\) is not an atom, so there is \(A_1\subseteq A\) with
\(0<\lambda(A_1)<\lambda(A)\). One of \(A_1\) and \(A\setminus A_1\) has
measure in \((0,\lambda(A)/2]\); call it \(A'\). Repeating the step
\(n\) times yields a subset with measure in \((0,\lambda(A)/2^n]\);
choose \(n\) with \(\lambda(A)/2^n<\varepsilon\). \(\blacksquare\)

\textbf{Theorem 4.3 (Sierpiński 1922).} Let \(\lambda\) be a finite
atomless measure, \(A\) measurable and \(t\in[0,\lambda(A)]\). Then
there is a measurable \(B\subseteq A\) with \(\lambda(B)=t\).

\emph{Proof.} The cases \(t=0\) and \(t=\lambda(A)\) are trivial; let
\(0<t<\lambda(A)\). Put \(B_0=\varnothing\). Given \(B_n\subseteq A\)
with \(\lambda(B_n)\le t\), let \[
s_n=\sup\{\lambda(C):\ C\subseteq A\setminus B_n\ \text{measurable},\ \lambda(B_n)+\lambda(C)\le t\}\ (\ge0),
\] choose \(C_n\subseteq A\setminus B_n\) with
\(\lambda(B_n)+\lambda(C_n)\le t\) and \(\lambda(C_n)\ge s_n/2\)
(\(C_n=\varnothing\) if \(s_n=0\)), and put \(B_{n+1}=B_n\cup C_n\).
Then \(\lambda(B_n)\le t\) increases, and \(B=\bigcup_nB_n\) satisfies
\(\lambda(B)=\lim_n\lambda(B_n)\le t\). The \(C_n\) are pairwise
disjoint with \(\sum_n\lambda(C_n)=\lambda(B)<\infty\), so
\(\lambda(C_n)\to0\) and hence \(s_n\le2\lambda(C_n)\to0\). If
\(\lambda(B)<t\), then \(\lambda(A\setminus B)>\lambda(A)-t>0\), and
Lemma 4.2 gives \(C\subseteq A\setminus B\) with
\(0<\lambda(C)<t-\lambda(B)\). For every \(n\),
\(C\subseteq A\setminus B_n\) and
\(\lambda(B_n)+\lambda(C)\le\lambda(B)+\lambda(C)<t\), so
\(\lambda(C)\le s_n\); letting \(n\to\infty\) gives \(\lambda(C)=0\), a
contradiction. Hence \(\lambda(B)=t\). \(\blacksquare\)

\textbf{Corollary 4.4 (equipartition).} Under the assumptions of Theorem
4.3, if \(\lambda(E)>0\) then for every positive integer \(n\), \(E\)
can be partitioned into \(n\) measurable sets each of measure
\(\lambda(E)/n\).

\emph{Proof.} By Theorem 4.3 choose \(E_1\subseteq E\) with
\(\lambda(E_1)=\lambda(E)/n\), then \(E_2\subseteq E\setminus E_1\) with
\(\lambda(E_2)=\lambda(E)/n\), and so on; the last set is the remainder.
\(\blacksquare\)

\textbf{Lemma 4.5.}

\begin{enumerate}
\def\labelenumi{(\alph{enumi})}
\item
  If \(\nu\) is a finite atomless measure and \(\lambda\ll\nu\) is a
  finite measure, then \(\lambda\) is atomless.
\item
  The restriction of an atomless measure to a measurable subset is
  atomless.
\item
  A finite sum of finite atomless measures is atomless.
\item
  If \(\nu\) is a finite signed measure and \(h\in L^1(|\nu|)\) (for
  instance \(h\) bounded measurable), then \(h\nu\) is a finite signed
  measure with \(|h\nu|=|h|\,|\nu|\); hence if \(|\nu|\) is atomless, so
  is \(h\nu\).
\end{enumerate}

\emph{Proof.} (a) Let \(A\) be an atom of \(\lambda\). Since
\(\lambda\ll\nu\), \(\nu(A)>0\). By Theorem 4.3 choose \(B\subseteq A\)
with \(\nu(B)=\nu(A)/2\). \(\lambda(B)\) is \(0\) or \(\lambda(A)\): in
the first case put \(A_1=A\setminus B\), in the second \(A_1=B\); in
both cases \(\nu(A_1)=\nu(A)/2\), \(\lambda(A_1)=\lambda(A)\), and
\(A_1\) is again an atom of \(\lambda\). This produces a decreasing
sequence \(A_n\) with \(\nu(A_n)=2^{-n}\nu(A)\) and
\(\lambda(A_n)=\lambda(A)\). By continuity from above,
\(\nu(\bigcap_nA_n)=0\) while \(\lambda(\bigcap_nA_n)=\lambda(A)>0\),
contradicting \(\lambda\ll\nu\).

\begin{enumerate}
\def\labelenumi{(\alph{enumi})}
\setcounter{enumi}{1}
\item
  Immediate from the definition.
\item
  Suppose \(\lambda=\sum_{i\le m}\lambda_i\) has an atom \(A\); pick
  \(i\) with \(\lambda_i(A)>0\) and, by Theorem 4.3, \(B\subseteq A\)
  with \(\lambda_i(B)=\lambda_i(A)/2\). If \(\lambda(B)=0\) then
  \(\lambda_i(B)=0\); if \(\lambda(B)=\lambda(A)\) then
  \(\lambda(A\setminus B)=0\), so \(\lambda_i(B)=\lambda_i(A)\). Both
  are contradictions.
\item
  Take a Hahn decomposition \(\Omega=P\cup N\) of \(\nu\) and put
  \(P_h=(P\cap\{h\ge0\})\cup(N\cap\{h<0\})\),
  \(N_h=\Omega\setminus P_h\). For every measurable \(E\),
  \((h\nu)(E\cap P_h)=\int_{E\cap P_h}|h|\,d|\nu|\ge0\) and
  \((h\nu)(E\cap N_h)=-\int_{E\cap N_h}|h|\,d|\nu|\le0\), so
  \((P_h,N_h)\) is a Hahn decomposition of \(h\nu\) and
  \(|h\nu|(E)=(h\nu)(E\cap P_h)-(h\nu)(E\cap N_h)=\int_E|h|\,d|\nu|\).
  Hence \(|h\nu|\ll|\nu|\), and (a) gives the claim. \(\blacksquare\)
\end{enumerate}

\begin{figure}
\centering
\includegraphics[width=1\textwidth,height=\textheight]{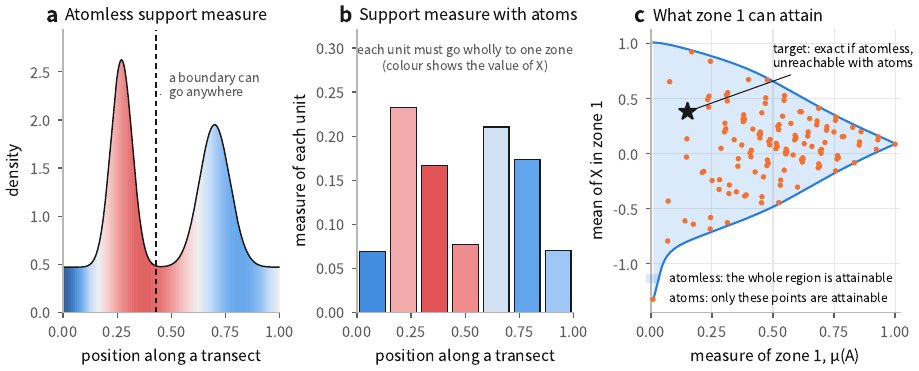}
\caption{Atoms and atomless measures. (a) An atomless support measure
distributed continuously along a transect; zone boundaries can be placed
anywhere. Color shows the value of the variable \(X\). (b) The same
measure lumped into 7 indivisible units (atoms); each unit can only be
assigned to a zone as a whole. (c) The combinations of (measure, zone
mean) attainable by zone 1. Without atoms it is the whole shaded region,
whose upper and lower boundaries come from the two extreme ways of
filling the zone: taking the part with the largest values of \(X\)
first, or the smallest first. With atoms only 127 points are attainable,
namely all non-empty combinations of the 7 units. The target marked by
the star is attained exactly without atoms but cannot be attained with
atoms. This is precisely the difference captured by Lyapunov's convexity
theorem (Theorem 4.6).}
\end{figure}

\hypertarget{lyapunovs-convexity-theorem}{%
\subsection{4.2 Lyapunov's convexity
theorem}\label{lyapunovs-convexity-theorem}}

\textbf{Theorem 4.6 (strong form of Lyapunov's convexity theorem;
Lyapunov 1940, proof after Lindenstrauss 1966).} Let
\(\nu=(\nu_1,\dots,\nu_m)\) be a finite signed vector measure on
\((\Omega,\mathcal{F})\) with every \(|\nu_i|\) atomless. Then for every
measurable \(g_0:\Omega\to[0,1]\) there is \(A\in\mathcal{F}\) with \[
\nu(A)=\int g_0\,d\nu .
\] In particular, the range \(\{\nu(A):A\in\mathcal{F}\}\) of \(\nu\) is
a convex subset of \(\mathbb{R}^m\).

\emph{Proof.}

\begin{enumerate}
\def\labelenumi{\arabic{enumi}.}
\tightlist
\item
  Let \(\lambda=\sum_i|\nu_i|\). By Lemma 4.5(c), \(\lambda\) is a
  finite atomless measure. Since \(|\nu_i|\le\lambda\), the
  Radon--Nikodym theorem gives \(\nu_i=\rho_i\lambda\) with
  \(|\rho_i|\le1\); for bounded measurable \(g\),
  \(\int g\,d\nu_i=\int g\rho_i\,d\lambda\).
\item
  Equip \(L^\infty(\lambda)=L^1(\lambda)^*\) with the weak-* topology
  and define \(T(g)=\big(\int g\rho_i\,d\lambda\big)_{i\le m}\). Since
  \(\rho_i\in L^1(\lambda)\), \(T\) is a weak-* continuous linear map.
  Let \[
  \mathcal{K}=\{g\in L^\infty(\lambda):\ 0\le g\le1\ \text{a.e.},\ T(g)=T(g_0)\}.
  \] \(\mathcal{K}\) contains \(g_0\) and is convex. The condition
  ``\(0\le g\le1\) a.e.'' is equivalent to
  \(0\le\int gu\,d\lambda\le\int u\,d\lambda\) for every
  \(0\le u\in L^1(\lambda)\), an intersection of weak-* closed
  conditions; \(T(g)=T(g_0)\) is weak-* closed as well. \(\mathcal{K}\)
  lies in the unit ball, so it is weak-* compact by the Banach--Alaoglu
  theorem. The weak-* topology is Hausdorff and locally convex, so by
  the Krein--Milman theorem \(\mathcal{K}\) has an extreme point
  \(g^*\).
\item
  \textbf{Extreme points are indicators almost everywhere.} Otherwise
  \(\lambda(0<g^*<1)>0\). Since
  \(\{0<g^*<1\}=\bigcup_n\{1/n\le g^*\le1-1/n\}\), there are
  \(\varepsilon>0\) and a set \(E\) with \(\lambda(E)>0\) on which
  \(\varepsilon\le g^*\le1-\varepsilon\). By Corollary 4.4 (\(\lambda\)
  is atomless on \(E\)), partition \(E\) into \(m+1\) sets
  \(E_0,\dots,E_m\) of positive measure. The linear map
  \(c\mapsto T\big(\sum_{i=0}^mc_i\mathbf{1}_{E_i}\big)\) from
  \(\mathbb{R}^{m+1}\) to \(\mathbb{R}^m\) has a non-zero kernel vector
  \(c\); put \(h=\sum_ic_i\mathbf{1}_{E_i}\) and rescale so that
  \(\|h\|_\infty=\max_i|c_i|\le\varepsilon\). Then \(h\neq0\) (each
  \(E_i\) has positive measure), \(T(h)=0\), and \(h\) vanishes outside
  \(E\), so \(0\le g^*\pm h\le1\) and \(T(g^*\pm h)=T(g_0)\),
  i.e.~\(g^*\pm h\in\mathcal{K}\). But
  \(g^*=\tfrac12(g^*+h)+\tfrac12(g^*-h)\) with \(g^*+h\neq g^*-h\),
  contradicting extremality.
\item
  Hence there is a measurable \(A\) with \(g^*=\mathbf{1}_A\) a.e., and
  \(\nu_i(A)=\int\mathbf{1}_A\rho_i\,d\lambda=\int g_0\rho_i\,d\lambda=\int g_0\,d\nu_i\)
  for all \(i\).
\item
  Convexity of the range: for \(A,B\in\mathcal{F}\) and \(t\in[0,1]\),
  take \(g_0=t\mathbf{1}_A+(1-t)\mathbf{1}_B\in[0,1]\); there is \(C\)
  with \(\nu(C)=t\nu(A)+(1-t)\nu(B)\). \(\blacksquare\)
\end{enumerate}

\hypertarget{fuzzy-partitions-and-the-dvoretzkywaldwolfowitz-theorem}{%
\subsection{4.3 Fuzzy partitions and the Dvoretzky--Wald--Wolfowitz
theorem}\label{fuzzy-partitions-and-the-dvoretzkywaldwolfowitz-theorem}}

\textbf{Definition 4.7 (fuzzy partition).} A \(k\)-part fuzzy partition
is a family of measurable functions \(f=(f_1,\dots,f_k)\) with
\(f_j\ge0\) and \(\sum_jf_j=1\). The indicator family
\((\mathbf{1}_{A_1},\dots,\mathbf{1}_{A_k})\) is the fuzzy partition of
the partition \(\{A_j\}\).

Fuzzy partitions are much more flexible than partitions: the zone masses
\(\int f_j\,d\mu\) and zone totals \(\int f_jZ\,d\mu\) are linear in
\(f_j\) and can be tuned continuously. The following theorem says that,
for atomless measures, this flexibility brings no additional attainable
values.

\textbf{Theorem 4.8 (partition version; Dvoretzky, Wald \& Wolfowitz
1951).} Let \(\nu=(\nu_1,\dots,\nu_m)\) satisfy the assumptions of
Theorem 4.6. Then for every \(k\)-part fuzzy partition \(f\) there is a
measurable partition \(\{A_1,\dots,A_k\}\) with
\(\nu(A_j)=\int f_j\,d\nu\) for all \(j\).

\emph{Proof.} Induction on \(k\). For \(k=1\), \(f_1\equiv1\) and
\(A_1=\Omega\). Assume the claim for \(k-1\) (for every measurable space
and every vector measure satisfying the assumptions of Theorem 4.6) and
consider \(k\).

\begin{enumerate}
\def\labelenumi{\arabic{enumi}.}
\tightlist
\item
  Let \(g=f_1+\dots+f_{k-1}=1-f_k\in[0,1]\) and \(\Omega'=\{g>0\}\). For
  \(j<k\), \(0\le f_j\le g\), so \(f_j=0\) on \(\{g=0\}\); on
  \(\Omega'\), \(0\le f_j/g\le1\) and \(\sum_{j<k}f_j/g=1\).
\item
  On the measurable subsets of \(\Omega'\) define the \(mk\)-dimensional
  vector measure \(\lambda\) with components \(\nu_i|_{\Omega'}\) and
  \((f_j/g)\nu_i|_{\Omega'}\) (\(i\le m\), \(j<k\)). By Lemma 4.5(b)(d),
  all components are finite and atomless.
\item
  Apply Theorem 4.6 on \(\Omega'\) to \(\lambda\) and the function
  \(g|_{\Omega'}\in[0,1]\): there is a measurable \(B\subseteq\Omega'\)
  with \(\lambda(B)=\int_{\Omega'}g\,d\lambda\). Componentwise, \[
  \nu_i(B)=\int_{\Omega'}g\,d\nu_i=\int g\,d\nu_i,\qquad
  \int_B\frac{f_j}{g}\,d\nu_i=\int_{\Omega'}g\cdot\frac{f_j}{g}\,d\nu_i=\int f_j\,d\nu_i\quad(j<k),
  \] where removing \(\Omega'\) uses \(g=f_j=0\) on \(\{g=0\}\).
\item
  On \(B\), \((f_1/g,\dots,f_{k-1}/g)\) is a \((k-1)\)-part fuzzy
  partition, and the restriction of \(\nu\) to \(B\) satisfies the
  assumptions of Theorem 4.6 (Lemma 4.5(b)). By the induction hypothesis
  there is a measurable partition \(A_1,\dots,A_{k-1}\) of \(B\) with
  \(\nu(A_j)=\int_B(f_j/g)\,d\nu=\int f_j\,d\nu\) (\(j<k\)).
\item
  Put \(A_k=\Omega\setminus B\). Then
  \(\nu(A_k)=\nu(\Omega)-\nu(B)=\int(1-g)\,d\nu=\int f_k\,d\nu\), and
  \(\{A_1,\dots,A_k\}\) is a measurable partition of \(\Omega\).
\end{enumerate}

If \(\Omega'=\varnothing\) (i.e.~\(f_k\equiv1\)), take
\(B=\varnothing\); the steps remain valid. \(\blacksquare\)

The theorem does not require the \(A_j\) to have positive measure; in
the applications below positivity comes from
\(\mu(A_j)=\int f_j\,d\mu\).

\textbf{Corollary 4.9 (zonings with equal means).} Let \(\mu\) be
atomless and \(Z\in L^1(\mu;\mathbb{R}^d)\). For every \(k\) and every
\(p\in\Delta_k^\circ\) there is a \(k\)-zone zoning with zone weights
\(p\) in which every zone mean equals \(\mathbb{E}Z\).

\emph{Proof.} The components of the vector measure
\(\nu=(\mu,Z_1\mu,\dots,Z_d\mu)\) are finite signed measures; by Lemma
4.5(d), \(|Z_i\mu|=|Z_i|\,\mu\ll\mu\), which is atomless by Lemma
4.5(a). The constant functions \(f_j\equiv p_j\) form a fuzzy partition.
By Theorem 4.8 there is a partition with \(\mu(A_j)=p_j\) and
\(\int_{A_j}Z_i\,d\mu=p_j\mathbb{E}Z_i\); since \(p_j>0\),
\(\mathbb{E}[Z\mid A_j]=\mathbb{E}Z\). \(\blacksquare\)

\textbf{Corollary 4.10 (convexity of the set of zonal configurations).}
Let \(\mu\) be atomless, \(Z\in L^1(\mu;\mathbb{R}^d)\) and
\(p\in\Delta_k^\circ\). Let
\(\mathcal{C}_k(p)\subseteq(\mathbb{R}^d)^k\) be the set of vectors of
zone means \((\bar z_1,\dots,\bar z_k)\) of all \(k\)-zone zonings with
zone weights \(p\). Then \(\mathcal{C}_k(p)\) is convex, and it equals
the set of ``fuzzy zone means'' \(\big(p_j^{-1}\int f_jZ\,d\mu\big)_j\)
of all fuzzy partitions with \(\int f_j\,d\mu=p_j\).

\emph{Proof.} Every partition is a fuzzy partition, so
\(\mathcal{C}_k(p)\) is contained in the latter set. Conversely, given a
fuzzy partition with \(\int f_j\,d\mu=p_j\), apply Theorem 4.8 to
\(\nu=(\mu,Z_1\mu,\dots,Z_d\mu)\) (atomless, as in the proof of
Corollary 4.9): there is a partition with \(\mu(A_j)=p_j\) and
\(\int_{A_j}Z\,d\mu=\int f_jZ\,d\mu\), so its zone means equal the fuzzy
zone means. The latter set is convex: if \(f,f'\) satisfy the
constraints, so does the fuzzy partition \(tf+(1-t)f'\), and its fuzzy
zone means are the convex combination of the two. \(\blacksquare\)

Corollary 4.10 is the geometric core of the indeterminacy theorems of
Section 6: to show that a zonal configuration is attainable, it suffices
to construct a fuzzy partition.

\hypertarget{variance-decomposition-simpson-reversal-and-an-impossibility-theorem}{%
\section{5 Variance decomposition, Simpson reversal and an impossibility
theorem}\label{variance-decomposition-simpson-reversal-and-an-impossibility-theorem}}

The fine-scale correlation is a weighted combination of the zonal
correlation and the pooled within-zone correlation, with weights
determined solely by the between-zone variance shares (Theorem 5.2). The
necessary and sufficient condition for Simpson reversal follows at once
(Corollary 5.4). Combined with Section 4, a stronger negative result
follows: under an atomless support measure, a zonal statistic that takes
the same value for all zonings can only be a function of the overall
mean (Theorem 5.5).

\hypertarget{the-withinbetween-decomposition}{%
\subsection{5.1 The within/between
decomposition}\label{the-withinbetween-decomposition}}

\textbf{Notation.} Let \(X,Y\in L^2(\mu)\) with
\(\sigma_X^2,\sigma_Y^2>0\), and let \(\mathcal{G}\) be any
sub-\(\sigma\)-algebra. Write \(\tilde X=\mathbb{E}[X\mid\mathcal{G}]\)
and \(\tilde Y=\mathbb{E}[Y\mid\mathcal{G}]\) (both in \(L^2\) by Lemma
3.5), and put

\begin{itemize}
\tightlist
\item
  between-zone variance \(B_X=\operatorname{Var}(\tilde X)\),
  within-zone variance \(W_X=\mathbb{E}\big[(X-\tilde X)^2\big]\),
  correlation ratio \(\eta_X^2=B_X/\sigma_X^2\); likewise for \(Y\);
\item
  between-zone covariance \(C_B=\operatorname{Cov}(\tilde X,\tilde Y)\),
  within-zone covariance
  \(C_W=\mathbb{E}\big[(X-\tilde X)(Y-\tilde Y)\big]\);
\item
  the zonal correlation (when \(B_XB_Y>0\)), the pooled within-zone
  correlation (when \(W_XW_Y>0\)) and the fine-scale correlation \[
  \rho_B=\frac{C_B}{\sqrt{B_XB_Y}},\qquad \rho_W=\frac{C_W}{\sqrt{W_XW_Y}},\qquad \rho=\frac{\operatorname{Cov}(X,Y)}{\sigma_X\sigma_Y}.
  \]
\end{itemize}

For a \(k\)-zone zoning, Lemma 3.1 gives
\(B_X=\sum_jp_j(\bar x_j-\mathbb{E}X)^2\) and
\(W_X=\sum_jp_j\operatorname{Var}(X\mid A_j)\), and \(\rho_B\) is the
correlation of the zone means weighted by the zone weights.

\textbf{Lemma 5.1 (law of total variance and covariance).}
\(\operatorname{Cov}(X,Y)=C_B+C_W\); in particular
\(\sigma_X^2=B_X+W_X\) and \(\eta_X^2\in[0,1]\).

\emph{Proof.} Taking \(G=\Omega\) in the definition gives
\(\mathbb{E}\tilde X=\mathbb{E}X\) and
\(\mathbb{E}\tilde Y=\mathbb{E}Y\). Write \[
\begin{aligned}(X-\mathbb{E}X)(Y-\mathbb{E}Y)={}&(X-\tilde X)(Y-\tilde Y)+(X-\tilde X)(\tilde Y-\mathbb{E}Y)\\&+(\tilde X-\mathbb{E}X)(Y-\tilde Y)+(\tilde X-\mathbb{E}X)(\tilde Y-\mathbb{E}Y);\end{aligned}
\] all four terms are products of \(L^2\) functions, hence integrable.
\(\tilde Y-\mathbb{E}Y\) is a \(\mathcal{G}\)-measurable \(L^2\)
function, so by Lemma 3.5 the second term has expectation zero; the
third likewise. The first term has expectation \(C_W\), the fourth
\(\mathbb{E}[(\tilde X-\mathbb{E}\tilde X)(\tilde Y-\mathbb{E}\tilde Y)]=C_B\).
Taking expectations gives \(\operatorname{Cov}(X,Y)=C_B+C_W\); with
\(Y=X\) we get \(\sigma_X^2=B_X+W_X\) with both terms non-negative, so
\(\eta_X^2\in[0,1]\). \(\blacksquare\)

\textbf{Theorem 5.2 (fundamental identity; an equivalent form appears in
Robinson 1950).} \[
\rho=\eta_X\,\eta_Y\,\rho_B+\sqrt{(1-\eta_X^2)(1-\eta_Y^2)}\;\rho_W,\qquad |\rho_B|\le1,\quad|\rho_W|\le1 .
\] Convention: when \(B_XB_Y=0\), \(C_B=0\) and the first term is \(0\);
when \(W_XW_Y=0\), \(C_W=0\) and the second term is \(0\).

\emph{Proof.} By Lemma 5.1, \(\rho\,\sigma_X\sigma_Y=C_B+C_W\), and
\(B_X=\eta_X^2\sigma_X^2\), \(W_X=(1-\eta_X^2)\sigma_X^2\), likewise for
\(Y\). The first identity below holds when \(B_XB_Y>0\), the second when
\(W_XW_Y>0\): \[
\begin{aligned}C_B&=\rho_B\sqrt{B_XB_Y}=\rho_B\,\eta_X\eta_Y\,\sigma_X\sigma_Y,\\C_W&=\rho_W\sqrt{W_XW_Y}=\rho_W\sqrt{(1-\eta_X^2)(1-\eta_Y^2)}\,\sigma_X\sigma_Y .\end{aligned}
\] If \(B_X=0\), then \(\operatorname{Var}\tilde X=0\),
\(\tilde X=\mathbb{E}X\) a.e., and \(C_B=0\); likewise if \(B_Y=0\). If
\(W_XW_Y=0\), Cauchy--Schwarz gives \(|C_W|\le\sqrt{W_XW_Y}=0\).
Substituting and dividing by \(\sigma_X\sigma_Y>0\) gives the identity.
\(|\rho_B|\le1\) is Cauchy--Schwarz for \(\tilde X-\mathbb{E}X\) and
\(\tilde Y-\mathbb{E}Y\), and \(|\rho_W|\le1\) is Cauchy--Schwarz for
\(X-\tilde X\) and \(Y-\tilde Y\). \(\blacksquare\)

The zonal correlation is not an amplified or attenuated version of the
fine-scale correlation but a different quantity: the fine-scale
correlation is a weighted sum of it and the pooled within-zone
correlation, and the contribution of the zonal correlation is multiplied
by \(\eta_X\eta_Y\). For example, if only 10\% of the variance of each
variable lies between zones (\(\eta_X^2=\eta_Y^2=0.1\)), \(\rho_W=0.1\)
and the zone means are nearly collinear (\(\rho_B=0.9\)), then
\(\rho=0.1\times0.9+0.9\times0.1=0.18\): the zonal correlation is high
while the fine-scale correlation is weak.

\textbf{Corollary 5.3 (the overall slope is a convex combination of two
slopes).} Let \(0<\eta_X^2<1\). With the fine-scale slope
\(\beta=\operatorname{Cov}(X,Y)/\sigma_X^2\), the zonal regression slope
\(\beta_B=C_B/B_X\) and the pooled within-zone slope
\(\beta_W=C_W/W_X\), \[
\beta=\eta_X^2\,\beta_B+(1-\eta_X^2)\,\beta_W .
\]

\emph{Proof.} \(0<\eta_X^2<1\) ensures \(B_X>0\) and \(W_X>0\), so both
slopes are defined. Substitute \(C_B=\beta_BB_X\) and \(C_W=\beta_WW_X\)
into Lemma 5.1 and divide by \(\sigma_X^2=B_X+W_X\). \(\blacksquare\)

\hypertarget{simpson-reversal}{%
\subsection{5.2 Simpson reversal}\label{simpson-reversal}}

\textbf{Corollary 5.4 (necessary and sufficient condition for Simpson
reversal).} Let \(B_XB_YW_XW_Y>0\). The fine-scale and pooled
within-zone correlations have opposite signs (\(\rho\,\rho_W<0\)) if and
only if \(\rho_B\rho_W<0\) and \[
\eta_X\eta_Y\,|\rho_B|>\sqrt{(1-\eta_X^2)(1-\eta_Y^2)}\;|\rho_W| .
\]

\emph{Proof.} The assumption ensures that \(\rho_B,\rho_W\) are defined.
Write \(u=\eta_X\eta_Y\rho_B\) and
\(v=\sqrt{(1-\eta_X^2)(1-\eta_Y^2)}\,\rho_W\); by Theorem 5.2,
\(\rho=u+v\), and \(v\) has the sign of \(\rho_W\)
(\(1-\eta_X^2,1-\eta_Y^2>0\)). If \(\rho\rho_W<0\), then \(v\neq0\) and
\(\rho\) has the sign of \(-v\), so \(u=\rho+(-v)\) is a sum of two
numbers of the same sign, has the sign of \(\rho\), and
\(|u|=|\rho|+|v|>|v|\); hence \(u,v\) have opposite signs,
i.e.~\(\rho_B\rho_W<0\), and \(|u|>|v|\). Conversely, if \(u,v\) have
opposite signs and \(|u|>|v|\), then \(\rho=u+v\) has the sign of \(u\)
and the opposite sign of \(v\), i.e.~\(\rho\rho_W<0\). \(\blacksquare\)

In the example of Simpson (1951), an association present within every
subgroup vanishes in the pooled table, i.e.~the two components of
Theorem 5.2 cancel exactly; a reversal of sign, as in Corollary 5.4,
occurs when they point in opposite directions and the between-zone
component dominates. These phenomena and the MAUP are manifestations of
one identity under different \(\mathcal{G}\).

\begin{figure}
\centering
\includegraphics[width=1\textwidth,height=\textheight]{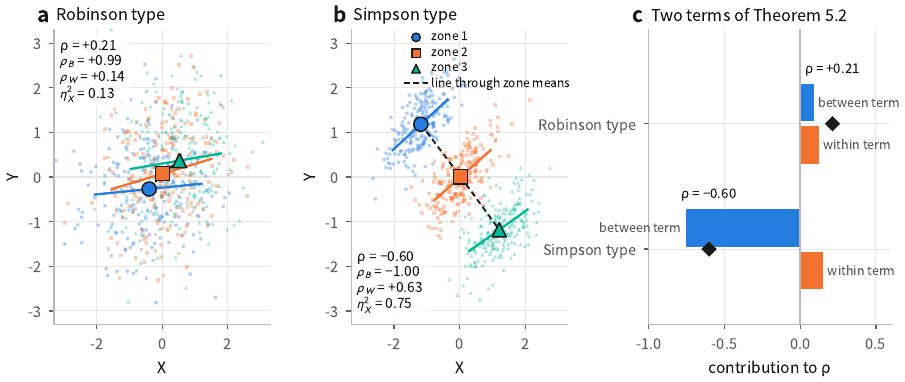}
\caption{The two typical cases of the fundamental identity. (a) Robinson
type: within each of three zones (distinguished by color and marker) the
correlation is weak (\(\rho_W=+0.14\); the within-zone regression lines
are flat), yet the zone means (large markers) are nearly collinear: the
zonal correlation is \(\rho_B=+0.99\) while the fine-scale correlation
is only \(\rho=+0.21\), because only about 13\% of the variance of \(X\)
and 6\% of that of \(Y\) lie between zones (\(\eta_X^2=0.13\),
\(\eta_Y^2=0.06\)). (b) Simpson type: each zone is positively correlated
internally (\(\rho_W=+0.63\)), yet the line through the zone means
(dashed) has a negative slope (\(\rho_B=-1.00\)); about 75\% of the
variance of each variable lies between zones
(\(\eta_X^2\approx\eta_Y^2\approx0.75\)), so the fine-scale correlation
is negative (\(\rho=-0.60\)). (c) The two terms of Theorem 5.2: in blue
the between-zone term \(\eta_X\eta_Y\rho_B\), in orange the within-zone
term \(\sqrt{(1-\eta_X^2)(1-\eta_Y^2)}\,\rho_W\); the diamond marks
their sum, the fine-scale correlation.}
\end{figure}

\hypertarget{an-impossibility-theorem}{%
\subsection{5.3 An impossibility
theorem}\label{an-impossibility-theorem}}

\textbf{Theorem 5.5 (only means are invariant to zoning).} Let \(\mu\)
be atomless, \(Z\in L^1(\mu;\mathbb{R}^d)\) and \(k\ge1\). Let
\(\mathcal{L}_k\) be the set of probability distributions on
\(\mathbb{R}^d\) supported on at most \(k\) points and
\(T:\mathcal{L}_k\to\mathcal{Y}\) any map (\(\mathcal{Y}\) any set). If
\(T\big(\operatorname{law}(\mathbb{E}[Z\mid\sigma(\mathcal{Z})])\big)\)
takes the same value for all \(k\)-zone zonings \(\mathcal{Z}\), then
that value is \(T(\delta_{\mathbb{E}Z})\). Consequently, if a zonal
statistic is invariant to \(k\)-zone zonings on every atomless support
measure space, its value depends only on the overall mean
\(\mathbb{E}Z\).

\emph{Proof.} (i) For a \(k\)-zone zoning \(\mathcal{Z}\), by Lemma
3.1(d) the law \(\sum_jp_j\delta_{\bar z_j}\) of
\(\mathbb{E}[Z\mid\sigma(\mathcal{Z})]\) is supported on at most \(k\)
points and belongs to \(\mathcal{L}_k\), so \(T\) is defined on it. (ii)
By Corollary 4.9 (with any \(p\in\Delta_k^\circ\)) there is a \(k\)-zone
zoning in which every zone mean equals \(\mathbb{E}Z\), whose law is
\(\delta_{\mathbb{E}Z}\). If \(T\) takes the same value on all
\(k\)-zone zonings, that value is \(T(\delta_{\mathbb{E}Z})\). (iii)
With \(f(m)=T(\delta_m)\), the invariant value on every atomless support
measure space is \(f(\mathbb{E}Z)\). \(\blacksquare\)

\textbf{Example.} Let \(k\ge2\) and let \(\mu\) and \(Z\) satisfy
assumption (H) of Section 6. The zonal variance vanishes under an
equal-mean zoning (Corollary 4.9), whereas Theorem 6.2 applied to any
\(V\neq0\) with \(p^\top V=0\) yields a zoning with zonal covariance
\(\varepsilon^2V^\top D_pV\), whose trace
\(\varepsilon^2\sum_jp_j|v_j|^2\) is positive. So the zonal variance is
not invariant to zoning. Zonal correlations and regression slopes are
not defined at \(\delta_{\mathbb{E}Z}\), and over the zonings on which
they are defined they vary as well: for \(d=2\) the zonal slope takes
every real value (Corollary 6.4) and, when \(k\ge3\), the zonal
correlation every value in \([-1,1]\) (Corollary 6.3). On the other
hand, \(\sum_jp_j\bar z_j=\mathbb{E}Z\) for every zoning (Lemma 3.1);
the zoning invariance of average effects in Section 11 is this same
fact.

Theorem 5.5 shows that the search for ``frame-independent'' association
statistics is futile. The meaningful question is how far zonal
statistics deviate and under what conditions, which is what Sections 7
and 8 answer.

\hypertarget{data-dependent-zonings-indeterminacy-theorems}{%
\section{6 Data-dependent zonings: indeterminacy
theorems}\label{data-dependent-zonings-indeterminacy-theorems}}

This section deals with zonings whose boundaries may be drawn with
reference to the values of the variables under study, for example
regionalizations driven by within-zone homogeneity or by a target
statistic. The conclusion is that once the zoning may depend on the
data, the zonal correlation is completely unconstrained by the
fine-scale correlation. For \(k\ge3\) every value can be attained
exactly with any prescribed zone weights (Theorem 6.2, Corollary 6.3);
in the unit square, requiring every zone to be connected only turns
``attained exactly'' into ``approached arbitrarily closely'' (Theorem
6.5, Corollary 6.6).

\hypertarget{assumption-and-a-basic-lemma}{%
\subsection{6.1 Assumption and a basic
lemma}\label{assumption-and-a-basic-lemma}}

\textbf{Assumption (H).} \(\mu\) is atomless;
\(Z\in L^1(\mu;\mathbb{R}^d)\); non-degeneracy: for every
\(a\in\mathbb{R}^d\setminus\{0\}\), \(a^\top Z\) is not almost surely
constant. Write \(D_p=\operatorname{diag}(p)\).

The non-degeneracy condition rules out only one situation: an exact
affine relation \(a^\top Z=c\) almost surely (for example \(Y=2X+1\), or
one variable almost surely constant), in which case zone means can only
lie in a fixed affine subspace.

\textbf{Lemma 6.1 (zone means can be moved in any direction).} Under
(H), the linear map \(L(h)=\mathbb{E}[hZ]\) maps
\(B_0=\{h\ \text{bounded measurable}:\ \mathbb{E}h=0\}\) onto
\(\mathbb{R}^d\).

\emph{Proof.} For bounded measurable \(h\), \(|hZ|\le\|h\|_\infty|Z|\),
so \(L(h)\) is defined and linear in \(h\); \(S=L(B_0)\) is a linear
subspace of \(\mathbb{R}^d\). If \(S\neq\mathbb{R}^d\), choose
\(a\neq0\) orthogonal to \(S\); then \(\mathbb{E}[hU]=a^\top L(h)=0\)
for all \(h\in B_0\), where \(U=a^\top Z\) is integrable and not almost
surely constant.

\emph{There is \(t\) with \(q=\mu(U>t)\in(0,1)\).} Let
\(F(t)=\mu(U\le t)\), which is non-decreasing and right-continuous, with
limits \(0\) at \(-\infty\) and \(1\) at \(+\infty\). If
\(F(t)\in\{0,1\}\) for all \(t\), let \(t^*=\inf\{t:F(t)=1\}\), a finite
number; \(F(t)=0\) for \(t<t^*\) and \(F(t^*)=1\) by right-continuity,
so \(\mu(U=t^*)=1\), contradicting the assumption that \(U\) is not
almost surely constant.

\emph{Computation.} Let \(h=\mathbf{1}\{U>t\}-q\in B_0\),
\(m_+=\mathbb{E}[U\mathbf{1}\{U>t\}]/q\) and
\(m_-=\mathbb{E}[U\mathbf{1}\{U\le t\}]/(1-q)\). Then
\(\mathbb{E}U=qm_++(1-q)m_-\) and \[
\mathbb{E}[hU]=qm_+-q\,\mathbb{E}U=qm_+-q^2m_+-q(1-q)m_-=q(1-q)(m_+-m_-).
\] On \(\{U>t\}\) we have \(U-t>0\), and this set has positive measure,
so \(\mathbb{E}[(U-t)\mathbf{1}\{U>t\}]>0\), i.e.~\(m_+>t\); likewise
\(m_-\le t\). Hence \(\mathbb{E}[hU]>0\), a contradiction.
\(\blacksquare\)

\hypertarget{unconstrained-zonings}{%
\subsection{6.2 Unconstrained zonings}\label{unconstrained-zonings}}

\textbf{Theorem 6.2 (exact indeterminacy for unconstrained zonings).}
Assume (H), \(k\ge2\), \(p\in\Delta_k^\circ\), and
\(V\in\mathbb{R}^{k\times d}\) with \(p^\top V=0\) (rows
\(v_1,\dots,v_k\)). Then there is \(\varepsilon_0>0\) such that for
every \(\varepsilon\in[0,\varepsilon_0]\) there is a \(k\)-zone zoning
\(\{A_j\}\) with zone weights \(p\) and \[
\mathbb{E}[Z\mid A_j]=\mathbb{E}Z+\varepsilon v_j\qquad(j=1,\dots,k).
\]

\emph{Proof.}

\begin{enumerate}
\def\labelenumi{\arabic{enumi}.}
\tightlist
\item
  By Lemma 6.1, choose \(u_1,\dots,u_d\in B_0\) with \(L(u_l)=e_l\) (the
  \(l\)-th unit vector). Put \(h_j=\sum_lV_{jl}u_l\). Since \(B_0\) is a
  linear space, \(h_j\in B_0\); \(L(h_j)=\sum_lV_{jl}e_l=v_j\); and
  pointwise
  \(\sum_jp_jh_j=\sum_l\big(\sum_jp_jV_{jl}\big)u_l=\sum_l(p^\top V)_lu_l=0\).
\item
  Let \(M=\max_j\|h_j\|_\infty\); take \(\varepsilon_0=1\) if \(M=0\)
  and \(\varepsilon_0=1/M\) otherwise. For
  \(\varepsilon\in[0,\varepsilon_0]\) put
  \(f_j=p_j(1+\varepsilon h_j)\). Then
  \(f_j\ge p_j(1-\varepsilon|h_j|)\ge p_j(1-\varepsilon M)\ge0\) and
  \(\sum_jf_j=\sum_jp_j+\varepsilon\sum_jp_jh_j=1\), so \(f\) is a fuzzy
  partition; moreover \[
  \int f_j\,d\mu=p_j(1+\varepsilon\,\mathbb{E}h_j)=p_j,\qquad \int f_jZ\,d\mu=p_j\mathbb{E}Z+\varepsilon p_j\,\mathbb{E}[h_jZ]=p_j(\mathbb{E}Z+\varepsilon v_j).
  \]
\item
  The vector measure \(\nu=(\mu,Z_1\mu,\dots,Z_d\mu)\) satisfies the
  assumptions of Theorem 4.6 (see the proof of Corollary 4.9). By
  Theorem 4.8 there is a partition \(\{A_j\}\) with
  \(\mu(A_j)=\int f_j\,d\mu=p_j>0\) and
  \(\int_{A_j}Z_i\,d\mu=\int f_jZ_i\,d\mu\). Dividing by \(p_j\) gives
  \(\mathbb{E}[Z\mid A_j]=\mathbb{E}Z+\varepsilon v_j\). The case
  \(\varepsilon=0\) is Corollary 4.9. \(\blacksquare\)
\end{enumerate}

\textbf{Corollary 6.3 (complete characterization of zonal correlation
matrices).} Assume (H), \(k\ge2\), and let \(p\in\Delta_k^\circ\) be
arbitrary. Over \(k\)-zone zonings with zone weights \(p\) and all zonal
variances positive, the (weighted) zonal correlation matrix ranges over
exactly the \(d\times d\) correlation matrices of rank at most \(k-1\).
In particular, for \(d=2\): with \(k=2\) the zonal correlation can only
be \(\pm1\), and both values occur; with \(k\ge3\) every value in
\([-1,1]\) is attained exactly, regardless of the fine-scale correlation
\(\rho\).

\emph{Proof.} We first record three facts of linear algebra. (1) For a
\(k\times d\) matrix \(M\),
\(\operatorname{rank}(M^\top D_pM)=\operatorname{rank}(M)\): if
\(M^\top D_pMx=0\) then \(x^\top M^\top D_pMx=\sum_jp_j(Mx)_j^2=0\), so
\(Mx=0\), and the two kernels coincide. (2) For an invertible diagonal
matrix \(S^{-1/2}\),
\(\operatorname{rank}(S^{-1/2}CS^{-1/2})=\operatorname{rank}(C)\). (3)
\(H=\{u\in\mathbb{R}^k:p^\top u=0\}\) is a \((k-1)\)-dimensional
subspace.

\emph{Necessity.} Let \(M_c\) be the \(k\times d\) matrix whose \(j\)-th
row is \(\bar z_j-\mathbb{E}Z\). By Lemma 3.1(d), \(p^\top M_c=0\),
i.e.~every column of \(M_c\) lies in \(H\), so
\(\operatorname{rank}M_c\le k-1\). The zonal covariance is
\(C=M_c^\top D_pM_c\); by (1) and (2) the zonal correlation matrix
\(S^{-1/2}CS^{-1/2}\) (\(S=\operatorname{diag}C\) with positive
diagonal) has rank \(\operatorname{rank}M_c\le k-1\).

\emph{Sufficiency.} Let the correlation matrix \(R\) (positive
semi-definite with unit diagonal) have rank \(r\le k-1\). By the
spectral decomposition \(R=\sum_{i\le r}\lambda_iw_iw_i^\top\)
(\(\lambda_i>0\), \(w_i\) orthonormal), and with
\(B\in\mathbb{R}^{r\times d}\) having \(i\)-th row
\(\sqrt{\lambda_i}\,w_i^\top\), \(R=B^\top B\). Equip \(\mathbb{R}^k\)
with the inner product \(\langle u,w\rangle_p=\sum_jp_ju_jw_j\)
(positive definite since \(p_j>0\)). Apply Gram--Schmidt with respect to
\(\langle\cdot,\cdot\rangle_p\) to a basis of \(H\) and take the first
\(r\) vectors as the columns of \(Q\in\mathbb{R}^{k\times r}\):
\(p^\top Q=0\) and \(Q^\top D_pQ=I_r\). Put \(V=QB\); then
\(p^\top V=0\) and \(V^\top D_pV=B^\top Q^\top D_pQB=R\). Applying
Theorem 6.2 with this \(V\) and any \(\varepsilon\in(0,\varepsilon_0]\):
the zone means are \(\mathbb{E}Z+\varepsilon v_j\), their weighted mean
is \(\mathbb{E}Z+\varepsilon V^\top p=\mathbb{E}Z\), and the zonal
covariance is
\(\sum_jp_j(\varepsilon v_j)(\varepsilon v_j)^\top=\varepsilon^2V^\top D_pV=\varepsilon^2R\)
with diagonal entries \(\varepsilon^2>0\), so the zonal correlation
matrix is \(R\).

\emph{\(d=2\).} The determinant of
\(\begin{pmatrix}1&\rho_B\\\rho_B&1\end{pmatrix}\) is \(1-\rho_B^2\), so
it has rank 1 if and only if \(\rho_B=\pm1\). For \(k=2\) the zonal
correlation matrix has rank at most 1, so \(\rho_B=\pm1\), and both
values are attained by sufficiency; for \(k\ge3\) every \(2\times2\)
correlation matrix has rank at most \(2\le k-1\). \(\blacksquare\)

\textbf{Corollary 6.4 (arbitrary zonal regression slopes).} Let \(d=2\),
\(Z=(X,Y)\), and assume (H). For any \(k\ge2\), any
\(p\in\Delta_k^\circ\) and any real \(b\) there is a \(k\)-zone zoning
with zone weights \(p\) whose weighted zonal regression of \(Y\) on
\(X\) has slope exactly \(b\).

\emph{Proof.} \(c=(p_2,-p_1,0,\dots,0)^\top\) satisfies \(p^\top c=0\)
and \(c\neq0\). Let the \(j\)-th row of \(V\) be \(v_j=c_j(1,b)\), so
\(p^\top V=0\). By Theorem 6.2 (\(\varepsilon\in(0,\varepsilon_0]\)),
the zone means are
\((\mathbb{E}X+\varepsilon c_j,\ \mathbb{E}Y+\varepsilon bc_j)\) with
weighted mean \(\mathbb{E}Z\), and the weighted zonal regression slope
is \[
\frac{\sum_jp_j(\bar x_j-\mathbb{E}X)(\bar y_j-\mathbb{E}Y)}{\sum_jp_j(\bar x_j-\mathbb{E}X)^2}=\frac{\varepsilon^2b\sum_jp_jc_j^2}{\varepsilon^2\sum_jp_jc_j^2}=b,
\] the denominator being positive since \(c\neq0\) and \(p_j>0\).
\(\blacksquare\)

\textbf{Remark.} The case \(k=2\) of Corollary 6.3 shows that the zonal
correlation of any two-zone comparison is \(\pm1\): two points always
lie on a line, which has nothing to do with fine-scale relationships.
With more than two zones, a data-dependent zoning can give the zonal
correlation any value. The zone means constructed in Theorem 6.2 can be
taken arbitrarily close to the overall mean (the offsets
\(\varepsilon v_j\) shrink with \(\varepsilon\)): manufacturing a strong
zonal correlation does not require large between-zone differences, only
that the differences line up in a prescribed direction.

\hypertarget{connected-zonings}{%
\subsection{6.3 Connected zonings}\label{connected-zonings}}

Let \(D=(0,1)^2\) be the open unit square with Lebesgue measure
\(\lambda\), let \(\mu=q\lambda\) be a probability measure on \(D\) with
density \(q\in L^1(\lambda)\) (hence atomless), and let
\(Z\in L^1(\mu;\mathbb{R}^d)\) satisfy the non-degeneracy condition of
(H). A \textbf{polygonal connected zoning} consists of pairwise disjoint
path-connected open polygonal sets \(A_1,\dots,A_k\) in \(D\) with
\(\mu(A_j)>0\) for every \(j\) and
\(\lambda\big(D\setminus\bigcup_jA_j\big)=0\), the remainder lying in a
finite union of segments (hence of \(\mu\)-measure zero).

\textbf{Theorem 6.5 (connectivity does not remove indeterminacy).} Let
\(k\ge2\), \(p\in\Delta_k^\circ\), \(V\) with \(p^\top V=0\), and
\(\varepsilon\in(0,\varepsilon_0)\), where \(\varepsilon_0>0\) is such
that the conclusion of Theorem 6.2 holds for this \(V\). Then for every
\(\eta>0\) there is a polygonal connected \(k\)-zone zoning with, for
all \(j\), \[
|\mu(A_j)-p_j|<\eta,\qquad\big\|\mathbb{E}[Z\mid A_j]-(\mathbb{E}Z+\varepsilon v_j)\big\|<\eta .
\]

\textbf{Corollary 6.6.} For \(d=2\) and \(k\ge3\), the zonal
correlations of polygonal connected \(k\)-zone zonings are dense in
\([-1,1]\); for \(k=2\), both \(\pm1\) are attained exactly.

The idea of the proof is to approximate the fuzzy partition of Theorem
6.2 by a step function on a grid, and then to use a ``\(k\)-lane
serpentine path'' through the whole square to allocate every small cell
to the \(k\) zones in proportion; each lane is connected from end to
end, so each zone is connected. Figure 4(a)--(c) shows the construction:
for the same pair of variables (fine-scale correlation 0.878), without
connectivity the three zones are interleaved, with connectivity they are
long thin strips winding side by side through the domain, and in both
cases the zonal correlation can be pushed to \(\pm0.9\). The complete
proof, in nine steps, is given in Section 6.4.

\begin{figure}
\centering
\includegraphics[width=1\textwidth,height=\textheight]{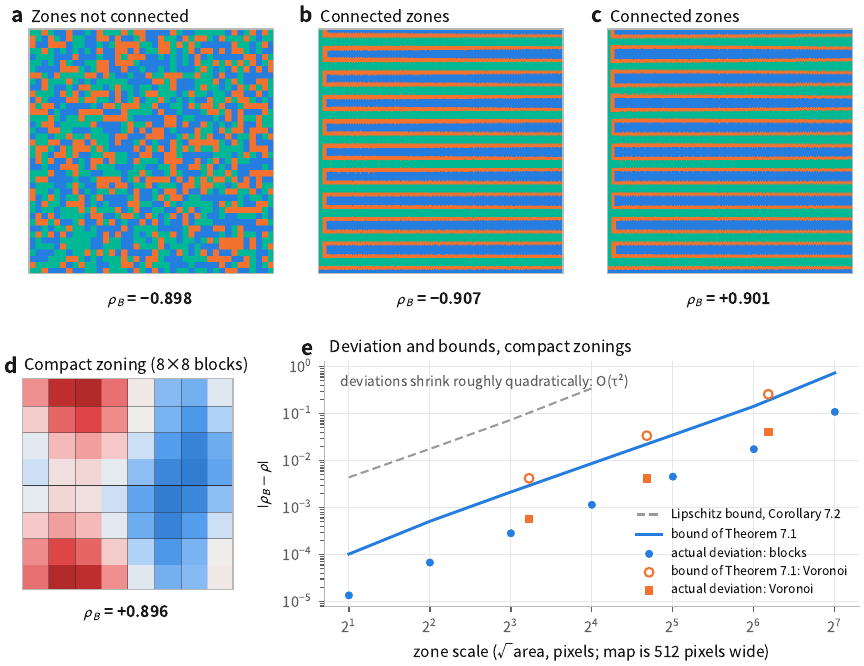}
\caption{Data-dependent versus compact zonings. The same pair of
variables \(X,Y\) (fine-scale correlation \(\rho=0.878\)); in (a)--(c)
the three zones are shown in three colors. (a) Data-dependent,
connectivity not required: obtained by randomized rounding of the fuzzy
partition of Theorem 6.2; the three zones are interleaved (only a
\(40\times40\)-pixel window is shown); zonal correlation
\(\rho_B=-0.898\). (b), (c) Data-dependent, every zone connected:
obtained with the serpentine-lane construction of Theorem 6.5 (only the
south-west 1/36 of the map is shown); \(\rho_B=-0.907\) and \(+0.901\).
(d) A compact zoning: an \(8\times8\) grid of square blocks (side 64
pixels), colored by zone mean; \(\rho_B=+0.896\), almost equal to the
fine-scale correlation. (e) Deviations and bounds for compact zonings
(log--log axes); the horizontal axis is the zone scale (square root of
zone area in pixels; the whole map is \(512\times512\) pixels). Blue
dots and the blue line are the actual deviations \(|\rho_B-\rho|\) for
square blocks and the bound of Theorem 7.1; orange squares and open
circles are the actual deviations and bounds for random Voronoi zonings
(50, 400 and 3000 zones); the grey dashed line is the Lipschitz bound of
Corollary 7.2 (defined only for \(\tau<1\) and clearly looser).
Deviations and bounds both grow roughly quadratically with zone scale,
and the bounds always lie above the actual deviations.}
\end{figure}

\hypertarget{proof-of-theorem-6.5}{%
\subsection{6.4 Proof of Theorem 6.5}\label{proof-of-theorem-6.5}}

Notation is as in Section 6.3: \(D=(0,1)^2\), \(\lambda\) Lebesgue
measure, \(\mu=q\lambda\), \(Z\in L^1(\mu;\mathbb{R}^d)\)
non-degenerate, \(\varepsilon\in(0,\varepsilon_0)\).

\textbf{Step 0: reduction to finitely many integrals.} For a measurable
set \(A\) write \(m(A)=\mu(A)=\int\mathbf{1}_Aq\,d\lambda\) and
\(s(A)=\int_AZ\,d\mu=\int\mathbf{1}_AZq\,d\lambda\). The targets are
\(m_j^*=p_j\) and \(s_j^*=p_j(\mathbb{E}Z+\varepsilon v_j)\). If
\(|m(A_j)-m_j^*|\le\delta_1\le p_j/2\) and
\(\|s(A_j)-s_j^*\|_\infty\le\delta_1\), then from \[
\frac{s}{m}-\frac{s^*}{m^*}=\frac{s-s^*}{m}+s^*\,\frac{m^*-m}{m\,m^*}
\] and \(m\ge p_j/2\) we get
\(\big\|\mathbb{E}[Z\mid A_j]-(\mathbb{E}Z+\varepsilon v_j)\big\|_\infty\le\frac{2\delta_1}{p_j}+\frac{2\|s_j^*\|_\infty\,\delta_1}{p_j^2}\).
Hence, given \(\eta\), we can choose \(\eta'>0\) so small that \[
\Big|\int(\mathbf{1}_{A_j}-f_j)\,\psi\,d\lambda\Big|<\eta'\qquad\text{for all }j\text{ and all }\psi\in\Psi=\{q,Z_1q,\dots,Z_dq\}\tag{6.1}
\] implies both conclusions of the theorem (the Euclidean and sup norms
differ by at most a factor \(\sqrt d\), which is absorbed into the
choice of \(\eta'\)). Here we used \(\int f_jq\,d\lambda=m_j^*\) and
\(\int f_jZq\,d\lambda=s_j^*\) (Step 1).

\textbf{Step 1: the target fuzzy partition.} By Theorem 6.2 (with
\(\varepsilon_0\)) there is a zoning \(\{B_j\}\) with \(\mu(B_j)=p_j\)
and \(\mathbb{E}[Z\mid B_j]=\mathbb{E}Z+\varepsilon_0v_j\). Put
\(t=\varepsilon/\varepsilon_0\in(0,1)\) and
\(f_j=t\,\mathbf{1}_{B_j}+(1-t)\,p_j\). Then \(f=(f_1,\dots,f_k)\) is a
fuzzy partition, \(\int f_j\,d\mu=tp_j+(1-t)p_j=p_j\),
\(\int f_jZ\,d\mu=tp_j(\mathbb{E}Z+\varepsilon_0v_j)+(1-t)p_j\,\mathbb{E}Z=p_j(\mathbb{E}Z+\varepsilon v_j)\),
and \(c_0:=(1-t)\min_jp_j\le f_j\le1\). Put
\(\phi=(1+|Z|)q\in L^1(\lambda)\). Every \(\psi\in\Psi\) satisfies
\(|\psi|\le\phi\), so for bounded measurable \(u,w\), \[
\Big|\int(u-w)\psi\,d\lambda\Big|\le\int|u-w|\,\phi\,d\lambda=:\|u-w\|.\tag{6.2}
\]

\textbf{Step 2: step-function approximation on a grid.} Put
\(\eta''=c_0\eta'/8\). \(\phi\lambda\) is a finite Borel measure on the
open set \(D\), so by (S5) \(C_c(D)\) is dense in \(L^1(\phi\lambda)\);
choose \(F_j\in C_c(D)\) with \(\|f_j-F_j\|<\eta''\). Divide \(D\) into
\(n^2\) half-open grid squares of side \(1/n\) (their boundaries are
\(\lambda\)-null) and let \(G_j\) take, on each grid square, the value
of \(F_j\) at its center. Every point of a grid square is within
\(\sqrt2/(2n)\) of its center, so
\(|F_j-G_j|\le\omega_{F_j}(\sqrt2/(2n))\) (\(\omega\) the modulus of
continuity) and
\(\|F_j-G_j\|\le\omega_{F_j}(\sqrt2/(2n))\int\phi\,d\lambda<\eta''\) for
\(n\) large. Put \(\hat g_j=\min(\max(G_j,c_0),1)\). The clipping map
\(x\mapsto\min(\max(x,c_0),1)\) is 1-Lipschitz and fixes \([c_0,1]\),
and \(f_j\in[c_0,1]\), so \(|\hat g_j-f_j|\le|G_j-f_j|\). Let
\(S=\sum_i\hat g_i\in[kc_0,k]\) and \(g_j=\hat g_j/S\). Since
\(\sum_if_i=1\), \[
\begin{aligned}g_j-f_j&=\frac{\hat g_j-f_jS}{S}=\frac{(\hat g_j-f_j)-f_j\sum_i(\hat g_i-f_i)}{S},\\|g_j-f_j|&\le\frac{|\hat g_j-f_j|+\sum_i|\hat g_i-f_i|}{kc_0}\le\frac{2}{kc_0}\sum_i|G_i-f_i| .\end{aligned}
\] Hence
\(\|g_j-f_j\|\le\frac{2}{kc_0}\sum_i(\|G_i-F_i\|+\|F_i-f_i\|)<\frac{2}{kc_0}\cdot k\cdot2\eta''=\frac{4\eta''}{c_0}\).
Moreover \(g_j\ge c_0/S\ge c_0/k=:c_1\), \(\sum_jg_j=1\), and each
\(g_j\) is constant on every grid square. By the choice of \(\eta''\),
(6.2) gives \[
\Big|\int(g_j-f_j)\psi\,d\lambda\Big|<\eta'/2\qquad(\psi\in\Psi).\tag{6.3}
\] Write \(\gamma(Q)=(g_1,\dots,g_k)|_Q\) for the value of \(g\) on a
grid square \(Q\); its components are at least \(c_1\) and sum to \(1\).

\textbf{Step 3: the boustrophedon path.} Take an integer \(M'\ge1\), put
\(N=mM'\ge2\), and divide \(D\) into \(N^2\) cells of side \(1/N\) (each
lying in a unique grid square). Row \(r\) (\(r=1,\dots,N\) from south to
north) is traversed west to east for odd \(r\) and east to west for even
\(r\); joining the rows in order gives a sequence of cells
\(C_1,C_2,\dots,C_{N^2}\) in which consecutive cells share an edge:
within a row a vertical edge, and between the last cell of one row and
the first cell of the next (same column) a horizontal edge. The last
cell of every row except the last and the first cell of every row except
the first are \textbf{turning cells}, \(2(N-1)\) in all; the others are
\textbf{straight cells}. There are four kinds of turning cells:

\begin{itemize}
\tightlist
\item
  T1: end of an odd (eastbound) row, entered through the west edge,
  exited through the north edge;
\item
  T2: start of an even (westbound) row, entered through the south edge,
  exited through the west edge;
\item
  T3: end of an even row, entered through the east edge, exited through
  the north edge;
\item
  T4: start of an odd row (\(r>1\)), entered through the south edge,
  exited through the east edge.
\end{itemize}

Let \(e_t\) be the edge shared by \(C_t\) and \(C_{t+1}\); the entry of
\(C_1\) is its west edge \(e_0\), and the exit \(e_{N^2}\) of
\(C_{N^2}\) is the vertical edge of the last cell of the last row lying
on \(\partial D\). \textbf{Proportion vectors}:
\(\pi(e_0)=\gamma(C_1)\); for \(t\ge1\), \(\pi(e_t)=\gamma(C_t)\) if
\(C_t\) is a straight cell and \(\pi(e_t)=\pi(e_{t-1})\) if it is a
turning cell. By induction every \(\pi(e_t)\) is a value of \(\gamma\),
with components at least \(c_1\). Write \(E_j(e)=\sum_{i\le j}\pi(e)_i\)
(\(E_0=0\), \(E_k=1\)).

\textbf{Step 4: the lateral coordinate.} In each cell use local
coordinates \((x,y)\in[0,1]^2\) (\(x\) eastward, \(y\) northward) and
define the lateral coordinate \(\ell\), measured from the right-hand
side of the direction of travel:

\begin{longtable}[]{@{}
  >{\raggedright\arraybackslash}p{(\columnwidth - 6\tabcolsep) * \real{0.2286}}
  >{\raggedright\arraybackslash}p{(\columnwidth - 6\tabcolsep) * \real{0.2571}}
  >{\raggedright\arraybackslash}p{(\columnwidth - 6\tabcolsep) * \real{0.2571}}
  >{\raggedright\arraybackslash}p{(\columnwidth - 6\tabcolsep) * \real{0.2571}}@{}}
\caption{Lateral coordinates of the six cell types}\tabularnewline
\toprule\noalign{}
\begin{minipage}[b]{\linewidth}\raggedright
cell type
\end{minipage} & \begin{minipage}[b]{\linewidth}\raggedright
\(\ell(x,y)\)
\end{minipage} & \begin{minipage}[b]{\linewidth}\raggedright
\(\ell\) on the entry edge
\end{minipage} & \begin{minipage}[b]{\linewidth}\raggedright
\(\ell\) on the exit edge
\end{minipage} \\
\midrule\noalign{}
\endfirsthead
\toprule\noalign{}
\begin{minipage}[b]{\linewidth}\raggedright
cell type
\end{minipage} & \begin{minipage}[b]{\linewidth}\raggedright
\(\ell(x,y)\)
\end{minipage} & \begin{minipage}[b]{\linewidth}\raggedright
\(\ell\) on the entry edge
\end{minipage} & \begin{minipage}[b]{\linewidth}\raggedright
\(\ell\) on the exit edge
\end{minipage} \\
\midrule\noalign{}
\endhead
\bottomrule\noalign{}
\endlastfoot
eastbound straight & \(y\) & west: \(y\) & east: \(y\) \\
westbound straight & \(1-y\) & east: \(1-y\) & west: \(1-y\) \\
T1 (west in, north out) & \(1-\max(x,1-y)\) & west: \(y\) & north:
\(1-x\) \\
T2 (south in, west out) & \(1-\max(x,y)\) & south: \(1-x\) & west:
\(1-y\) \\
T3 (east in, north out) & \(\max(1-x,1-y)\) & east: \(1-y\) & north:
\(1-x\) \\
T4 (south in, east out) & \(\max(1-x,y)\) & south: \(1-x\) & east:
\(y\) \\
\end{longtable}

The edge values follow by substitution, e.g.~on the west edge \(x=0\) of
T1, \(1-\max(0,1-y)=y\), and on its north edge \(y=1\),
\(1-\max(x,0)=1-x\). \textbf{Adjacent cells agree on their shared edge}:
the two cells on either side of a shared vertical edge belong to the
same row and travel in the same direction, and the table gives \(y\) on
both sides (eastbound) or \(1-y\) on both sides (westbound); a shared
horizontal edge joins the end of one row (T1 or T3, exit north,
\(\ell=1-x\)) to the start of the next (T2 or T4, entry south,
\(\ell=1-x\)), and the two cells have the same \(x\) coordinate (same
column). On every entry and exit edge, in particular on every \(e_t\),
\(\ell\) is an affine parameter of the edge ranging over \([0,1]\); on
the other edges of a cell \(\ell\) is constant, equal to \(0\) or \(1\).
In the four kinds of turning cells, \(\ell\) is a maximum of two affine
functions (or \(1\) minus such a maximum), and its level sets
\(\{\ell=c\}\) are L-shaped polylines: in T1, for example,
\(\{\max(x,1-y)=1-c\}\) consists of the vertical segment
\(\{x=1-c,\ y\ge c\}\) and the horizontal segment \(\{y=c,\ x\le1-c\}\),
with corner at \((1-c,c)\).

\textbf{Step 5: lanes.} Let \(C=C_t\) have entry proportions
\(\pi^{\mathrm{in}}=\pi(e_{t-1})\) and exit proportions
\(\pi^{\mathrm{out}}=\pi(e_t)\), with partial sums \(E^{\mathrm{in}}_j\)
and \(E^{\mathrm{out}}_j\). In a straight cell put
\(E_j(\theta)=(1-\theta)E_j^{\mathrm{in}}+\theta E_j^{\mathrm{out}}\),
where \(\theta\) is the coordinate along the direction of travel
(\(\theta=x\) in eastbound cells, \(\theta=1-x\) in westbound ones); in
a turning cell \(\pi^{\mathrm{out}}=\pi^{\mathrm{in}}\) and
\(E_j\equiv E^{\mathrm{in}}_j\). Lane \(j\) is \[
\Lambda_j(C)=\{(x,y)\in(0,1)^2:\ E_{j-1}<\ell(x,y)<E_j\}\quad(\text{mapped back to the cell}).
\] (i) In a straight cell
\(E_j(\theta)-E_{j-1}(\theta)=(1-\theta)\pi^{\mathrm{in}}_j+\theta\pi^{\mathrm{out}}_j\ge c_1>0\),
so the lane is an open trapezoid bounded by two straight lines and the
two vertical edges; it is convex, hence path-connected, and its area is
\(\frac12(\pi^{\mathrm{in}}_j+\pi^{\mathrm{out}}_j)\lambda(C)\). (ii) In
a turning cell the lane is the open region between two L-shaped level
curves: the lane containing the inner corner of the turn is an open
square (lane \(k\) in T1 and T2, where \(\ell=1\) only at that corner;
lane \(1\) in T3 and T4, where \(\ell=0\) only at that corner), and the
others are open L-shaped bands; all are polygons and path-connected (an
L-shaped band is the union of two open rectangles meeting in an open
square). (iii) On the shared edge \(e_t\), the exit values of \(C_t\)
and the entry values of \(C_{t+1}\) are both \(E_j(e_t)\), and \(\ell\)
agrees, so lane \(j\) of both cells has the same open interval
\(I_j(e_t)=\{z\in e_t:E_{j-1}(e_t)<\ell(z)<E_j(e_t)\}\) as part of its
boundary, of length \(\pi(e_t)_j/N\ge c_1/N\).

\textbf{Step 6: the zones and their properties.} Let
\(U_j=\bigcup_t\Lambda_j(C_t)\) (an open set) and
\(A_j=\operatorname{int}\overline{U_j}\).

\begin{enumerate}
\def\labelenumi{(\alph{enumi})}
\item
  \textbf{\(U_j\cup\bigcup_{t=1}^{N^2-1}I_j(e_t)\subseteq A_j\).}
  \(\ell\) and \(E_j\) are continuous on each closed cell. For
  \(1\le t\le N^2-1\) and \(z\in I_j(e_t)\) (in the relative interior of
  the shared edge) there is \(r>0\) with
  \(B(z,r)\subseteq\overline{C_t}\cup\overline{C_{t+1}}\) such that
  \(E_{j-1}<\ell<E_j\) on both \(B(z,r)\cap\overline{C_t}\) and
  \(B(z,r)\cap\overline{C_{t+1}}\) (each according to its own cell).
  Hence
  \(B(z,r)\subseteq\Lambda_j(C_t)\cup I_j(e_t)\cup\Lambda_j(C_{t+1})\subseteq\overline{U_j}\),
  and \(z\in A_j\).
\item
  \textbf{\(A_j\) is path-connected.} By the argument in (a),
  \(V_t=\Lambda_j(C_t)\cup I_j(e_t)\cup\Lambda_j(C_{t+1})\) is open and
  path-connected (each lane is path-connected and meets a disk around
  any point of \(I_j(e_t)\)). Chaining \(V_1,\dots,V_{N^2-1}\) shows
  that \(W_j=U_j\cup\bigcup_{t=1}^{N^2-1}I_j(e_t)\) is path-connected.
  For any \(z\in A_j\) there is a disk
  \(B(z,r)\subseteq\overline{U_j}\); the open set \(B(z,r)\) lies in the
  interior of \(\overline{U_j}\), i.e.~\(B(z,r)\subseteq A_j\), and it
  meets \(U_j\) (it lies in the closure of \(U_j\)). So \(z\) can be
  joined to \(W_j\) by a segment inside \(A_j\), and \(A_j\) is
  path-connected.
\item
  \textbf{Pairwise disjointness.} Different lanes are disjoint (within a
  cell they are separated by the values of \(\ell\), and different cells
  have disjoint interiors), so \(U_i\cap U_j=\varnothing\) for
  \(i\neq j\). If \(z\in A_i\cap A_j\), take a disk
  \(B\subseteq\overline{U_i}\cap\overline{U_j}\); \(O=B\cap U_i\) is a
  non-empty open set (\(B\) lies in the closure of \(U_i\)) with
  \(O\subseteq\overline{U_j}\), so \(O\cap U_j\neq\varnothing\),
  contradicting \(O\subseteq U_i\).
\item
  \textbf{Polygonality and full measure.} \(\overline{U_j}\) is a finite
  union of closed polygons, and the boundary of \(A_j\) lies in a finite
  union of segments. \(D\setminus\bigcup_jU_j\) lies in the union of the
  cell edges and the level curves \(\{\ell=E_j\}\), a finite union of
  segments, which is \(\lambda\)-null and, since \(\mu\ll\lambda\),
  \(\mu\)-null. Also \(\overline{U_j}\subseteq[0,1]^2\), so
  \(A_j\subseteq(0,1)^2=D\).
\item
  \textbf{Within each cell, \(A_j\) and lane \(j\) differ by a null
  set.} \(A_j\supseteq\Lambda_j(C)\), while \(A_j\) is disjoint from
  \(A_i\supseteq\Lambda_i(C)\) for \(i\neq j\), and the interior of the
  cell minus its lanes consists of finitely many segments; hence
  \(\lambda\big((A_j\cap C)\,\triangle\,\Lambda_j(C)\big)=0\).
\end{enumerate}

\textbf{Step 7: weak convergence.} Call a straight cell \(C\)
\textbf{regular} if both its entry and exit proportions equal the
\(\gamma\) of its grid square. By step 5(i) and step 6(e), for a regular
cell \[
\int_C\mathbf{1}_{A_j}\,d\lambda=\lambda(\Lambda_j(C))=\gamma_j\,\lambda(C)=\int_Cg_j\,d\lambda .
\] The exit proportions of a straight cell are always its own
\(\gamma\); its entry proportions are the exit proportions of the
preceding cell, which equal \(\gamma\) if the preceding cell is a
straight cell in the same grid square. So non-regular cells are of three
kinds only: turning cells (\(2(N-1)\) of them); the first straight cell
after a turn (at most one per row, at most \(N\)); and straight cells
whose predecessor is a straight cell in another grid square (each row
crosses \(n-1\) vertical grid lines, at most \(N(n-1)\)). The entry of
\(C_1\) is \(\gamma(C_1)\) by definition. The total area of non-regular
cells is at most \(\big(2(N-1)+N+N(n-1)\big)/N^2\le(n+2)/N\).

For \(\psi\in C_c(D)\): on a regular cell \(C\), the display above gives
\(\int_C(\mathbf{1}_{A_j}-g_j)\,d\lambda=0\), so \[
\Big|\int_C(\mathbf{1}_{A_j}-g_j)\psi\,d\lambda\Big|=\Big|\int_C(\mathbf{1}_{A_j}-g_j)(\psi-\psi(x_C))\,d\lambda\Big|\le\omega_\psi\big(\tfrac{\sqrt2}{2N}\big)\lambda(C),
\] with \(x_C\) the center of the cell; the sum over regular cells is at
most \(\omega_\psi(\sqrt2/(2N))\). Non-regular cells contribute at most
\(\|\psi\|_\infty(n+2)/N\). Both tend to zero as \(M'\to\infty\) (\(n\)
fixed, \(N\to\infty\)). For general \(\psi\in L^1(\lambda)\) and any
\(\delta>0\), choose \(\psi_c\in C_c(D)\) with
\(\|\psi-\psi_c\|_{L^1(\lambda)}<\delta\); since
\(|\mathbf{1}_{A_j}-g_j|\le1\), \[
\limsup_{M'\to\infty}\Big|\int(\mathbf{1}_{A_j}-g_j)\psi\,d\lambda\Big|\le\limsup_{M'\to\infty}\Big|\int(\mathbf{1}_{A_j}-g_j)\psi_c\,d\lambda\Big|+\delta=\delta .
\]

\textbf{Step 8: conclusion.} First fix \(n\) and \(g\) as in step 2 (so
that (6.3) holds); then, by step 7 applied to the finite set \(\Psi\)
and \(j=1,\dots,k\), choose \(M'\) so large that
\(|\int(\mathbf{1}_{A_j}-g_j)\psi\,d\lambda|<\eta'/2\). Adding (6.3)
gives (6.1), and step 0 completes the proof. \(\blacksquare\)

\hypertarget{proof-of-corollary-6.6}{%
\subsection{6.5 Proof of Corollary 6.6}\label{proof-of-corollary-6.6}}

\textbf{\(k\ge3\).} For a target \(\rho^*\in[-1,1]\),
\(R=\begin{pmatrix}1&\rho^*\\\rho^*&1\end{pmatrix}\) has rank at most
\(2\le k-1\); take the \(V\) of Corollary 6.3. Write the zonal
correlation as a function of the zone masses and zone means: \[
\mathsf{R}\big((m_j,\bar z_j)_j\big)=\frac{\sum_jm_j(\bar x_j-\bar x)(\bar y_j-\bar y)}{\sqrt{\sum_jm_j(\bar x_j-\bar x)^2\,\sum_jm_j(\bar y_j-\bar y)^2}},\qquad \bar x=\frac{\sum_jm_j\bar x_j}{\sum_jm_j},\ \ \bar y=\frac{\sum_jm_j\bar y_j}{\sum_jm_j},
\] which is continuous wherever the denominator is positive. Fix
\(\varepsilon\in(0,\varepsilon_0)\) and \(\delta>0\). At the target
configuration \((p_j,\mathbb{E}Z+\varepsilon v_j)\) both weighted zonal
variances equal \(\varepsilon^2\) (the diagonal of \(V^\top D_pV=R\) is
\(1\)), so the denominator is \(\varepsilon^2>0\), and the value is
\(\rho^*\). By continuity there is \(\eta\in(0,\min_jp_j)\) such that
the denominator is positive and \(|\mathsf{R}-\rho^*|<\delta\) whenever
\(|m_j-p_j|<\eta\) and
\(\|\bar z_j-(\mathbb{E}Z+\varepsilon v_j)\|<\eta\) for all \(j\).
Theorem 6.5 supplies a polygonal connected zoning with this property
(for an actual zoning \(\sum_jm_j=1\)), and its zonal correlation is
within \(\delta\) of \(\rho^*\).

\textbf{\(k=2\).} For a two-zone configuration put
\(\Delta x=\bar x_1-\bar x_2\) and \(\Delta y=\bar y_1-\bar y_2\). Since
\(\bar x_1-\bar x=m_2\Delta x\) and \(\bar x_2-\bar x=-m_1\Delta x\)
(\(m_1+m_2=1\)), the zonal covariance is
\(m_1m_2^2\Delta x\Delta y+m_2m_1^2\Delta x\Delta y=m_1m_2\Delta x\Delta y\),
and the two zonal variances are \(m_1m_2\Delta x^2\) and
\(m_1m_2\Delta y^2\). When the variances are positive, \[
\rho_B=\frac{\Delta x\Delta y}{|\Delta x|\,|\Delta y|}=\operatorname{sign}(\Delta x\Delta y).
\] Take \(c=(p_2,-p_1)\) and \(v_j=c_j(1,\pm1)\); at the target
configuration \(\Delta x=\varepsilon(c_1-c_2)=\varepsilon\) and
\(\Delta y=\pm\varepsilon\). By Theorem 6.5 take a connected zoning
whose zone means are within \(\varepsilon/4\) of the targets; then
\(\Delta x\) and \(\Delta y\) change by less than \(\varepsilon/2\),
keep their signs and remain non-zero, so \(\rho_B=\pm1\) exactly.
\(\blacksquare\)

The connected zones constructed in Theorem 6.5 are long and thin: each
is made of strips running across the whole domain, and its diameter is
close to that of the domain itself. What constrains the zonal
correlation is compactness (Section 7), not connectivity. The case of
discrete basic units (each unit assigned to a zone as a whole) is
discussed in Section 13.3.

\hypertarget{compact-zonings-an-exact-error-bound-and-a-second-order-rate}{%
\section{7 Compact zonings: an exact error bound and a second-order
rate}\label{compact-zonings-an-exact-error-bound-and-a-second-order-rate}}

``Smaller zones, smaller aggregation error'' only gives a direction.
This section answers how well the error can be controlled. How far the
zonal correlation can deviate from the fine-scale correlation \(\rho\)
is bounded in terms of \(\rho\) and the two within-zone variance shares
\(a=1-\eta_X^2\) and \(b=1-\eta_Y^2\) alone (Theorem 7.1); the bound is
attained whenever it does not exceed \(1-|\rho|\) (Proposition 7.4), and
it can be computed from fine-scale data alone, without any model or
smoothness assumption. For spatially Lipschitz fields the within-zone
variance shares are in turn controlled by the square of the ratio of
zone diameter to variation length, so the error bound is of second order
in the zone diameter: halving the largest zone diameter reduces the
leading term of the bound to one quarter (Corollaries 7.2 and 7.3), and
the second order cannot be improved (Example 7.5).

\hypertarget{an-exact-error-bound}{%
\subsection{7.1 An exact error bound}\label{an-exact-error-bound}}

\textbf{Theorem 7.1 (exact error bound).} In the notation of Section 5,
let \(a=W_X/\sigma_X^2<1\), \(b=W_Y/\sigma_Y^2<1\) and
\(s=\sqrt{(1-a)(1-b)}\). Then \[
\rho_B=\frac{\rho-\sqrt{ab}\,\rho_W}{s},\qquad
|\rho_B-\rho|\le\frac{|\rho|\,(1-s)+\sqrt{ab}}{s}=:\mathcal{B}(a,b;\rho),
\] with the convention \(\sqrt{ab}\,\rho_W=0\) when \(ab=0\).
\(\mathcal{B}\) is non-decreasing in both \(a\) and \(b\), so if
\(a\le\alpha<1\) and \(b\le\beta<1\) are known, \(\alpha\) and \(\beta\)
may be substituted.

\emph{Proof.} \(a<1\) is equivalent to \(B_X=(1-a)\sigma_X^2>0\), and
likewise \(B_Y>0\), so \(\rho_B\) is defined. By Theorem 5.2
(\(\eta_X^2=1-a\), \(\eta_Y^2=1-b\)),
\(\rho=s\rho_B+\sqrt{ab}\,\rho_W\), the second term vanishing by the
convention of Theorem 5.2 when \(ab=0\). Solving for \(\rho_B\) gives
the first formula. Hence \[
\rho_B-\rho=\rho\Big(\frac1s-1\Big)-\frac{\sqrt{ab}}{s}\,\rho_W,
\] and \(|\rho_W|\le1\) gives
\(|\rho_B-\rho|\le|\rho|(1/s-1)+\sqrt{ab}/s=\mathcal{B}\). Monotonicity:
\(1/s=\big((1-a)(1-b)\big)^{-1/2}\) increases in \(a\) and \(b\), and
\(\sqrt{ab}/s=\sqrt{a/(1-a)}\cdot\sqrt{b/(1-b)}\) is a product of two
non-negative increasing functions, hence non-decreasing.
\(\blacksquare\)

\textbf{Remark (when the bound is informative).} The trivial bound is
\(|\rho_B-\rho|\le1+|\rho|\). For \(a=b\),
\(\mathcal{B}=a(1+|\rho|)/(1-a)\), which is smaller than the trivial
bound if and only if \(a<1/2\). So Theorem 7.1 is uninformative when
more than half of the variation of both variables lies within zones; the
more homogeneous the zones, the tighter the bound. Table 3 lists the
bound for \(a=b\).

\begin{longtable}[]{@{}
  >{\centering\arraybackslash}p{(\columnwidth - 8\tabcolsep) * \real{0.3571}}
  >{\centering\arraybackslash}p{(\columnwidth - 8\tabcolsep) * \real{0.1607}}
  >{\centering\arraybackslash}p{(\columnwidth - 8\tabcolsep) * \real{0.1607}}
  >{\centering\arraybackslash}p{(\columnwidth - 8\tabcolsep) * \real{0.1607}}
  >{\centering\arraybackslash}p{(\columnwidth - 8\tabcolsep) * \real{0.1607}}@{}}
\caption{Upper bound \(\mathcal{B}\) on \(|\rho_B-\rho|\) from Theorem
7.1 when both variables have the same within-zone variance share \(a\).
By Proposition 7.4 an entry is attained if and only if \(\rho\le1-2a\),
i.e.~if it does not exceed \(1-\rho\); the other entries are strict
upper bounds}\tabularnewline
\toprule\noalign{}
\begin{minipage}[b]{\linewidth}\centering
Within-zone variance share \(a=b\)
\end{minipage} & \begin{minipage}[b]{\linewidth}\centering
\(\rho=0\)
\end{minipage} & \begin{minipage}[b]{\linewidth}\centering
\(\rho=0.3\)
\end{minipage} & \begin{minipage}[b]{\linewidth}\centering
\(\rho=0.6\)
\end{minipage} & \begin{minipage}[b]{\linewidth}\centering
\(\rho=0.9\)
\end{minipage} \\
\midrule\noalign{}
\endfirsthead
\toprule\noalign{}
\begin{minipage}[b]{\linewidth}\centering
Within-zone variance share \(a=b\)
\end{minipage} & \begin{minipage}[b]{\linewidth}\centering
\(\rho=0\)
\end{minipage} & \begin{minipage}[b]{\linewidth}\centering
\(\rho=0.3\)
\end{minipage} & \begin{minipage}[b]{\linewidth}\centering
\(\rho=0.6\)
\end{minipage} & \begin{minipage}[b]{\linewidth}\centering
\(\rho=0.9\)
\end{minipage} \\
\midrule\noalign{}
\endhead
\bottomrule\noalign{}
\endlastfoot
0.02 & 0.020 & 0.027 & 0.033 & 0.039 \\
0.05 & 0.053 & 0.068 & 0.084 & 0.100 \\
0.10 & 0.111 & 0.144 & 0.178 & 0.211 \\
0.20 & 0.250 & 0.325 & 0.400 & 0.475 \\
0.30 & 0.429 & 0.557 & 0.686 & 0.814 \\
0.50 & 1.000 & 1.300 & 1.600 & 1.900 \\
\end{longtable}

Theorem 7.1 needs only \(a\), \(b\) and \(\rho\). With fine-scale data,
\(a=\sum_jp_j\operatorname{Var}(X\mid A_j)/\sigma_X^2\) can be computed
directly, so for any given zoning one can compute in advance how far the
zonal correlation can deviate from the fine-scale correlation, without
knowing how fast the variables change in space.

\hypertarget{lipschitz-fields-and-the-second-order-rate}{%
\subsection{7.2 Lipschitz fields and the second-order
rate}\label{lipschitz-fields-and-the-second-order-rate}}

\textbf{Corollary 7.2 (Lipschitz fields and zone diameter).} Let
\(D\subset\mathbb{R}^q\), let \(\mu\) be a probability measure on \(D\),
and let \(X,Y:D\to\mathbb{R}\) be \(L_X\)- and \(L_Y\)-Lipschitz
(Euclidean distance) with \(\sigma_X,\sigma_Y>0\). If every zone has
diameter at most \(\delta\), put \[
\tau_X=\frac{L_X\delta}{2\sigma_X},\qquad \tau_Y=\frac{L_Y\delta}{2\sigma_Y},\qquad \tau=\max(\tau_X,\tau_Y),
\] and suppose \(\tau<1\). Then \[
|\rho_B-\rho|\le\frac{|\rho|\Big(1-\sqrt{(1-\tau_X^2)(1-\tau_Y^2)}\Big)+\tau_X\tau_Y}{\sqrt{(1-\tau_X^2)(1-\tau_Y^2)}}=:\mathcal{B}_{\mathrm{L}}(\tau_X,\tau_Y;\rho).
\] If moreover \(\tau\le1/\sqrt2\), then \[
0\le\mathcal{B}_{\mathrm{L}}-\Big[\tau_X\tau_Y+\frac{|\rho|}{2}\big(\tau_X^2+\tau_Y^2\big)\Big]\le5\,\tau^4 .\tag{7.1}
\]

\emph{Proof.} Any two points of a zone \(A_j\) are at distance at most
\(\delta\), so the values of \(X\) on \(A_j\) lie in an interval of
length at most \(L_X\delta\). Popoviciu's inequality (Popoviciu 1935): a
random variable \(U\) with values in an interval \([c-\ell/2,c+\ell/2]\)
of length \(\ell\) satisfies
\(\operatorname{Var}U\le\mathbb{E}(U-c)^2\le(\ell/2)^2\). Hence
\(\operatorname{Var}(X\mid A_j)\le(L_X\delta)^2/4\),
\(W_X=\sum_jp_j\operatorname{Var}(X\mid A_j)\le(L_X\delta)^2/4\),
i.e.~\(a\le\tau_X^2<1\); likewise \(b\le\tau_Y^2<1\). By the
monotonicity in Theorem 7.1,
\(|\rho_B-\rho|\le\mathcal{B}(\tau_X^2,\tau_Y^2;\rho)=\mathcal{B}_{\mathrm{L}}\).

For (7.1), write \(u=\tau_X^2\), \(w=\tau_Y^2\) and
\(f(x)=(1-x)^{-1/2}=\sum_{n\ge0}c_nx^n\) with
\(c_n=\binom{2n}{n}4^{-n}\) (\(c_0=1\), \(c_1=1/2\), \(c_2=3/8\), and
\(c_n\) decreasing). For \(0\le x\le1/2\), \[
0\le f(x)-1-\frac x2=\sum_{n\ge2}c_nx^n\le x^2\sum_{n\ge2}c_2\,2^{-(n-2)}=\frac34x^2 .
\] Hence \(1/\sqrt{(1-u)(1-w)}=f(u)f(w)\) satisfies \[
0\le f(u)f(w)-1-\frac{u+w}{2}\le\frac{uw}{4}+\frac34(u^2+w^2)+\frac38(u^2w+uw^2)+\frac{9}{16}u^2w^2\le2.27\,\tau^4,
\] using \(u,w\le\tau^2\le1/2\). So \(|\rho|\big(f(u)f(w)-1\big)\)
differs from the leading term \(|\rho|(u+w)/2\) by an amount in
\([0,2.27\tau^4]\), while
\(\tau_X\tau_Yf(u)f(w)-\tau_X\tau_Y=\tau_X\tau_Y\big(f(u)f(w)-1\big)\in[0,\ \tau^2(\tau^2+2.27\tau^4)]\subseteq[0,\ 2.14\tau^4]\).
Adding and using \(|\rho|\le1\) gives (7.1). \(\blacksquare\)

\(\ell_X=\sigma_X/L_X\) is the \textbf{variation length} of \(X\): the
distance over which \(X\) changes by one standard deviation at its
maximal rate; \(\tau_X=\delta/(2\ell_X)\).

\textbf{Corollary 7.3 (second-order rate).} Under the assumptions of
Corollary 7.2, fix \(L_X,L_Y,\sigma_X,\sigma_Y\) and \(\rho\), and
regard the error bound as a function of the largest zone diameter
\(\delta\): \[
\begin{aligned}\mathcal{B}_{\mathrm{L}}(\delta)&=\kappa\,\delta^2+r(\delta),\qquad 0\le r(\delta)\le5\tau(\delta)^4\quad(\tau(\delta)\le1/\sqrt2),\\\kappa&=\frac{L_XL_Y}{4\sigma_X\sigma_Y}+\frac{|\rho|}{8}\Big(\frac{L_X^2}{\sigma_X^2}+\frac{L_Y^2}{\sigma_Y^2}\Big).\end{aligned}
\] Hence the error bound is of second order in \(\delta\): halving
\(\delta\) reduces the leading term \(\kappa\delta^2\) exactly to one
quarter, and
\(\mathcal{B}_{\mathrm{L}}(\delta/2)/\mathcal{B}_{\mathrm{L}}(\delta)\to1/4\)
as \(\delta\to0\) (\(\kappa>0\), because \(\sigma_X,\sigma_Y>0\) forces
\(L_X,L_Y>0\)).

\emph{Proof.} Substituting \(\tau_X=L_X\delta/(2\sigma_X)\) and
\(\tau_Y=L_Y\delta/(2\sigma_Y)\), the leading term
\(\tau_X\tau_Y+\tfrac{|\rho|}{2}(\tau_X^2+\tau_Y^2)\) equals
\(\kappa\delta^2\); the remainder is controlled by (7.1), and
\(\tau(\delta)^4=O(\delta^4)\). Hence
\(\mathcal{B}_{\mathrm{L}}(\delta/2)/\mathcal{B}_{\mathrm{L}}(\delta)=\big(\kappa\delta^2/4+O(\delta^4)\big)/\big(\kappa\delta^2+O(\delta^4)\big)\to1/4\).
\(\blacksquare\)

The second order arises because within-zone differences grow at most
linearly with the diameter, while aggregation loses within-zone
variances and covariances, which are quadratic in those differences.
Corollary 7.2 relies on Lipschitz constants, which are set by the
steepest part of the field and are often very conservative (in Figure
4(e) the bound of Corollary 7.2 exceeds that of Theorem 7.1 by more than
an order of magnitude); in practice one should compute \(a,b\) directly
from fine-scale data and use Theorem 7.1. The role of Corollary 7.2 is
to explain why the error is of second order. For stationary random
fields and data-independent designs, Proposition 8.16 gives the
design-level counterpart for the within-zone variance shares and hence
for the bound: they vanish at second order in the zone scale for
mean-square differentiable fields, but only at order \(\alpha<2\) for
rough fields with \(C(0)-C(h)\sim c|h|^\alpha\) (first order for the
exponential covariance), and the deviation itself can then be of first
order; the rate is set by the smoothness of the fields.

\hypertarget{sharpness}{%
\subsection{7.3 Sharpness}\label{sharpness}}

\textbf{Proposition 7.4 (when the bound of Theorem 7.1 is attained).}
For any \(a,b\in(0,1)\) and \(\rho_B^*\in[-1,1]\) there are an atomless
support measure space and a 3-zone zoning with within-zone variance
shares \(a,b\), zonal correlation \(\rho_B^*\), pooled within-zone
correlation \(\rho_W=-1\) and fine-scale correlation
\(\rho^*=s\rho_B^*-\sqrt{ab}\). When \(\rho^*\ge0\), the inequality of
Theorem 7.1 holds with equality. Consequently, for \(a,b\in(0,1)\) and
\(\rho\in[-1,1]\) the bound \(\mathcal{B}(a,b;\rho)\) is attained,
i.e.~\(|\rho_B-\rho|=\mathcal{B}(a,b;\rho)\) for some support measure
space, zoning and variables with within-zone variance shares \(a,b\) and
fine-scale correlation \(\rho\), if and only if
\(|\rho|+\sqrt{ab}\le s\), equivalently
\(\mathcal{B}(a,b;\rho)\le1-|\rho|\); otherwise
\(|\rho_B-\rho|<\mathcal{B}(a,b;\rho)\) in every such configuration.

\emph{Proof.} Let \(\Omega=\{1,2,3\}\times[0,1]\), let \(\mu\) be the
product of the uniform distribution on \(\{1,2,3\}\) and Lebesgue
measure on \([0,1]\), and let the zones be \(A_j=\{j\}\times[0,1]\).
\(\mu\) is atomless: if \(\mu(A)>0\), some section
\(A^{(j)}=\{t:(j,t)\in A\}\) has positive Lebesgue measure, and
\(t\mapsto\lambda(A^{(j)}\cap[0,t])\) is continuous, so \(A\) has
subsets of measure strictly between \(0\) and \(\mu(A)\). With
\(p=(1/3,1/3,1/3)\), the linear-algebra construction in the sufficiency
part of the proof of Corollary 6.3 (which does not use (H)) gives
\(V\in\mathbb{R}^{3\times2}\) with \(p^\top V=0\) and
\(V^\top D_pV=\begin{pmatrix}1&\rho_B^*\\\rho_B^*&1\end{pmatrix}\);
putting \(u_{1j}=\sqrt{1-a}\,V_{j1}\) and \(u_{2j}=\sqrt{1-b}\,V_{j2}\),
the three points \((u_{1j},u_{2j})\) have equally weighted mean \(0\),
variances \(1-a\) and \(1-b\), and correlation \(\rho_B^*\). Let
\(w(t)=\mathbf{1}\{t<1/2\}-\mathbf{1}\{t\ge1/2\}\) (mean \(0\) and
variance \(1\) on \([0,1]\)) and define \[
X(j,t)=u_{1j}+\sqrt{a}\,w(t),\qquad Y(j,t)=u_{2j}-\sqrt{b}\,w(t).
\] Then \(\mathbb{E}[X\mid A_j]=u_{1j}\),
\(\mathbb{E}[Y\mid A_j]=u_{2j}\), \(\operatorname{Var}(X\mid A_j)=a\),
\(\operatorname{Var}(Y\mid A_j)=b\), and the within-zone covariance is
\(-\sqrt{ab}\), so \(\rho_W=-1\); \(B_X=1-a\), \(W_X=a\),
\(\sigma_X^2=1\), so the within-zone variance share of \(X\) is \(a\),
and that of \(Y\) is \(b\); by Lemma 5.1,
\(\operatorname{Cov}(X,Y)=C_B+C_W=\rho_B^*\sqrt{(1-a)(1-b)}-\sqrt{ab}=\rho^*\).
If \(\rho^*\ge0\), then \(s\rho_B^*\ge\sqrt{ab}\), so \(\rho_B^*\ge0\),
and \[
\rho_B-\rho=\rho_B^*-s\rho_B^*+\sqrt{ab}=\rho_B^*(1-s)+\sqrt{ab}\ge0,
\] while
\(\mathcal{B}=\big(\rho^*(1-s)+\sqrt{ab}\big)/s=\big(s\rho_B^*(1-s)-\sqrt{ab}(1-s)+\sqrt{ab}\big)/s=\rho_B^*(1-s)+\sqrt{ab}\);
the two coincide.

\emph{The last statement.} Replacing \(Y\) by \(-Y\) changes the signs
of \(\rho\), \(\rho_B\) and \(\rho_W\) and leaves \(a\), \(b\) and
\(|\rho_B-\rho|\) unchanged, so we may assume \(\rho\ge0\). By the proof
of Theorem 7.1, \(\rho_B-\rho=\rho(1/s-1)-\sqrt{ab}\,\rho_W/s\), where
\(1/s-1>0\), \(\sqrt{ab}>0\) and \(|\rho_W|\le1\). Hence
\(|\rho_B-\rho|=\mathcal{B}\) forces \(\rho_W=-1\) if \(\rho>0\), and
\(|\rho_W|=1\) if \(\rho=0\); in both cases
\(|\rho_B|=(\rho+\sqrt{ab})/s\), and \(|\rho_B|\le1\) gives
\(\rho+\sqrt{ab}\le s\). Conversely, if \(\rho+\sqrt{ab}\le s\), put
\(\rho_B^*=(\rho+\sqrt{ab})/s\in(0,1]\); the construction above has
fine-scale correlation \(s\rho_B^*-\sqrt{ab}=\rho\ge0\) and attains the
bound. Finally, \(\mathcal{B}\le1-\rho\) is equivalent to
\(\rho(1-s)+\sqrt{ab}\le s-s\rho\), i.e.~to \(\rho+\sqrt{ab}\le s\).
\(\blacksquare\)

\textbf{Example 7.5 (the second order cannot be improved).} Take
\(D=[0,1)\), \(\mu\) Lebesgue measure, an integer \(n\ge2\),
\(\delta=1/n\), and the zoning into \(n\) intervals
\(A_j=[(j-1)\delta,j\delta)\) of length \(\delta\) with midpoints
\(c_j\). Let \(X(x)=x\) and \(Y(x)=x+h(x)\), where
\(h(x)=\delta/2-|x-c_j|\) for \(x\in A_j\) is a triangle wave of period
\(\delta\). \(X\) is 1-Lipschitz and \(Y\) is 2-Lipschitz. Within each
zone \(h\) is symmetric about \(c_j\), with zone mean \(\delta/4\) and
within-zone variance \(\delta^2/48\), and it is uncorrelated with
\(x-c_j\); so \(\bar y_j=\bar x_j+\delta/4\), the zone means are
collinear and \(\rho_B=1\); and since \(\operatorname{Cov}(X,h)=0\), \[
\rho=\frac{\operatorname{Var}X}{\sqrt{\operatorname{Var}X\,(\operatorname{Var}X+\operatorname{Var}h)}}=\Big(1+\frac{\delta^2/48}{1/12}\Big)^{-1/2}=\Big(1+\frac{\delta^2}{4}\Big)^{-1/2},
\] so \(\rho_B-\rho=1-(1+\delta^2/4)^{-1/2}=\delta^2/8+O(\delta^4)\).
The fields depend on \(n\), but \(L_X=1\), \(L_Y=2\),
\(\sigma_X^2=1/12\) and \(\sigma_Y^2=1/12+\delta^2/48\) are bounded
above and away from zero uniformly in \(n\). Thus for Lipschitz fields
and zones of diameter \(\delta\) the deviation can be of exact order
\(\delta^2\), and no bound that depends only on these constants,
\(\rho\) and \(\delta\) can be of smaller order than \(\delta^2\): the
order in Corollary 7.3 cannot be improved.

Figure 4(d)(e) illustrates this section: for the same pair of variables,
the zonal correlation over an \(8\times8\) grid of square blocks (side
64 pixels) is \(+0.896\), almost equal to the fine-scale correlation of
\(0.878\); as the blocks grow from 2 to 128 pixels, the actual deviation
always stays below the bound of Theorem 7.1, and both grow roughly as
the square of the zone scale.

\hypertarget{data-independent-zoning-designs-spectral-theory}{%
\section{8 Data-independent zoning designs: spectral
theory}\label{data-independent-zoning-designs-spectral-theory}}

The zonings of Section 6 may be drawn by looking at the data. This
section considers the other kind of zoning: the shapes and positions of
the zones are fixed in advance, independently of the variables under
study, as with regular grids, fixed administrative partitions, or zones
of given shapes laid out at random. The zonal second moments are then
determined by the spectral matrix of the two fields, and the conclusions
are the opposite of those of Section 6:

\begin{itemize}
\tightlist
\item
  the zonal correlation is a weighted average of the real coherence
  times a factor at most 1, so it always lies between \(0\) and the
  extreme values of the coherence, i.e.~in
  \([\min(0,\operatorname{ess\,inf}\operatorname{coh}),\)
  \(\max(0,\operatorname{ess\,sup}\operatorname{coh})]\) (Theorem 8.9);
\item
  the zonal correlation takes the same non-zero value for all connected
  zone shapes and their finite mixtures if and only if the two fields
  follow the intrinsic correlation model, and it vanishes for all of
  them if and only if the real part of the cross-spectrum vanishes
  (Theorem 8.11); when this condition fails and the fine-scale
  correlation is non-zero, some design made of at most two connected
  zone shapes has a zonal correlation different from the fine-scale
  correlation (Corollary 8.12); and equal correlations at any finite
  number of designs do not imply the condition (Proposition 8.13);
\item
  the aggregated regression slope of a linear model is unbiased for all
  connected designs if and only if the symmetrized cross-covariance
  between the explanatory variable and the error vanishes identically
  (Theorem 8.14);
\item
  as the zones shrink, the zonal correlation tends to the fine-scale
  correlation, and the within-zone variance shares that bound the
  deviation vanish at a rate set by the smoothness of the fields; as the
  zones grow without bound, the zonal correlation tends to the
  zero-frequency coherence, whatever the zone shape, for
  square-integrable kernels and spectral densities that are bounded and
  continuous at the origin with \(f_{XX}(0)f_{YY}(0)>0\) (Theorem 8.15,
  Proposition 8.16).
\end{itemize}

Section 8.1 derives the spectral representation of jointly stationary
fields from scratch (uniqueness theorem, Bochner's theorem, Cramér's
theorem); Section 8.2 establishes the second moments of block means and
the design-level covariance; Sections 8.3--8.6 prove the results above.

\hypertarget{spectral-representation-of-jointly-stationary-fields}{%
\subsection{8.1 Spectral representation of jointly stationary
fields}\label{spectral-representation-of-jointly-stationary-fields}}

\textbf{Setting.} The index space is \(\mathbb{R}^q\) (\(q\ge1\)), with
Lebesgue measure \(ds\); the volume of a Borel set \(A\) is \(|A|\).
\((X,Y)\) is a pair of real-valued random fields indexed by
\(\mathbb{R}^q\) such that:

\begin{itemize}
\tightlist
\item
  \textbf{joint measurability}:
  \((s,\text{sample point})\mapsto(X(s),Y(s))\) is measurable with
  respect to the product of \(\mathcal{B}(\mathbb{R}^q)\) and the
  \(\sigma\)-algebra of the sample space;
\item
  \textbf{finite second moments and zero means}:
  \(\mathbb{E}X(s)^2<\infty\), \(\mathbb{E}Y(s)^2<\infty\),
  \(\mathbb{E}X(s)=\mathbb{E}Y(s)=0\);
\item
  \textbf{joint weak stationarity}: for \(U,V\in\{X,Y\}\),
  \(\mathbb{E}[U(t+h)V(t)]=C_{UV}(h)\) does not depend on \(t\);
\item
  \textbf{continuous covariance functions}: \(C_{XX}\), \(C_{YY}\),
  \(C_{XY}\) are continuous on \(\mathbb{R}^q\).
\end{itemize}

When the means are not zero, subtracting the constant means changes no
block covariance, so zero means are no loss of generality. The
fine-scale correlation is \(\rho=C_{XY}(0)/\sqrt{C_{XX}(0)C_{YY}(0)}\)
(when \(C_{XX}(0)C_{YY}(0)>0\)).

For a finite complex Borel measure \(\nu\) on \(\mathbb{R}^q\), its
\textbf{Fourier--Stieltjes transform} is
\(\hat\nu(h)=\int e^{ih\cdot\omega}\,\nu(d\omega)\); by dominated
convergence (S2), \(\hat\nu\) is continuous and
\(|\hat\nu|\le|\nu|(\mathbb{R}^q)\), where \(|\nu|\) is the total
variation measure of \(\nu\). A function
\(\varphi:\mathbb{R}^q\to\mathbb{C}\) is \textbf{positive definite} if
for every \(n\ge1\), \(s_1,\dots,s_n\in\mathbb{R}^q\) and
\(c\in\mathbb{C}^n\), \(\sum_{j,k}c_j\bar c_k\,\varphi(s_j-s_k)\) is a
non-negative real number.

\textbf{Theorem 8.1 (Fourier--Stieltjes uniqueness).} Let \(\nu\) be a
finite complex Borel measure on \(\mathbb{R}^q\) with \(\hat\nu(h)=0\)
for all \(h\). Then \(\nu=0\). Consequently two finite complex measures
with the same Fourier--Stieltjes transform are equal.

\emph{Proof.} Let \(\psi\in\mathcal{S}(\mathbb{R}^q)\) and put
\(\phi(\omega)=\int\psi(h)e^{ih\cdot\omega}\,dh\). Since
\(\int\!\!\int|\psi(h)|\,dh\,|\nu|(d\omega)=\|\psi\|_1|\nu|(\mathbb{R}^q)<\infty\),
Fubini's theorem (S1) gives \[
\int\phi\,d\nu=\int\psi(h)\Big(\int e^{ih\cdot\omega}\,\nu(d\omega)\Big)dh=\int\psi(h)\,\hat\nu(h)\,dh=0 .
\] \(\phi(\omega)\) is the Fourier transform of \(\psi\) evaluated at
\(-\omega\); by (S6), as \(\psi\) ranges over \(\mathcal{S}\), so does
\(\phi\). Hence \(\int\phi\,d\nu=0\) for all
\(\phi\in\mathcal{S}\supseteq C_c^\infty\). For
\(\phi\in C_0(\mathbb{R}^q)\), choose by (S5) \(\phi_n\in C_c^\infty\)
with \(\|\phi_n-\phi\|_\infty\to0\); then
\(|\int\phi\,d\nu|=|\int(\phi-\phi_n)\,d\nu|\le\|\phi-\phi_n\|_\infty|\nu|(\mathbb{R}^q)\to0\).
So \(\nu\) and the zero measure define the same linear functional on
\(C_0(\mathbb{R}^q)\). Finite Borel measures on \(\mathbb{R}^q\) are
regular, and by uniqueness of the representing measure in the Riesz
representation theorem (S5), \(\nu=0\). The second sentence follows by
applying the first to the difference of the two measures.
\(\blacksquare\)

\textbf{Theorem 8.2 (Bochner).} Let
\(\varphi:\mathbb{R}^q\to\mathbb{C}\) be continuous and positive
definite. Then there is a unique finite positive Borel measure \(\mu\)
with \[
\varphi(h)=\int e^{ih\cdot\omega}\,\mu(d\omega)\qquad(h\in\mathbb{R}^q),
\] and \(\mu(\mathbb{R}^q)=\varphi(0)\). Conversely, the
Fourier--Stieltjes transform of a finite positive measure is continuous
and positive definite.

\emph{Proof.}

\emph{Step 0: the converse.} If \(\varphi=\hat\mu\) with \(\mu\) finite
and positive, then \(\varphi\) is continuous and \[
\sum_{j,k}c_j\bar c_k\varphi(s_j-s_k)=\int\sum_{j,k}c_je^{is_j\cdot\omega}\,\overline{c_ke^{is_k\cdot\omega}}\,\mu(d\omega)=\int\Big|\sum_jc_je^{is_j\cdot\omega}\Big|^2\mu(d\omega)\ge0 .
\]

\emph{Step 1: basic properties.} Taking \(n=1\) gives
\(\varphi(0)\ge0\). Taking \(n=2\), \(s_1=h\), \(s_2=0\): for all
\(c\in\mathbb{C}^2\), \[
(|c_1|^2+|c_2|^2)\varphi(0)+c_1\bar c_2\,\varphi(h)+c_2\bar c_1\,\varphi(-h)\ \text{is a non-negative real number}.
\] \(c=(1,1)\) gives \(\varphi(h)+\varphi(-h)\in\mathbb{R}\), and
\(c=(1,i)\) gives \(-i\big(\varphi(h)-\varphi(-h)\big)\in\mathbb{R}\);
together, \(\varphi(-h)=\overline{\varphi(h)}\). So the matrix
\(\begin{pmatrix}\varphi(0)&\varphi(h)\\\overline{\varphi(h)}&\varphi(0)\end{pmatrix}\)
is Hermitian positive semi-definite, and its determinant is
non-negative: \(|\varphi(h)|\le\varphi(0)\). In particular \(\varphi\)
is bounded.

\emph{Step 2: positive definiteness in integral form.} For
\(\psi\in L^1(\mathbb{R}^q)\) let \[
J(\psi)=\iint\varphi(s-t)\,\psi(s)\,\overline{\psi(t)}\,ds\,dt ;
\] since \(|\varphi|\le\varphi(0)\), the integrand is absolutely
integrable. First let \(\psi\in C_c(\mathbb{R}^q)\). The integrand is
continuous with compact support in \(\mathbb{R}^{2q}\), so \(J(\psi)\)
is a limit of Riemann sums: dividing \(\mathbb{R}^q\) into cubes of side
\(1/n\) with centres \(s_1,s_2,\dots\) (only finitely many of which meet
the support), \[
J(\psi)=\lim_{n\to\infty}\sum_{j,k}\varphi(s_j-s_k)\,\frac{\psi(s_j)}{n^q}\,\overline{\frac{\psi(s_k)}{n^q}}\ \ge0,
\] each sum being a quadratic form as in the definition of positive
definiteness. For general \(\psi\in L^1\) choose by (S5)
\(\psi_n\in C_c\) with \(\|\psi_n-\psi\|_1\to0\). From \[
\psi(s)\overline{\psi(t)}-\psi_n(s)\overline{\psi_n(t)}=\big(\psi(s)-\psi_n(s)\big)\overline{\psi(t)}+\psi_n(s)\big(\overline{\psi(t)}-\overline{\psi_n(t)}\big)
\] we get
\(|J(\psi)-J(\psi_n)|\le\varphi(0)\big(\|\psi\|_1+\|\psi_n\|_1\big)\|\psi-\psi_n\|_1\to0\),
so \(J(\psi)\ge0\).

\emph{Step 3: Gaussian regularization.} For \(\varepsilon>0\) define \[
f_\varepsilon(\omega)=(2\pi)^{-q}\int\varphi(h)\,e^{-\varepsilon|h|^2/2}\,e^{-ih\cdot\omega}\,dh .
\] Fix \(\omega\) and take
\(\psi(s)=e^{-\varepsilon|s|^2}e^{-is\cdot\omega}\in L^1\) in Step 2, so
that
\(\psi(s)\overline{\psi(t)}=e^{-\varepsilon(|s|^2+|t|^2)}e^{-i(s-t)\cdot\omega}\).
Change variables to \(h=s-t\), \(u=s+t\): the absolute value of the
Jacobian determinant is \(2^q\), so \(ds\,dt=2^{-q}\,dh\,du\); moreover
\(|s|^2+|t|^2=\tfrac14\big(|u+h|^2+|u-h|^2\big)=\tfrac12\big(|h|^2+|u|^2\big)\).
By Fubini's theorem, \[
0\le J(\psi)=2^{-q}\int\varphi(h)e^{-\varepsilon|h|^2/2}e^{-ih\cdot\omega}\,dh\cdot\int e^{-\varepsilon|u|^2/2}\,du=2^{-q}\Big(\frac{2\pi}{\varepsilon}\Big)^{q/2}(2\pi)^q\,f_\varepsilon(\omega).
\] Hence \(f_\varepsilon(\omega)\) is a non-negative real number.

\emph{Step 4: \(\int f_\varepsilon=\varphi(0)\).} For \(\delta>0\), the
Gaussian integral
\(\int e^{-ih\cdot\omega}e^{-\delta|\omega|^2/2}\,d\omega=(2\pi/\delta)^{q/2}e^{-|h|^2/(2\delta)}\)
and Fubini's theorem give \[
\int f_\varepsilon(\omega)\,e^{-\delta|\omega|^2/2}\,d\omega=\int\varphi(h)\,e^{-\varepsilon|h|^2/2}\,p_\delta(h)\,dh,\qquad p_\delta(h)=(2\pi\delta)^{-q/2}e^{-|h|^2/(2\delta)} .
\] \(p_\delta\) is the density of \(N(0,\delta I_q)\), so the right-hand
side equals
\(\mathbb{E}\big[\varphi(\sqrt\delta\,G)\,e^{-\varepsilon\delta|G|^2/2}\big]\)
with \(G\sim N(0,I_q)\); as \(\delta\downarrow0\) it tends to
\(\varphi(0)\) by continuity of \(\varphi\) at the origin,
\(|\varphi|\le\varphi(0)\) and dominated convergence. The integrand on
the left is non-negative and increases to \(f_\varepsilon\) as
\(\delta\downarrow0\), so by monotone convergence the left-hand side
tends to \(\int f_\varepsilon\). Hence
\(\int f_\varepsilon=\varphi(0)<\infty\) and \(f_\varepsilon\in L^1\).

\emph{Step 5: inversion.}
\(\varphi_\varepsilon(h)=\varphi(h)e^{-\varepsilon|h|^2/2}\) is
continuous and integrable, and so is its Fourier transform
\(\int\varphi_\varepsilon(h)e^{-ih\cdot\omega}\,dh=(2\pi)^qf_\varepsilon(\omega)\).
By the inversion formula (S6), for every \(h\), \[
\varphi_\varepsilon(h)=(2\pi)^{-q}\int(2\pi)^qf_\varepsilon(\omega)\,e^{ih\cdot\omega}\,d\omega=\int e^{ih\cdot\omega}f_\varepsilon(\omega)\,d\omega .
\] So \(\varphi_\varepsilon\) is the Fourier--Stieltjes transform of the
finite positive measure
\(\mu_\varepsilon(d\omega)=f_\varepsilon(\omega)\,d\omega\), with
\(\mu_\varepsilon(\mathbb{R}^q)=\varphi(0)\).

\emph{Step 6: passage to the limit.} If \(\varphi(0)=0\), then
\(\varphi\equiv0\) by Step 1 and \(\mu=0\) works. Otherwise
\(P_n=\mu_{1/n}/\varphi(0)\) are probability measures whose
characteristic functions \(\varphi_{1/n}/\varphi(0)\) converge pointwise
to \(\varphi/\varphi(0)\), which is continuous at the origin. By Lévy's
continuity theorem (S8), \(P_n\) converges weakly to a probability
measure \(P\) with characteristic function \(\varphi/\varphi(0)\). Put
\(\mu=\varphi(0)P\); then \(\hat\mu=\varphi\) and
\(\mu(\mathbb{R}^q)=\hat\mu(0)=\varphi(0)\). Uniqueness follows from
Theorem 8.1. \(\blacksquare\)

\textbf{Theorem 8.3 (Cramér 1940: matrix spectral representation).} Let
\((X,Y)\) satisfy the setting of this section. There is a unique family
\(F=(F_{UV})_{U,V\in\{X,Y\}}\) of finite complex Borel measures on
\(\mathbb{R}^q\) with \[
C_{UV}(h)=\int e^{ih\cdot\omega}\,F_{UV}(d\omega)\qquad(h\in\mathbb{R}^q,\ U,V\in\{X,Y\}),
\] and for every Borel set \(A\) the \(2\times2\) matrix
\(F(A)=\big(F_{UV}(A)\big)\) is Hermitian positive semi-definite. In
particular \(F_{XX}\) and \(F_{YY}\) are finite positive measures with
\(F_{XX}(\mathbb{R}^q)=C_{XX}(0)\) and
\(F_{YY}(\mathbb{R}^q)=C_{YY}(0)\).

\(F\) is called the \textbf{spectral matrix} of \((X,Y)\), and
\(F_{XY}\) the \textbf{cross-spectral measure}.

\emph{Proof.} For \(c=(c_X,c_Y)\in\mathbb{C}^2\) let
\(Z_c(s)=c_XX(s)+c_YY(s)\). Then \[
\mathbb{E}\big[Z_c(t+h)\overline{Z_c(t)}\big]=\sum_{U,V}c_U\bar c_V\,C_{UV}(h)=:C_c(h)
\] does not depend on \(t\) and is continuous in \(h\) (note that
\(C_{YX}(h)=C_{XY}(-h)\) by stationarity, so \(C_{YX}\) is continuous).
\(C_c\) is positive definite: by stationarity,
\(\mathbb{E}[Z_c(s_j)\overline{Z_c(s_k)}]=C_c(s_j-s_k)\), so \[
\sum_{j,k}a_j\bar a_k\,C_c(s_j-s_k)=\mathbb{E}\Big|\sum_ja_jZ_c(s_j)\Big|^2\ge0 .
\] By Theorem 8.2 there is a unique finite positive measure \(\mu_c\)
with \(C_c=\hat\mu_c\).

\emph{Polarization.} Let \(F_{UU}=\mu_{e_U}\) and, for \(U\neq V\), \[
F_{UV}=\frac14\sum_{k=0}^{3}i^k\,\mu_{e_U+i^ke_V}.
\] For \(c=e_U+i^ke_V\) we have
\(C_c=C_{UU}+C_{VV}+i^{-k}C_{UV}+i^{k}C_{VU}\). Since
\(\sum_{k=0}^3i^k=\sum_{k=0}^3i^{2k}=0\), \[
\frac14\sum_{k=0}^{3}i^k\,C_{e_U+i^ke_V}=\frac14\Big[\sum_ki^k\,(C_{UU}+C_{VV})+4\,C_{UV}+\sum_ki^{2k}\,C_{VU}\Big]=C_{UV}.
\] By linearity of the transform, \(\hat F_{UV}=C_{UV}\); for \(U=V\)
this is the definition of \(\mu_{e_U}\).

\emph{Positive semi-definiteness.} For \(c\in\mathbb{C}^2\), the finite
complex measure \(\nu_c=\sum_{U,V}c_U\bar c_VF_{UV}\) has transform
\(\sum_{U,V}c_U\bar c_VC_{UV}=C_c=\hat\mu_c\), so \(\nu_c=\mu_c\ge0\) by
Theorem 8.1. Hence for every Borel set \(A\) and \(d=\bar c\), \[
d^*F(A)\,d=\sum_{U,V}\bar d_U\,F_{UV}(A)\,d_V=\nu_c(A)=\mu_c(A)\ge0 .
\] As \(d\) ranges over \(\mathbb{C}^2\), \(F(A)\) is Hermitian positive
semi-definite (a complex matrix \(M\) with \(d^*Md\) real for all \(d\)
is Hermitian).

\emph{Uniqueness.} If \(F'\) also represents the covariances, each
\(F_{UV}-F'_{UV}\) has zero transform and vanishes by Theorem 8.1. The
last sentence follows from
\(C_{UU}(0)=\hat F_{UU}(0)=F_{UU}(\mathbb{R}^q)\). \(\blacksquare\)

\textbf{Lemma 8.4 (symmetries of real fields).} Let \((X,Y)\) satisfy
the setting of this section with spectral matrix \(F\).

\begin{enumerate}
\def\labelenumi{(\alph{enumi})}
\item
  For \(U,V\in\{X,Y\}\) and Borel sets \(A\),
  \(F_{UV}(-A)=\overline{F_{UV}(A)}\) and
  \(F_{VU}(A)=\overline{F_{UV}(A)}\). In particular \(F_{XX}\) and
  \(F_{YY}\) are even measures, \(\Gamma:=\operatorname{Re}F_{XY}\) is a
  finite even signed measure, and \(\operatorname{Im}F_{XY}\) is an odd
  signed measure.
\item
  If \(g\) is a bounded measurable real-valued even function, then
  \(\int g\,d(\operatorname{Im}F_{XY})=0\), hence
  \(\int g\,dF_{XY}=\int g\,d\Gamma\in\mathbb{R}\).
\item
  For every Borel set \(A\), the real symmetric matrix
  \(\begin{pmatrix}F_{XX}(A)&\Gamma(A)\\\Gamma(A)&F_{YY}(A)\end{pmatrix}\)
  is positive semi-definite; in particular
  \(\Gamma(A)^2\le F_{XX}(A)\,F_{YY}(A)\).
\item
  \(C_{UU}(h)=\int\cos(h\cdot\omega)\,F_{UU}(d\omega)\), and the even
  and odd parts of the cross-covariance are \[
  \begin{aligned}
  \tfrac12\big(C_{XY}(h)+C_{XY}(-h)\big)&=\int\cos(h\cdot\omega)\,\Gamma(d\omega),\\
  \tfrac12\big(C_{XY}(h)-C_{XY}(-h)\big)&=-\int\sin(h\cdot\omega)\,(\operatorname{Im}F_{XY})(d\omega).
  \end{aligned}
  \] Hence \(\Gamma\) is uniquely determined by the even part of the
  cross-covariance, and \(\operatorname{Im}F_{XY}\) by its odd part.
\end{enumerate}

\emph{Proof.} (a) Let \(R(\omega)=-\omega\) and let \(R_\#\nu\) denote
the image of \(\nu\) under \(R\), \((R_\#\nu)(A)=\nu(-A)\). Since
\(C_{UV}\) is real-valued, \[
C_{UV}(h)=\overline{C_{UV}(h)}=\int e^{-ih\cdot\omega}\,\overline{F_{UV}}(d\omega)=\int e^{ih\cdot\omega}\,\big(R_\#\overline{F_{UV}}\big)(d\omega),
\] so by Theorem 8.1 \(F_{UV}=R_\#\overline{F_{UV}}\),
i.e.~\(F_{UV}(A)=\overline{F_{UV}(-A)}\). Next, with \(t'=t+h\),
stationarity gives
\(C_{VU}(h)=\mathbb{E}[V(t+h)U(t)]=\mathbb{E}[U(t'-h)V(t')]=C_{UV}(-h)\),
so \(\hat F_{VU}(h)=\hat F_{UV}(-h)=\widehat{R_\#F_{UV}}(h)\), and by
Theorem 8.1 \(F_{VU}=R_\#F_{UV}=\overline{F_{UV}}\). For \(U=V\),
\(F_{UU}\) is a positive measure and \(F_{UU}(-A)=F_{UU}(A)\). Splitting
\(F_{XY}(-A)=\overline{F_{XY}(A)}\) into real and imaginary parts gives
\(\Gamma(-A)=\Gamma(A)\) and
\(\operatorname{Im}F_{XY}(-A)=-\operatorname{Im}F_{XY}(A)\).

\begin{enumerate}
\def\labelenumi{(\alph{enumi})}
\setcounter{enumi}{1}
\item
  Let \(\nu=\operatorname{Im}F_{XY}\), so that \(R_\#\nu=-\nu\). By the
  change-of-variables formula for image measures and \(g\circ R=g\), \[
  \int g\,d\nu=\int g\circ R\,d\nu=\int g\,d(R_\#\nu)=-\int g\,d\nu ,
  \] so \(\int g\,d\nu=0\).
\item
  By (a), \(F_{YX}(A)=\overline{F_{XY}(A)}\), so
  \(\operatorname{Re}F(A)\) is real symmetric with off-diagonal entry
  \(\Gamma(A)\). For \(c\in\mathbb{R}^2\),
  \(c^\top\big(\operatorname{Re}F(A)\big)c=\operatorname{Re}\big(c^*F(A)c\big)\ge0\)
  (Theorem 8.3). A positive semi-definite real symmetric \(2\times2\)
  matrix has non-negative determinant.
\item
  Substitute
  \(e^{ih\cdot\omega}=\cos(h\cdot\omega)+i\sin(h\cdot\omega)\) and
  \(F_{XY}=\Gamma+i\operatorname{Im}F_{XY}\) into the representation. By
  (b), the integral of the even function \(\cos(h\cdot\omega)\) against
  the odd measure vanishes; likewise the integral of an odd function
  against an even measure vanishes (same argument as in (b)). Hence \[
  C_{XY}(h)=\int\cos(h\cdot\omega)\,\Gamma(d\omega)-\int\sin(h\cdot\omega)\,(\operatorname{Im}F_{XY})(d\omega),
  \] and adding and subtracting the same identity at \(-h\) gives the
  claim. The case of \(F_{UU}\) is identical. For the last sentence: if
  two even measures \(\nu_1,\nu_2\) satisfy
  \(\int\cos(h\cdot\omega)\,d(\nu_1-\nu_2)=0\) for all \(h\), then
  \(\nu_1-\nu_2\) has zero transform (the sine part integrates to zero
  against an even measure), and \(\nu_1=\nu_2\) by Theorem 8.1; the odd
  case is the same. \(\blacksquare\)
\end{enumerate}

When the spectral matrix has densities \(f_{XX},f_{YY},f_{XY}\) with
respect to Lebesgue measure, these are the usual spectral and
cross-spectral densities; \(\operatorname{Re}f_{XY}\) is the co-spectrum
and \(\operatorname{Im}f_{XY}\) corresponds to the quadrature spectrum
(Yaglom 1987). Lemma 8.4(d) shows that the even part of the
cross-covariance enters only \(\Gamma\), and the odd part only
\(\operatorname{Im}F_{XY}\).

\hypertarget{block-means-and-design-level-second-moments}{%
\subsection{8.2 Block means and design-level second
moments}\label{block-means-and-design-level-second-moments}}

\textbf{Aggregation kernels.} \(K\in L^1(\mathbb{R}^q)\) is an
\textbf{aggregation kernel} if \(K\ge0\) and \(\int K=1\). The uniform
density \(u_A=\mathbf{1}_A/|A|\) on a bounded Borel set \(A\) with
\(|A|>0\) is an aggregation kernel and corresponds to ``equally weighted
averaging over the zone \(A\)''; general aggregation kernels also allow
weighting within the zone. The \textbf{block mean} of the field over
\(K\) is \(X_K=\int K(s)X(s)\,ds\), and \(Y_K\) likewise. Write
\(\hat K(\omega)=\int e^{-is\cdot\omega}K(s)\,ds\). The translate of
\(K\) to position \(c\) is \(K(\cdot-c)\), and its rescaling by \(T>0\)
is \(K_T(s)=T^{-q}K(s/T)\).

\textbf{Lemma 8.5 (block means).} Let \((X,Y)\) satisfy the setting of
Section 8.1 and let \(K,L\) be aggregation kernels.

\begin{enumerate}
\def\labelenumi{(\alph{enumi})}
\item
  Almost surely \(\int K(s)|X(s)|\,ds<\infty\), so \(X_K\) is almost
  surely defined, and it is a random variable.
\item
  \(\mathbb{E}X_K=0\) and \(\mathbb{E}X_K^2\le C_{XX}(0)\).
\item
  \(\hat K\) is continuous, \(|\hat K|\le1\), \(\hat K(0)=1\) and
  \(\hat K(-\omega)=\overline{\hat K(\omega)}\), so \(|\hat K|^2\) is a
  continuous even function with values in \([0,1]\); moreover
  \(\sup_\omega|\hat K(\omega)-\hat L(\omega)|\le\|K-L\|_1\).
\item
  \(\widehat{K(\cdot-c)}(\omega)=e^{-ic\cdot\omega}\hat K(\omega)\) and
  \(\hat K_T(\omega)=\hat K(T\omega)\); in particular, translation does
  not change \(|\hat K|^2\).
\end{enumerate}

\emph{Proof.} (a) By the Cauchy--Schwarz inequality,
\(\mathbb{E}|X(s)|\le\big(\mathbb{E}X(s)^2\big)^{1/2}=C_{XX}(0)^{1/2}\).
By joint measurability and Tonelli's theorem (S1), \[
\mathbb{E}\int K(s)|X(s)|\,ds=\int K(s)\,\mathbb{E}|X(s)|\,ds\le C_{XX}(0)^{1/2}<\infty ,
\] so \(\int K|X|\,ds\) is finite almost surely. On this event of
probability one, applying Tonelli's theorem to the positive and negative
parts of \(KX\) shows that \(X_K\) is measurable.

\begin{enumerate}
\def\labelenumi{(\alph{enumi})}
\setcounter{enumi}{1}
\item
  By Fubini's theorem and the integrability in (a),
  \(\mathbb{E}X_K=\int K(s)\,\mathbb{E}X(s)\,ds=0\). Since \(K(s)\,ds\)
  is a probability measure, Jensen's inequality gives
  \(X_K^2\le\int K(s)X(s)^2\,ds\); taking expectations and using
  Tonelli's theorem,
  \(\mathbb{E}X_K^2\le\int K(s)\,C_{XX}(0)\,ds=C_{XX}(0)\).
\item
  \(|e^{-is\cdot\omega}K(s)|\le K(s)\), so \(\hat K\) is continuous by
  dominated convergence, and \(|\hat K(\omega)|\le\int K=1=\hat K(0)\).
  \(K\) is real-valued, so \(\hat K(-\omega)=\overline{\hat K(\omega)}\)
  and \(|\hat K(-\omega)|^2=|\hat K(\omega)|^2\). Finally
  \(|\hat K(\omega)-\hat L(\omega)|\le\int|K-L|\).
\item
  Substitute \(s\mapsto s+c\) and \(s=Tu\) respectively.
  \(\blacksquare\)
\end{enumerate}

\textbf{Lemma 8.6 (block covariance formula; change of support).} For
aggregation kernels \(K,L\), \[
\mathbb{E}\big[X_KY_L\big]=\iint K(s)\,L(t)\,C_{XY}(s-t)\,ds\,dt=\int\overline{\hat K(\omega)}\,\hat L(\omega)\,F_{XY}(d\omega).
\] In particular \[
\operatorname{Cov}(X_K,Y_K)=\int|\hat K|^2\,d\Gamma,\qquad \operatorname{Var}(X_K)=\int|\hat K|^2\,dF_{XX},\qquad \operatorname{Var}(Y_K)=\int|\hat K|^2\,dF_{YY},
\] and these three quantities do not change when \(K\) is translated to
\(K(\cdot-c)\).

\emph{Proof.} By the Cauchy--Schwarz inequality,
\(\mathbb{E}|X(s)Y(t)|\le\sqrt{C_{XX}(0)C_{YY}(0)}\), so \[
\mathbb{E}\iint K(s)L(t)\,|X(s)Y(t)|\,ds\,dt\le\sqrt{C_{XX}(0)C_{YY}(0)}<\infty ,
\] and Fubini's theorem gives \[
\mathbb{E}[X_KY_L]=\iint K(s)L(t)\,\mathbb{E}[X(s)Y(t)]\,ds\,dt=\iint K(s)L(t)\,C_{XY}(s-t)\,ds\,dt .
\] Substituting the representation of Theorem 8.3 and using Fubini's
theorem again (\(K,L\in L^1\), \(|F_{XY}|\) finite), \[
\begin{aligned}\iint K(s)L(t)\int e^{i(s-t)\cdot\omega}F_{XY}(d\omega)\,ds\,dt&=\int\Big(\int K(s)e^{is\cdot\omega}ds\Big)\Big(\int L(t)e^{-it\cdot\omega}dt\Big)F_{XY}(d\omega)\\&=\int\overline{\hat K}\,\hat L\,dF_{XY}.\end{aligned}
\] For \(K=L\), \(|\hat K|^2\) is a bounded even function, and by Lemma
8.4(b) \(\int|\hat K|^2dF_{XY}=\int|\hat K|^2d\Gamma\). The means are
zero, so the second moments are covariances. Translation invariance
follows from Lemma 8.5(d). \(\blacksquare\)

The first equality of Lemma 8.6 is the basic change-of-support
(regularization) formula of geostatistics: the covariance between blocks
is the double average of the point covariance (Journel \& Huijbregts
1978; Cressie 1993; Chilès \& Delfiner 2012; Gotway \& Young 2002). The
second shows that aggregation acts in the frequency domain as
multiplication by the filter \(|\hat K|^2\): large zones keep only the
low-frequency (long-wave) part, small zones keep almost all frequencies.

\textbf{Zoning designs.} Let \(\mathcal{K}_0\) be a finite set of
aggregation kernels. A probability measure \(\pi\) on \(\mathcal{K}_0\)
is a \textbf{zoning design}; \(\pi(K)\) is the share of the zone weights
carried by zones of shape \(K\). The \textbf{filter} of the design is \[
g_\pi=\sum_{K}\pi(K)\,|\hat K|^2 ,
\] a continuous even function with values in \([0,1]\), \(g_\pi(0)=1\),
linear in \(\pi\). For a finite measure \(G\) write \(G(g)=\int g\,dG\).
The \textbf{design-level zonal second-moment matrix} is \[
\Sigma(\pi)=\sum_K\pi(K)\begin{pmatrix}\operatorname{Var}X_K&\operatorname{Cov}(X_K,Y_K)\\\operatorname{Cov}(X_K,Y_K)&\operatorname{Var}Y_K\end{pmatrix}
=\begin{pmatrix}F_{XX}(g_\pi)&\Gamma(g_\pi)\\\Gamma(g_\pi)&F_{YY}(g_\pi)\end{pmatrix},
\] the second equality by Lemma 8.6. When
\(F_{XX}(g_\pi)F_{YY}(g_\pi)>0\), the \textbf{design-level zonal
correlation} is \[
\rho(\pi)=\frac{\Gamma(g_\pi)}{\sqrt{F_{XX}(g_\pi)\,F_{YY}(g_\pi)}} .
\] No aggregation corresponds to \(g\equiv1\): by Lemma 8.4(b),
\(C_{XY}(0)=F_{XY}(\mathbb{R}^q)=\Gamma(\mathbb{R}^q)\), so the
fine-scale correlation is
\(\rho=\Gamma(\mathbb{R}^q)/\sqrt{F_{XX}(\mathbb{R}^q)F_{YY}(\mathbb{R}^q)}\).
The zonal correlation depends only on the spectral matrix and not on any
particular realization of the fields; this is what ``data-independent''
means.

\textbf{Proposition 8.7 (meaning of the design-level second moments).}
Let \((K_i,c_i,w_i)_{i=1}^n\) be random elements independent of
\((X,Y)\) (precisely: the design and the fields are defined on a product
probability space, the design depending only on the first and the fields
only on the second coordinate), with \(K_i\) taking values in a finite
set \(\mathcal{K}_0\) of aggregation kernels, \(c_i\in\mathbb{R}^q\),
\(w_i\ge0\) and \(\sum_iw_i=1\). Let the block means of the \(i\)-th
zone be \(X_i=\int K_i(s-c_i)X(s)\,ds\), and \(Y_i\) likewise. Then \[
\mathbb{E}\Big[\sum_iw_iX_iY_i\Big]=\Sigma_{XY}(\pi),\qquad \mathbb{E}\Big[\sum_iw_iX_i^2\Big]=\Sigma_{XX}(\pi),\qquad \mathbb{E}\Big[\sum_iw_iY_i^2\Big]=\Sigma_{YY}(\pi),
\] where \(\pi(K)=\mathbb{E}\big[\sum_iw_i\mathbf{1}\{K_i=K\}\big]\).

\emph{Proof.} Fix \(K\in\mathcal{K}_0\). \((s,c)\mapsto K(s-c)\) is
Borel measurable, and by Tonelli's theorem
\(\mathbb{E}\int K(s-c)|X(s)|\,ds\le C_{XX}(0)^{1/2}\) for every \(c\),
so \((c,\text{sample point})\mapsto X_{K(\cdot-c)}\) is jointly
measurable and, for any distribution of \(c\), defined almost everywhere
with respect to the product of that distribution and the probability
measure of the fields. By the product structure in the statement,
Fubini's theorem applies on the product space. By the Cauchy--Schwarz
inequality and Lemma 8.5(b),
\(\mathbb{E}\big|X_{K(\cdot-c)}Y_{K(\cdot-c)}\big|\le\sqrt{C_{XX}(0)C_{YY}(0)}\),
and \(0\le w_i\le1\). We may therefore integrate first over the field
and then over the design (Fubini's theorem): \[
\mathbb{E}\big[w_iX_iY_i\big]=\mathbb{E}\big[w_i\,m(K_i,c_i)\big],\qquad m(K,c)=\mathbb{E}\big[X_{K(\cdot-c)}Y_{K(\cdot-c)}\big]=\operatorname{Cov}(X_K,Y_K),
\] the last equality being the translation invariance in Lemma 8.6.
Hence
\(\mathbb{E}[w_iX_iY_i]=\sum_K\operatorname{Cov}(X_K,Y_K)\,\mathbb{E}[w_i\mathbf{1}\{K_i=K\}]\);
summing over \(i\) gives the first identity, and the other two are the
same. \(\blacksquare\)

A fixed zoning is the degenerate special case. Proposition 8.7 shows
that \(\Sigma(\pi)\) is the expectation, under the law of the fields, of
the zonal second moments weighted by zone weights, provided the shapes,
positions and weights of the zones do not depend on the fields. The
zonal correlation computed from one realization on a bounded study
region (Section 5, centred at the weighted overall mean) fluctuates
around \(\rho(\pi)\): the larger the study region relative to the
correlation range of the fields, the closer the weighted overall mean is
to zero (when the spectral measures have no atom at the origin) and,
typically, the weighted second moments to their expectations. The error
bars of the Monte Carlo results in Figure 5(f) show this fluctuation; we
do not quantify it here (Section 13.3).

\hypertarget{coherence-and-the-coherence-bound}{%
\subsection{8.3 Coherence and the coherence
bound}\label{coherence-and-the-coherence-bound}}

\textbf{Lemma 8.8 (real coherence).} Let \(\Lambda=F_{XX}+F_{YY}\). Then
\(|\Gamma|\le\Lambda/2\) (\(|\Gamma|\) the total variation measure), so
\(F_{XX}\), \(F_{YY}\) and \(\Gamma\) are absolutely continuous with
respect to \(\Lambda\). Versions of the Radon--Nikodym densities \[
f_X=\frac{dF_{XX}}{d\Lambda},\qquad f_Y=\frac{dF_{YY}}{d\Lambda},\qquad \gamma=\frac{d\Gamma}{d\Lambda}
\] can be chosen with \(f_X,f_Y\ge0\), \(f_X+f_Y=1\) and
\(\gamma^2\le f_Xf_Y\) everywhere. Define the \textbf{real coherence} \[
\operatorname{coh}=\begin{cases}\gamma\big/\sqrt{f_Xf_Y},& f_Xf_Y>0,\\[2pt] 0,& f_Xf_Y=0 .\end{cases}
\] Then \(\operatorname{coh}\) takes values in \([-1,1]\), and
\(\gamma=\operatorname{coh}\cdot\sqrt{f_Xf_Y}\) everywhere.

\emph{Proof.} For any finite measurable partition \(\{A_k\}\) of a Borel
set \(A\), Lemma 8.4(c) and the arithmetic--geometric mean inequality
give \[
\sum_k|\Gamma(A_k)|\le\sum_k\sqrt{F_{XX}(A_k)F_{YY}(A_k)}\le\sum_k\frac{F_{XX}(A_k)+F_{YY}(A_k)}{2}=\frac{\Lambda(A)}{2},
\] and taking the supremum over partitions gives
\(|\Gamma|(A)\le\Lambda(A)/2\). \(F_{XX},F_{YY}\le\Lambda\) is obvious.
By the Radon--Nikodym theorem (S3) the three densities exist; choose
versions with \(f_X,f_Y\ge0\). Since \(F_{XX}+F_{YY}=\Lambda\),
\(f_X+f_Y=1\) holds \(\Lambda\)-almost everywhere, and modifying on a
null set makes it hold everywhere.

For a rational pair \((\alpha,\beta)\) and a Borel set \(A\), Lemma
8.4(c) gives \[
\int_A\big(\alpha^2f_X+2\alpha\beta\gamma+\beta^2f_Y\big)\,d\Lambda=\alpha^2F_{XX}(A)+2\alpha\beta\,\Gamma(A)+\beta^2F_{YY}(A)\ge0 .
\] Taking \(A\) to be the set where the integrand is negative shows that
the integrand is non-negative \(\Lambda\)-almost everywhere; let
\(N_{\alpha\beta}\) be the exceptional null set, so that
\(N=\bigcup_{(\alpha,\beta)\in\mathbb{Q}^2}N_{\alpha\beta}\) is still
null. Off \(N\), the quadratic form
\(\alpha^2f_X+2\alpha\beta\gamma+\beta^2f_Y\) is non-negative for all
rational pairs, hence by continuity for all real pairs; that is, the
matrix \(\begin{pmatrix}f_X&\gamma\\\gamma&f_Y\end{pmatrix}\) is
positive semi-definite and \(\gamma^2\le f_Xf_Y\). Redefine \(\gamma=0\)
on \(N\) (still a density of \(\Gamma\)); then \(\gamma^2\le f_Xf_Y\)
everywhere. Hence \(|\operatorname{coh}|\le1\) where \(f_Xf_Y>0\), and
where \(f_Xf_Y=0\) we have \(\gamma=0\), so both sides of
\(\gamma=\operatorname{coh}\sqrt{f_Xf_Y}\) vanish. \(\blacksquare\)

Up to \(\Lambda\)-null sets, \(\operatorname{coh}\) does not depend on
the choice of dominating measure (changing the dominating measure
multiplies the three densities by a common factor that is positive
\(\Lambda\)-almost everywhere). When the spectral matrix has Lebesgue
densities,
\(\operatorname{coh}=\operatorname{Re}f_{XY}/\sqrt{f_{XX}f_{YY}}\)
\(\Lambda\)-almost everywhere on \(\{f_{XX}f_{YY}>0\}\), the real part
of the complex coherency: it is the correlation of the two fields at
frequency \(\omega\).

\textbf{Theorem 8.9 (coherence bound).} For any zoning design \(\pi\)
with \(F_{XX}(g_\pi)F_{YY}(g_\pi)>0\), \[
\begin{aligned}
\rho(\pi)&=\theta(\pi)\,\bar c(\pi),\qquad
\bar c(\pi)=\frac{\int g_\pi\,\operatorname{coh}\,\sqrt{f_Xf_Y}\,d\Lambda}{\int g_\pi\sqrt{f_Xf_Y}\,d\Lambda},\\
\theta(\pi)&=\frac{\int g_\pi\sqrt{f_Xf_Y}\,d\Lambda}{\sqrt{F_{XX}(g_\pi)\,F_{YY}(g_\pi)}}\in[0,1],
\end{aligned}
\] with the convention \(\bar c(\pi)=0\) when
\(\int g_\pi\sqrt{f_Xf_Y}\,d\Lambda=0\), in which case \(\rho(\pi)=0\).
\(\bar c(\pi)\) is a weighted average of \(\operatorname{coh}\) with
respect to the positive measure \(g_\pi\sqrt{f_Xf_Y}\,\Lambda\), so \[
\min\Big(0,\ \operatorname*{ess\,inf}_{\Lambda}\operatorname{coh}\Big)\le\rho(\pi)\le\max\Big(0,\ \operatorname*{ess\,sup}_{\Lambda}\operatorname{coh}\Big),\qquad |\rho(\pi)|\le\operatorname*{ess\,sup}_{\Lambda}|\operatorname{coh}| .
\] The fine-scale correlation (\(g\equiv1\)) satisfies the same bounds.
\(\theta(\pi)=1\) if and only if \(f_X\) and \(f_Y\) are proportional
(\(\Lambda\)-almost everywhere) on \(\{g_\pi>0\}\).

\emph{Proof.} Write \(g=g_\pi\). By Lemmas 8.6 and 8.8, \[
\rho(\pi)=\frac{\int g\gamma\,d\Lambda}{\sqrt{\int gf_X\,d\Lambda\cdot\int gf_Y\,d\Lambda}},\qquad \int g\gamma\,d\Lambda=\int g\,\operatorname{coh}\sqrt{f_Xf_Y}\,d\Lambda .
\] Let \(w=g\sqrt{f_Xf_Y}\,\Lambda\), a finite positive measure. By the
Cauchy--Schwarz inequality, \[
w(\mathbb{R}^q)=\int\sqrt{gf_X}\sqrt{gf_Y}\,d\Lambda\le\Big(\int gf_X\,d\Lambda\Big)^{1/2}\Big(\int gf_Y\,d\Lambda\Big)^{1/2},
\] so \(\theta\in[0,1]\). If \(w(\mathbb{R}^q)=0\), the numerator
\(\int\operatorname{coh}\,dw\) vanishes and \(\rho(\pi)=0\). Otherwise
\(\rho(\pi)=\theta\,\bar c\), where
\(\bar c=\int\operatorname{coh}\,dw/w(\mathbb{R}^q)\) lies between the
essential infimum and supremum of \(\operatorname{coh}\) with respect to
\(w\); since \(w\ll\Lambda\), these do not exceed the essential bounds
with respect to \(\Lambda\). Multiplying \(\bar c\) by
\(\theta\in[0,1]\) only shrinks it towards \(0\), which gives the first
pair of inequalities; the second follows from
\(|\bar c|\le\operatorname{ess\,sup}_\Lambda|\operatorname{coh}|\).
\(\theta=1\) is the equality case of the Cauchy--Schwarz inequality,
i.e.~\(\sqrt{gf_X}\) and \(\sqrt{gf_Y}\) are proportional in
\(L^2(\Lambda)\), which means that \(f_X\) and \(f_Y\) are proportional
on \(\{g>0\}\) (when neither vanishes identically there).
\(\blacksquare\)

\textbf{Remark (aggregation averages the correlations at all
frequencies).} Theorem 8.9 says that the zonal correlation is a weighted
average of the frequency-wise correlations
\(\operatorname{coh}(\omega)\), the weights being the filter of the
design \(g_\pi\) times the common power \(\sqrt{f_Xf_Y}\) of the two
fields at that frequency, multiplied by a factor \(\theta(\pi)\le1\)
that measures how different the two power spectra are on
\(\{g_\pi>0\}\). The larger the zones, the more the filter concentrates
at low frequencies and the more weight goes to the coherence at long
wavelengths. This yields a hard constraint opposite to Section 6: if
\(|\operatorname{coh}|\le c_0\) \(\Lambda\)-almost everywhere, then
every data-independent zoning design satisfies \(|\rho(\pi)|\le c_0\).
Data-dependent zonings can produce any zonal correlation;
data-independent designs cannot leave the interval between \(0\) and the
extreme values of the coherence.

\hypertarget{a-necessary-and-sufficient-condition-for-scale-invariance}{%
\subsection{8.4 A necessary and sufficient condition for scale
invariance}\label{a-necessary-and-sufficient-condition-for-scale-invariance}}

Let \(\mathcal{K}_c\) be the set of uniform densities on non-empty
bounded connected open sets in \(\mathbb{R}^q\) (\textbf{connected zone
shapes}). Let \(\mathcal{P}\) be the set of finitely supported designs
on \(\mathcal{K}_c\) with \(F_{XX}(g_\pi)F_{YY}(g_\pi)>0\), and
\(\mathcal{P}_2\subseteq\mathcal{P}\) the designs in \(\mathcal{P}\)
made of at most two zone shapes.

\textbf{Lemma 8.10 (anchor design and separation lemma).}

\begin{enumerate}
\def\labelenumi{(\alph{enumi})}
\item
  For \(r>0\) let \(u_r\) be the uniform density on the ball \(B(0,r)\).
  Then \(\hat u_r(\omega)=\Phi(r|\omega|)\), where
  \(\Phi:[0,\infty)\to\mathbb{R}\) is continuous with \(\Phi(0)=1\).
  Hence there is \(t_0>0\) such that \(|\hat u_r(\omega)|^2\ge1/4\)
  whenever \(|\omega|\le t_0/r\). If \(G\) is a non-zero finite positive
  measure, then \(G(|\hat u_r|^2)>0\) for all sufficiently small \(r\).
\item
  If \(\nu\) is a finite even signed measure on \(\mathbb{R}^q\) with
  \(\int|\hat K|^2\,d\nu=0\) for all \(K\in\mathcal{K}_c\), then
  \(\nu=0\).
\end{enumerate}

\emph{Proof.} (a) A ball is a bounded connected open set, so
\(u_r\in\mathcal{K}_c\). For an orthogonal matrix \(O\), substituting
\(s=Os'\) and using \(u_r(Os')=u_r(s')\), \[
\hat u_r(O\omega)=\int e^{-is\cdot O\omega}u_r(s)\,ds=\int e^{-iO^\top s\cdot\omega}u_r(s)\,ds=\int e^{-is'\cdot\omega}u_r(Os')\,ds'=\hat u_r(\omega),
\] so \(\hat u_r(\omega)\) depends only on \(|\omega|\). Also
\(u_r=(u_1)_r\), so by Lemma 8.5(d)
\(\hat u_r(\omega)=\hat u_1(r\omega)\). Let \(\Phi(t)=\hat u_1(te_1)\):
since \(u_1\) is even, \(\hat u_1\) is real-valued, and by Lemma 8.5(c)
\(\Phi\) is continuous with \(\Phi(0)=1\). Choose \(t_0>0\) with
\(\Phi(t)^2\ge1/4\) for \(t\le t_0\). Finally
\(G(|\hat u_r|^2)\ge\frac14G\big(B(0,t_0/r)\big)\), and
\(G(B(0,t_0/r))\to G(\mathbb{R}^q)>0\) as \(r\to0\).

\begin{enumerate}
\def\labelenumi{(\alph{enumi})}
\setcounter{enumi}{1}
\tightlist
\item
  \emph{The case \(q\ge2\).} Balls give \(\int|\hat u_r|^2\,d\nu=0\) for
  all \(r>0\). Fix \(a\in\mathbb{R}^q\) with \(|a|>2r\) and
  \(0<\varepsilon<r\). Let \(T_\varepsilon\) be the open
  \(\varepsilon\)-neighbourhood of the segment \([0,a]\) and
  \(D_\varepsilon=B(0,r)\cup B(a,r)\cup T_\varepsilon\).
  \(T_\varepsilon\) is convex and meets both balls, so \(D_\varepsilon\)
  is a bounded connected open set and
  \(K_\varepsilon=\mathbf{1}_{D_\varepsilon}/|D_\varepsilon|\in\mathcal{K}_c\).
  \(T_\varepsilon\) is contained in a cylinder of length
  \(|a|+2\varepsilon\) whose cross-section is a \((q-1)\)-dimensional
  ball of radius \(\varepsilon\), so \(|T_\varepsilon|\to0\) when
  \(q\ge2\). The two balls are disjoint; let
  \(K_0=\big(\mathbf{1}_{B(0,r)}+\mathbf{1}_{B(a,r)}\big)/(2|B(0,r)|)\).
  Then \(2|B(0,r)|\le|D_\varepsilon|\le2|B(0,r)|+|T_\varepsilon|\) and
  \[
  \|K_\varepsilon-K_0\|_1\le\Big|\frac{1}{|D_\varepsilon|}-\frac{1}{2|B(0,r)|}\Big|\cdot2|B(0,r)|+\frac{|T_\varepsilon|}{|D_\varepsilon|}\longrightarrow0 .
  \] By Lemma 8.5(c), \(\hat K_\varepsilon\to\hat K_0\) uniformly, hence
  \(|\hat K_\varepsilon|^2\to|\hat K_0|^2\) uniformly; \(\nu\) is
  finite, so
  \(\int|\hat K_0|^2\,d\nu=\lim_\varepsilon\int|\hat K_\varepsilon|^2\,d\nu=0\).
  By Lemma 8.5(d), \[
  \begin{aligned}\hat K_0(\omega)&=\frac12\big(1+e^{-ia\cdot\omega}\big)\hat u_r(\omega),\\|\hat K_0(\omega)|^2&=|\hat u_r(\omega)|^2\,\frac{|1+e^{-ia\cdot\omega}|^2}{4}=|\hat u_r(\omega)|^2\,\frac{1+\cos(a\cdot\omega)}{2}.\end{aligned}
  \] Subtracting the identity for the ball gives
  \(\int|\hat u_r(\omega)|^2\cos(a\cdot\omega)\,\nu(d\omega)=0\). Fix
  \(a\neq0\) and let \(r\to0\) (keeping \(r<|a|/2\)):
  \(|\hat u_r(\omega)|^2=\Phi(r|\omega|)^2\to1\) pointwise and is
  bounded by \(1\), so by dominated convergence (with respect to
  \(|\nu|\)) \(\int\cos(a\cdot\omega)\,\nu(d\omega)=0\). Likewise
  \(\int|\hat u_r|^2\,d\nu=0\) gives \(\nu(\mathbb{R}^q)=0\). Since
  \(\nu\) is even, \(\int\sin(a\cdot\omega)\,\nu(d\omega)=0\) (same
  argument as in Lemma 8.4(b)), so \(\hat\nu(a)=0\) for all \(a\), and
  \(\nu=0\) by Theorem 8.1.
\end{enumerate}

\emph{The case \(q=1\).} Here \(\mathcal{K}_c\) consists of the uniform
densities on intervals. The interval \([0,T]\) of length \(T\) gives \[
\hat K(\omega)=\frac{1-e^{-iT\omega}}{iT\omega},\qquad |\hat K(\omega)|^2=\frac{2(1-\cos T\omega)}{T^2\omega^2}\quad(\omega\neq0),\qquad |\hat K(0)|^2=1 .
\] Let \(h(T,\omega)=(1-\cos T\omega)/\omega^2\) (with
\(h(T,0)=T^2/2\)), so that \(|\hat K|^2=2h(T,\cdot)/T^2\) and the
hypothesis becomes \(H(T)=\int h(T,\omega)\,\nu(d\omega)=0\) for all
\(T>0\). For each \(\omega\),
\(\partial_Th(T,\omega)=\sin(T\omega)/\omega\) (equal to \(T\) at
\(\omega=0\)) and \(\partial_T^2h(T,\omega)=\cos(T\omega)\), with
\(|\partial_Th|\le T\) and \(|\partial_T^2h|\le1\). \(\nu\) is finite,
so by the domination condition for differentiation under the integral
sign (S2), \(H\) is twice differentiable with
\(H''(T)=\int\cos(T\omega)\,\nu(d\omega)\). Hence
\(\int\cos(T\omega)\,\nu(d\omega)=0\) for all \(T>0\). Since \(\nu\) is
even, \(\hat\nu(T)=\int\cos(T\omega)\,\nu(d\omega)\) is an even
function, so \(\hat\nu(T)=0\) for \(T\neq0\), and at \(T=0\) by
continuity of \(\hat\nu\). By Theorem 8.1, \(\nu=0\). \(\blacksquare\)

\textbf{Theorem 8.11 (scale invariance and the intrinsic correlation
model).} Let \(F_{XX}\neq0\), \(F_{YY}\neq0\) and \(\rho_0\in[-1,1]\).

(Necessity) If \(\rho(\pi)=\rho_0\) for all \(\pi\in\mathcal{P}_2\),
then:

\begin{enumerate}
\def\labelenumi{(\alph{enumi})}
\item
  if \(\rho_0\neq0\), there is \(\lambda>0\) with
  \(F_{YY}=\lambda F_{XX}\) and \(\Gamma=\rho_0\sqrt\lambda\,F_{XX}\);
  equivalently, \(C_{YY}=\lambda C_{XX}\) and
  \(\tfrac12\big(C_{XY}(h)+C_{XY}(-h)\big)=\rho_0\sqrt\lambda\,C_{XX}(h)\)
  for all \(h\);
\item
  if \(\rho_0=0\), then \(\Gamma=0\); equivalently, the even part of the
  cross-covariance vanishes identically.
\end{enumerate}

(Sufficiency) Conversely, if (a) or (b) holds, then \(\rho(\pi)=\rho_0\)
for every finite design \(\pi\) on arbitrary (not necessarily connected)
aggregation kernels with \(F_{XX}(g_\pi)F_{YY}(g_\pi)>0\), and the
fine-scale correlation is \(\rho=\rho_0\).

\emph{Proof.}

\emph{Sufficiency.} Under (a),
\(\rho(\pi)=\rho_0\sqrt\lambda F_{XX}(g)\big/\sqrt{F_{XX}(g)\cdot\lambda F_{XX}(g)}=\rho_0\);
under (b) the numerator vanishes. Taking \(g\equiv1\) gives
\(\rho=\rho_0\) (since \(F_{XX}(\mathbb{R}^q)F_{YY}(\mathbb{R}^q)>0\)).

\emph{Necessity, Step 1: the anchor.} By Lemma 8.10(a), choose \(r\) so
small that \(a_u=F_{XX}(|\hat u_r|^2)>0\) and
\(b_u=F_{YY}(|\hat u_r|^2)>0\); write \(u=u_r\) and
\(c_u=\Gamma(|\hat u|^2)\), so that \(\delta_u\in\mathcal{P}_2\). Let
\(\lambda=b_u/a_u>0\). From \(\rho(\delta_u)=\rho_0\),
\(c_u=\rho_0\sqrt{a_ub_u}=\rho_0\sqrt\lambda\,a_u\).

\emph{Step 2: mixing with the anchor.} Take any \(K\in\mathcal{K}_c\)
and write \(a=F_{XX}(|\hat K|^2)\), \(b=F_{YY}(|\hat K|^2)\),
\(c=\Gamma(|\hat K|^2)\) (\(a,b\) may vanish). For \(t\in(0,1]\) the
design \(\pi_t=t\delta_u+(1-t)\delta_K\) has at most two zone shapes and
\(g_{\pi_t}\ge t|\hat u|^2\), so \(F_{XX}(g_{\pi_t})\ge ta_u>0\),
\(F_{YY}(g_{\pi_t})\ge tb_u>0\), and \(\pi_t\in\mathcal{P}_2\). By
linearity of \(g\) in the design, \(\rho(\pi_t)=\rho_0\) reads \[
tc_u+(1-t)c=\rho_0\sqrt{\big(ta_u+(1-t)a\big)\big(tb_u+(1-t)b\big)} .
\] Dividing by \(t\) and putting \(\sigma=(1-t)/t\in[0,\infty)\), \[
c_u+\sigma c=\rho_0\sqrt{(a_u+\sigma a)(\lambda a_u+\sigma b)}\qquad(\sigma\ge0).\tag{8.1}
\]

\emph{Step 3: \(\rho_0\neq0\).} Squaring (8.1), both sides are
polynomials of degree at most two in \(\sigma\) that agree on
\([0,\infty)\), so their coefficients coincide: \[
\sigma^0:\ c_u^2=\rho_0^2\lambda a_u^2;\qquad
\sigma^1:\ 2c_uc=\rho_0^2a_u\,(b+\lambda a);\qquad
\sigma^2:\ c^2=\rho_0^2\,ab .
\] The \(\sigma^0\) term agrees with Step 1. Substituting
\(c_u=\rho_0\sqrt\lambda\,a_u\), the \(\sigma^1\) term gives
\(c=\rho_0(b+\lambda a)/(2\sqrt\lambda)\); substituting into the
\(\sigma^2\) term and cancelling \(\rho_0^2>0\), \[
\frac{(b+\lambda a)^2}{4\lambda}=ab\iff(b-\lambda a)^2=0 .
\] So \(b=\lambda a\) and \(c=\rho_0\sqrt\lambda\,a\). This holds for
every \(K\in\mathcal{K}_c\); that is, \(\nu_1=F_{YY}-\lambda F_{XX}\)
and \(\nu_2=\Gamma-\rho_0\sqrt\lambda\,F_{XX}\) satisfy
\(\int|\hat K|^2\,d\nu_i=0\) for \(K\in\mathcal{K}_c\). By Lemma 8.4(a),
\(\nu_1,\nu_2\) are finite even signed measures, and by Lemma 8.10(b)
\(\nu_1=\nu_2=0\).

\emph{Step 4: \(\rho_0=0\).} Now (8.1) reads \(c_u+\sigma c=0\) for
\(\sigma\ge0\). Taking \(\sigma=0\) gives \(c_u=0\), and then
\(\sigma=1\) gives \(c=0\). Hence \(\int|\hat K|^2\,d\Gamma=0\) for all
\(K\in\mathcal{K}_c\), and \(\Gamma=0\) by Lemma 8.10(b).

\emph{Step 5: covariance-function form.} By Lemma 8.4(d),
\(C_{UU}(h)=\int\cos(h\cdot\omega)\,dF_{UU}\) and
\(\tfrac12(C_{XY}(h)+C_{XY}(-h))=\int\cos(h\cdot\omega)\,d\Gamma\); by
the last sentence of Lemma 8.4(d) (an even measure is determined by its
cosine transform), the measure form and the covariance-function form are
equivalent. \(\blacksquare\)

Condition (a) is the \textbf{intrinsic correlation model} of
multivariate geostatistics (Wackernagel 2003): the auto-covariance
functions of the two variables (and hence their variograms) and the even
part of their cross-covariance function have the same shape and differ
only by constant factors. The classical intrinsic correlation model also
requires the cross-covariance to be symmetric; by Lemma 8.4(b),(d), the
odd part of the cross-covariance (for example a relative shift between
the two fields) is invisible to every design, so the condition of
Theorem 8.11 constrains only the even part. The sufficiency direction is
known; the necessity direction shows that ``the correlation does not
change with scale or zone shape'' characterizes exactly the intrinsic
correlation model when the common value is non-zero, and the vanishing
of the co-spectrum when it is zero, and that only designs made of at
most two connected zone shapes need to be examined.

\textbf{Corollary 8.12 (a deviating design must exist).} Let
\(F_{XX}\neq0\), \(F_{YY}\neq0\), and let \(u\) be the anchor ball of
Step 1 in the proof of Theorem 8.11.

\begin{enumerate}
\def\labelenumi{(\roman{enumi})}
\item
  \(\rho(\pi)\) is constant on \(\mathcal{P}\) if and only if (a) or (b)
  holds for some \(\rho_0\); the constant is then the fine-scale
  correlation \(\rho\).
\item
  If the fine-scale correlation \(\rho\neq0\) and the intrinsic
  correlation model (a) fails, then there are a connected zone shape
  \(K\in\mathcal{K}_c\) and \(t\in(0,1]\) such that the design
  \(t\delta_u+(1-t)\delta_K\) (at most two connected zone shapes)
  satisfies \(\rho(t\delta_u+(1-t)\delta_K)\neq\rho\).
\item
  If \(\rho=0\) but \(\Gamma\neq0\), there is likewise a design of this
  form whose zonal correlation is non-zero.
\end{enumerate}

\emph{Proof.} (i) If \(\rho(\pi)\equiv\rho_0\) on \(\mathcal{P}\), the
same holds on \(\mathcal{P}_2\), and necessity in Theorem 8.11 gives (a)
or (b); conversely, by sufficiency \(\rho(\pi)\equiv\rho_0\) and
\(\rho=\rho_0\).

\begin{enumerate}
\def\labelenumi{(\roman{enumi})}
\setcounter{enumi}{1}
\item
  Suppose that all designs \(\pi_t=t\delta_u+(1-t)\delta_K\)
  (\(K\in\mathcal{K}_c\), \(t\in(0,1]\)) satisfy \(\rho(\pi_t)=\rho\).
  The proof of necessity in Theorem 8.11 (Steps 1--3) uses only these
  designs; with \(\rho_0=\rho\neq0\) it gives (a), a contradiction.
\item
  Likewise, with \(\rho_0=0\), Step 4 gives \(\Gamma=0\), a
  contradiction. \(\blacksquare\)
\end{enumerate}

Corollary 8.12(ii) means: as soon as the two power spectra are not
proportional, or the co-spectrum is not proportional to the common power
spectrum, while the fine-scale correlation is non-zero, some design of
connected zone shapes (the small ball alone, or the small ball mixed in
some proportion with one other connected shape) has a zonal correlation
different from the fine-scale correlation. The proof uses mixtures of
two shapes; whether a single shape always suffices is not addressed
here.

\textbf{Example (linear model of coregionalization).} Let \(A\) and
\(B\) be independent zero-mean stationary fields with non-zero spectral
measures \(F_A\) and \(F_B\), and put \(X=A\), \(Y=\alpha A+\beta B\)
(\(\alpha\neq0\), \(\beta\neq0\)). By independence, \(C_{XX}=C_A\),
\(C_{YY}=\alpha^2C_A+\beta^2C_B\) and \(C_{XY}=\alpha C_A\), so by
Theorem 8.1 \[
F_{XX}=F_A,\qquad F_{YY}=\alpha^2F_A+\beta^2F_B,\qquad \Gamma=F_{XY}=\alpha F_A .
\] If \(F_B=\mu F_A\) (\(\mu>0\)), then (a) holds with
\(\lambda=\alpha^2+\beta^2\mu\) and
\(\rho_0=\alpha/\sqrt{\alpha^2+\beta^2\mu}\); for example \(\mu=1\),
\(\alpha=0.6\), \(\beta=0.8\) give \(\rho_0=0.6\), which is the setting
of check V5 in Section 12. Conversely, if \(F_{YY}=\lambda F_{XX}\),
then \(\beta^2F_B=(\lambda-\alpha^2)F_A\), i.e.~\(F_B\) is proportional
to \(F_A\). So the zonal correlation is invariant over all connected
designs if and only if the spectra of the two components are
proportional. When the spectra have densities, \[
\operatorname{coh}(\omega)=\frac{\alpha}{\sqrt{\alpha^2+\beta^2f_B(\omega)/f_A(\omega)}}\qquad(f_A(\omega)>0),
\] which varies with frequency unless \(f_B/f_A\) is constant. If \(B\)
is short-range (\(f_B/f_A\) small at low and large at high frequencies),
the zonal correlation typically strengthens as the zones grow and, by
Theorem 8.15 (when the densities are bounded and continuous at the
origin), converges to \(\alpha/\sqrt{\alpha^2+\beta^2f_B(0)/f_A(0)}\).

\textbf{Proposition 8.13 (finitely many designs cannot certify the
intrinsic correlation model).} Let \(\pi_1,\dots,\pi_n\) be finite
designs on arbitrary aggregation kernels, let \(F\) be a finite measure
on \(\mathbb{R}^q\) with an everywhere positive even density (for
example a Gaussian spectral density), and let \(\rho_0\in(-1,1)\). Then
there are a bounded measurable even function \(v\) (not \(F\)-almost
everywhere zero) and \(\varepsilon>0\) such that the spectral matrix \[
F_{XX}=F_{YY}=F,\qquad F_{XY}=\Gamma=(\rho_0+\varepsilon v)\,F\tag{8.2}
\] belongs to a pair of jointly stationary fields satisfying the setting
of Section 8.1, with \(\rho(\pi_i)=\rho_0\) (\(i=1,\dots,n\)) and
fine-scale correlation \(\rho=\rho_0\), while neither (a) nor (b) holds.
Consequently, by Corollary 8.12, some design made of at most two
connected zone shapes has a zonal correlation different from \(\rho_0\).

\emph{Proof.} \emph{Construction of \(v\).} Let
\(E_k=\{\omega:k-1\le|\omega|<k\}\) (\(k=1,\dots,n+2\)); these are
pairwise disjoint symmetric sets with \(F(E_k)>0\). For
\(t\in\mathbb{R}^{n+2}\) let \(v_t=\sum_kt_k\mathbf{1}_{E_k}\) and
consider the linear map \[
t\longmapsto\Big(\int g_{\pi_1}v_t\,dF,\ \dots,\ \int g_{\pi_n}v_t\,dF,\ \int v_t\,dF\Big)\in\mathbb{R}^{n+1}.
\] The domain has dimension \(n+2\), larger than that of the range, so
some \(t\neq0\) lies in the kernel; rescale it to \(\max_k|t_k|=1\) and
put \(v=v_t\). Then \(v\) is even, \(\|v\|_\infty=1\), and \(v\) takes a
non-zero constant value on some \(E_k\), so it is not \(F\)-almost
everywhere zero; moreover \(\int g_{\pi_i}v\,dF=0\) and
\(\int v\,dF=0\). Take \(\varepsilon=1-|\rho_0|>0\); then
\(c=\rho_0+\varepsilon v\) is a bounded even function with values in
\([-1,1]\).

\emph{Realization by random-phase fields.} \(F_1=\tfrac12(1+c)F\) and
\(F_2=\tfrac12(1-c)F\) are finite even positive measures. For \(i=1,2\)
with \(F_i\neq0\), let \(W_i\) be a random vector with law
\(F_i/F_i(\mathbb{R}^q)\) and \(\Phi_i\) a random variable uniform on
\([0,2\pi)\), all four independent, and put \[
S_i(s)=\sqrt{2F_i(\mathbb{R}^q)}\,\cos(s\cdot W_i+\Phi_i)\qquad(s\in\mathbb{R}^q);
\] if \(F_i=0\), put \(S_i\equiv0\). Each \(S_i\) is bounded and
continuous in \(s\), hence jointly measurable, and
\(\mathbb{E}S_i(s)=0\) because \(\Phi_i\) is uniform. By
\(\cos\alpha\cos\beta=\tfrac12\big(\cos(\alpha-\beta)+\cos(\alpha+\beta)\big)\)
and \(\mathbb{E}\cos(x+2\Phi_i)=0\) for every real \(x\), \[
\mathbb{E}\big[S_i(t+h)S_i(t)\big]=F_i(\mathbb{R}^q)\,\mathbb{E}\cos(h\cdot W_i)=\int\cos(h\cdot\omega)\,F_i(d\omega)=\hat F_i(h),
\] the last equality because \(F_i\) is even (as in Lemma 8.4(b), the
sine part integrates to zero). Put \(X=S_1+S_2\) and \(Y=S_1-S_2\). By
independence and zero means, \[
C_{XX}=C_{YY}=\hat F_1+\hat F_2=\hat F,\qquad C_{XY}=\hat F_1-\hat F_2=\widehat{cF}.
\] These do not depend on \(t\) and are continuous, so \((X,Y)\)
satisfies the setting of Section 8.1, and by the uniqueness in Theorem
8.3 its spectral matrix is (8.2). (Gaussian fields with this spectral
matrix would also do, but their construction needs Kolmogorov's
extension theorem and a measurable modification for the joint
measurability required in Section 8.1; the random-phase construction
needs neither.)

\emph{Verification.} The filters satisfy \(g_{\pi_i}(0)=1\) and are
continuous, and \(F\) has an everywhere positive density, so
\(F(g_{\pi_i})>0\). Hence \[
\rho(\pi_i)=\frac{\rho_0F(g_{\pi_i})+\varepsilon\int g_{\pi_i}v\,dF}{F(g_{\pi_i})}=\rho_0,\qquad
\rho=\frac{\rho_0F(\mathbb{R}^q)+\varepsilon\int v\,dF}{F(\mathbb{R}^q)}=\rho_0 .
\] If (a) held for some \(\rho_0'\neq0\) and \(\lambda>0\), then
\(F_{XX}=F_{YY}\) would force \(\lambda=1\), sufficiency in Theorem 8.11
would give \(\rho_0'=\rho=\rho_0\), and hence \(\Gamma=\rho_0F\); but
\(\Gamma-\rho_0F=\varepsilon vF\neq0\). If (b) held, then \(\Gamma=0\);
but for \(\rho_0=0\) we have \(\Gamma=\varepsilon vF\neq0\), and for
\(\rho_0\neq0\) (b) contradicts \(\rho=\rho_0\neq0\). The last sentence
follows from Corollary 8.12(ii)(iii). \(\blacksquare\)

Proposition 8.13 shows that observing similar (or even identical) zonal
correlations at several chosen scales or zone shapes, together with the
same fine-scale correlation, is still not enough to conclude the
intrinsic correlation model, and hence not enough to conclude that the
correlation is the same for other designs. Testing scale invariance
requires testing proportionality of the spectra (equivalently, of the
auto-covariances and the even part of the cross-covariance at all lags),
not comparing correlations at finitely many scales.

\hypertarget{a-necessary-and-sufficient-condition-for-aggregation-bias}{%
\subsection{8.5 A necessary and sufficient condition for aggregation
bias}\label{a-necessary-and-sufficient-condition-for-aggregation-bias}}

\textbf{Theorem 8.14 (aggregation bias = symmetrized lagged
cross-covariance).} Let \(Y(s)=\alpha+\beta X(s)+e(s)\), where \((X,e)\)
satisfies the setting of Section 8.1 (\(e\) has mean zero), with
cross-spectral measure \(F_{Xe}\) and \(F_{XX}\neq0\). For designs
\(\pi\) with \(F_{XX}(g_\pi)>0\), define the design-level zonal
regression slope \(\beta(\pi)=\Sigma_{XY}(\pi)/\Sigma_{XX}(\pi)\). Then
\[
\beta(\pi)=\beta+\frac{\int g_\pi\,d(\operatorname{Re}F_{Xe})}{F_{XX}(g_\pi)} .
\] For every \(q\ge1\) the following are equivalent: (i)
\(\beta(\pi)=\beta\) for all finite designs on \(\mathcal{K}_c\) with
\(F_{XX}(g_\pi)>0\); (ii) \(\operatorname{Re}F_{Xe}=0\); (iii)
\(C_{Xe}(h)+C_{Xe}(-h)=0\) for all \(h\). Moreover, (i) needs to be
checked only for one ball \(u=u_r\) with \(F_{XX}(|\hat u_r|^2)>0\)
(Lemma 8.10(a)) and the designs \(\tfrac12\delta_u+\tfrac12\delta_K\)
(\(K\in\mathcal{K}_c\)). In particular, being uncorrelated at the same
location, \(C_{Xe}(0)=0\), does not guarantee unbiased aggregation.

\emph{Proof.} \(Y-\alpha=\beta X+e\) has mean zero, \((X,Y-\alpha)\)
satisfies the setting of Section 8.1, and the zonal second moments are
computed with \(Y-\alpha\). Since \(\int K=1\),
\(Y_K=\alpha+\beta X_K+e_K\), so
\(\operatorname{Cov}(X_K,Y_K)=\beta\operatorname{Var}(X_K)+\operatorname{Cov}(X_K,e_K)\).
Applying Lemma 8.6 to \((X,e)\) and averaging with weights \(\pi\) gives
\(\Sigma_{XY}(\pi)=\beta F_{XX}(g_\pi)+\operatorname{Re}F_{Xe}(g_\pi)\),
which is the formula.

(ii)\(\Rightarrow\)(i) is immediate from the formula.
(i)\(\Rightarrow\)(ii): let \(u=u_r\) be such a ball. From
\(\beta(\delta_u)=\beta\), \(\operatorname{Re}F_{Xe}(|\hat u|^2)=0\);
for \(K\in\mathcal{K}_c\), \(\pi=\tfrac12\delta_u+\tfrac12\delta_K\) has
\(F_{XX}(g_\pi)>0\), and \(\beta(\pi)=\beta\) together with linearity
gives \(\operatorname{Re}F_{Xe}(|\hat K|^2)=0\). By Lemma 8.4(a)
(applied to \((X,e)\)), \(\operatorname{Re}F_{Xe}\) is a finite even
signed measure, and by Lemma 8.10(b) \(\operatorname{Re}F_{Xe}=0\).
(ii)\(\Leftrightarrow\)(iii): by Lemma 8.4(d),
\(\tfrac12\big(C_{Xe}(h)+C_{Xe}(-h)\big)=\int\cos(h\cdot\omega)\,d(\operatorname{Re}F_{Xe})\),
and an even measure is determined by its cosine transform. For the last
sentence:
\(C_{Xe}(0)=F_{Xe}(\mathbb{R}^q)=\operatorname{Re}F_{Xe}(\mathbb{R}^q)\)
constrains only the total mass of \(\operatorname{Re}F_{Xe}\), not its
distribution over frequencies. \(\blacksquare\)

\textbf{Example (spatial spillover).} Let \(C_{XX}(0)>0\), \(h_0\neq0\),
\(c=C_{XX}(h_0)/C_{XX}(0)\), and \[
e(s)=\gamma_0\big[X(s+h_0)+X(s-h_0)-2cX(s)\big].
\] \((X,e)\) satisfies the setting of Section 8.1, and \[
\begin{aligned}C_{Xe}(h)&=\gamma_0\big[C_{XX}(h-h_0)+C_{XX}(h+h_0)-2cC_{XX}(h)\big]\\&=\int e^{ih\cdot\omega}\,2\gamma_0\big(\cos(h_0\cdot\omega)-c\big)\,F_{XX}(d\omega).\end{aligned}
\] Hence \(C_{Xe}(0)=2\gamma_0\big(C_{XX}(h_0)-cC_{XX}(0)\big)=0\),
while by Theorem 8.1
\(F_{Xe}=\operatorname{Re}F_{Xe}=2\gamma_0(\cos(h_0\cdot\omega)-c)F_{XX}\).
This vanishes if and only if \(\gamma_0=0\) or \(F_{XX}\) is
concentrated on \(\{\omega:\cos(h_0\cdot\omega)=c\}\); when \(F_{XX}\)
has a density, that set is a countable union of hyperplanes and has
Lebesgue measure zero, so whenever \(\gamma_0\neq0\) the aggregated
slope is biased for some connected designs. If moreover \(K\in L^2\) and
\(F_{XX}\) has a bounded density that is continuous and positive at the
origin, then, as in the proof of Theorem 8.15(ii),
\(\beta(\delta_{K_T})\to\beta+2\gamma_0(1-c)\) as \(T\to\infty\), while
\(\beta(\delta_{K_T})\to\beta+C_{Xe}(0)/C_{XX}(0)=\beta\) as \(T\to0\):
unbiased at the finest scale, while at the coarsest scale the slope
includes the spillover effect. In check V9 of Section 12,
\(\operatorname{Cov}(X,e)=0\) at the same location, while the aggregated
slope rises from \(1.000\) to \(1.342\) as the block side grows from 1
to 64 pixels (Figure 8(c)).

Theorem 8.14 places aggregation bias in econometrics (Theil 1954) and
spatial confounding in spatial regression (Paciorek 2010; Guan et
al.~2023) under one spectral condition: unbiased aggregation requires
the symmetrized cross-covariance between error and explanatory variable
to vanish at all lags; checking exogeneity at the same location is not
enough.

\hypertarget{micro-and-macro-limits}{%
\subsection{8.6 Micro and macro limits}\label{micro-and-macro-limits}}

\textbf{Theorem 8.15 (micro and macro limits).} Let \(K\) be an
aggregation kernel and \(K_T(s)=T^{-q}K(s/T)\) (the same zone shape
enlarged by the factor \(T\)).

\begin{enumerate}
\def\labelenumi{(\roman{enumi})}
\item
  If \(C_{XX}(0)C_{YY}(0)>0\), then \(\rho(\delta_{K_T})\to\rho\) (the
  fine-scale correlation) as \(T\to0\).
\item
  If \(K\in L^2(\mathbb{R}^q)\), \(F_{XX}=f_{XX}\,d\omega\),
  \(F_{YY}=f_{YY}\,d\omega\) and \(\Gamma=c_{XY}\,d\omega\) with the
  three densities bounded and continuous at \(\omega=0\), and
  \(f_{XX}(0)f_{YY}(0)>0\), then \[
  \lim_{T\to\infty}\rho(\delta_{K_T})=\operatorname{coh}(0)=\frac{c_{XY}(0)}{\sqrt{f_{XX}(0)\,f_{YY}(0)}},
  \] whatever the shape of \(K\). Here \(\operatorname{coh}(0)\) is
  computed from the versions of the densities that are continuous at the
  origin; their values at \(0\) do not depend on the version.
\end{enumerate}

\emph{Proof.} By Lemma 8.5(d), \(\hat K_T(\omega)=\hat K(T\omega)\).

\begin{enumerate}
\def\labelenumi{(\roman{enumi})}
\item
  As \(T\to0\), \(\hat K(T\omega)\to\hat K(0)=1\) pointwise,
  \(|\hat K|\le1\) and \(F\) is finite, so by dominated convergence
  \(F_{UU}(|\hat K_T|^2)\to F_{UU}(\mathbb{R}^q)=C_{UU}(0)\) and
  \(\Gamma(|\hat K_T|^2)\to\Gamma(\mathbb{R}^q)=C_{XY}(0)\). For small
  \(T\) the denominator is positive, and the ratio tends to \(\rho\).
\item
  By Plancherel's theorem (S6),
  \(\int|\hat K(u)|^2\,du=(2\pi)^q\|K\|_2^2\in(0,\infty)\). For
  \(f\in\{f_{XX},f_{YY},c_{XY}\}\), substituting \(\omega=u/T\), \[
  T^{q}\int|\hat K(T\omega)|^2f(\omega)\,d\omega=\int|\hat K(u)|^2\,f(u/T)\,du\ \longrightarrow\ f(0)\int|\hat K(u)|^2\,du\quad(T\to\infty),
  \] since \(f(u/T)\to f(0)\) pointwise and
  \(|\hat K(u)|^2|f(u/T)|\le\|f\|_\infty|\hat K(u)|^2\in L^1\), so
  dominated convergence applies. For large \(T\),
  \(T^qF_{XX}(|\hat K_T|^2)\) and \(T^qF_{YY}(|\hat K_T|^2)\) are
  positive; substituting the three limits into \(\rho(\delta_{K_T})\),
  the common factors \(T^q\) and \(\int|\hat K|^2\) cancel.
  \(\blacksquare\)
\end{enumerate}

Theorem 8.15(ii) echoes the limit theory of temporal aggregation of time
series (Tiao 1972) and geostatistical work on how correlation changes
with support (Vargas-Guzmán, Warrick \& Myers 1999); here it is given in
general form under the spectral representation of stationary fields.
Both limits are independent of the zone shape; at intermediate scales
the zonal correlation is described by the weighted average of Theorem
8.9 and generally depends on the shape.

Theorem 8.15(i) says only that the zonal correlation tends to the
fine-scale correlation as the zones shrink. The next proposition gives
an explicit bound on the rate; it is the design-level counterpart of
Theorem 7.1 and Corollary 7.3, and it shows that the rate is set by the
smoothness of the fields.

\textbf{Proposition 8.16 (rate of the micro limit).} Let \((X,Y)\)
satisfy the setting of Section 8.1 with \(C_{XX}(0)C_{YY}(0)>0\). For a
finite design \(\pi\) let \[
a(\pi)=1-\frac{F_{XX}(g_\pi)}{C_{XX}(0)},\qquad b(\pi)=1-\frac{F_{YY}(g_\pi)}{C_{YY}(0)} .
\]

\begin{enumerate}
\def\labelenumi{(\roman{enumi})}
\item
  \(a(\pi)\) is the design-level within-zone variance share: for an
  aggregation kernel \(K\),
  \(\int K(s)\,\mathbb{E}\big[(X(s)-X_K)^2\big]\,ds=C_{XX}(0)-\operatorname{Var}(X_K)\).
  If \(a(\pi),b(\pi)<1\), then
  \(|\rho(\pi)-\rho|\le\mathcal{B}\big(a(\pi),b(\pi);\rho\big)\), where
  \(\mathcal{B}\) is the bound of Theorem 7.1.
\item
  Let the aggregation kernel \(K\) have finite second moments, with
  covariance matrix \(M_K=\int(s-m_K)(s-m_K)^\top K(s)\,ds\)
  (\(m_K=\int sK(s)\,ds\)), and let the spectral second moment
  \(\Lambda_X=\int\omega\omega^\top F_{XX}(d\omega)\) be finite (the
  spectral condition for mean-square differentiability of \(X\); Yaglom
  1987; Stein 1999). Then, as \(T\to0\), \[
  a(\delta_{K_T})=T^2\,\frac{\operatorname{tr}(M_K\Lambda_X)}{C_{XX}(0)}+o(T^2),
  \] and likewise for \(Y\). If both hold, write
  \(\alpha_X=\operatorname{tr}(M_K\Lambda_X)/C_{XX}(0)\) and
  \(\alpha_Y=\operatorname{tr}(M_K\Lambda_Y)/C_{YY}(0)\); the bound in
  (i) is then \[
  \mathcal{B}\big(a(\delta_{K_T}),b(\delta_{K_T});\rho\big)=T^2\Big(\sqrt{\alpha_X\alpha_Y}+\frac{|\rho|}{2}(\alpha_X+\alpha_Y)\Big)+o(T^2):
  \] halving the zone scale reduces the leading term to one quarter.
\item
  Rough fields: suppose that
  \(C_{XX}(0)-C_{XX}(h)=c\,|h|^\alpha\big(1+o(1)\big)\) as \(h\to0\),
  with \(c>0\) and \(0<\alpha<2\), and that \(K\) has bounded support.
  Then, as \(T\to0\), \[
  a(\delta_{K_T})=\frac{c\,\mathbb{E}|D|^\alpha}{C_{XX}(0)}\,T^\alpha\big(1+o(1)\big),
  \] where \(D=S-S'\) with \(S,S'\) independent with density \(K\). In
  particular, if \(q=1\), \(C_{XX}(h)=\sigma^2e^{-|h|/\ell}\) and \(K\)
  is the uniform density on \([0,1]\) (\(\alpha=1\),
  \(c=\sigma^2/\ell\), \(\mathbb{E}|D|=1/3\)), then
  \(a(\delta_{K_T})=T/(3\ell)+O(T^2)\).
\end{enumerate}

\emph{Proof.} (i) By Fubini's theorem, \[
\begin{aligned}\int K(s)\,\mathbb{E}\big[(X(s)-X_K)^2\big]ds&=\int K(s)\,\mathbb{E}X(s)^2\,ds-2\,\mathbb{E}\Big[X_K\int K(s)X(s)\,ds\Big]+\mathbb{E}X_K^2\\&=C_{XX}(0)-\mathbb{E}X_K^2 .\end{aligned}
\] Write \(g=g_\pi\); \(1-g\) is a non-negative bounded even function.
As in the proof of Proposition 9.7(b), Lemma 8.4(c) shows that the
matrix \[
\Sigma_W(\pi)=\begin{pmatrix}F_{XX}(1-g)&\Gamma(1-g)\\\Gamma(1-g)&F_{YY}(1-g)\end{pmatrix}
\] is positive semi-definite, and \(F_{XX}(1-g)=a\,C_{XX}(0)\),
\(F_{YY}(1-g)=b\,C_{YY}(0)\), \(F_{XX}(g)=(1-a)C_{XX}(0)\),
\(F_{YY}(g)=(1-b)C_{YY}(0)\). Since
\(\Gamma(\mathbb{R}^q)=\Gamma(g)+\Gamma(1-g)\), \[
\rho=\frac{\Gamma(g)+\Gamma(1-g)}{\sqrt{C_{XX}(0)C_{YY}(0)}}=\sqrt{(1-a)(1-b)}\;\rho(\pi)+\sqrt{ab}\;\rho_W(\pi),
\] where
\(\rho_W(\pi)=\Gamma(1-g)/\sqrt{F_{XX}(1-g)F_{YY}(1-g)}\in[-1,1]\)
(\(0\) if the denominator vanishes). This is exactly the identity of
Theorem 5.2; the proof of Theorem 7.1 uses only this identity and
\(|\rho_W|\le1\), so the same bound holds.

\begin{enumerate}
\def\labelenumi{(\roman{enumi})}
\setcounter{enumi}{1}
\item
  Let \(S,S'\) be independent with density \(K\) and \(D=S-S'\). Then
  \(\hat K(\omega)=\mathbb{E}e^{-iS\cdot\omega}\) and
  \(|\hat K(\omega)|^2=\mathbb{E}e^{-iD\cdot\omega}=\mathbb{E}\cos(D\cdot\omega)\)
  (the law of \(D\) is symmetric about the origin). By Lemma 8.5(d),
  \(|\hat K_T(\omega)|^2=\mathbb{E}\cos(TD\cdot\omega)\), so \[
  \frac{a(\delta_{K_T})\,C_{XX}(0)}{T^2}=\int\frac{1-|\hat K(T\omega)|^2}{T^2}\,F_{XX}(d\omega)=\int\mathbb{E}\Big[\frac{1-\cos(TD\cdot\omega)}{T^2}\Big]F_{XX}(d\omega).
  \] By \(0\le1-\cos x\le x^2/2\) and \((1-\cos x)/x^2\to1/2\) as
  \(x\to0\), the integrand converges pointwise to
  \(\tfrac12\mathbb{E}(D\cdot\omega)^2=\omega^\top M_K\,\omega\) (since
  \(\mathbb{E}DD^\top=2M_K\)) and is dominated by
  \(\tfrac12\mathbb{E}(D\cdot\omega)^2\le|\omega|^2\operatorname{tr}M_K\),
  which is \(F_{XX}\)-integrable because
  \(\int|\omega|^2\,dF_{XX}=\operatorname{tr}\Lambda_X<\infty\). Two
  applications of dominated convergence (first in \(D\), then in
  \(\omega\)) give the limit
  \(\int\omega^\top M_K\omega\,F_{XX}(d\omega)=\operatorname{tr}(M_K\Lambda_X)\).
  For the expansion of the bound: putting \(u=a\) and \(w=b\) in the
  proof of Corollary 7.2, the same estimate gives
  \(0\le\mathcal{B}(a,b;\rho)-\big[\sqrt{ab}+\tfrac{|\rho|}{2}(a+b)\big]\le5\max(a,b)^2\)
  for \(\max(a,b)\le1/2\), and \(\max(a,b)^2=O(T^4)\).
\item
  As in the proof of (ii), Fubini's theorem and Lemma 8.4(d) give \[
  a(\delta_{K_T})\,C_{XX}(0)=\mathbb{E}\int\big(1-\cos(TD\cdot\omega)\big)F_{XX}(d\omega)=\mathbb{E}\big[C_{XX}(0)-C_{XX}(TD)\big].
  \] Write \(C_{XX}(0)-C_{XX}(h)=c|h|^\alpha\big(1+r(h)\big)\) with
  \(r(h)\to0\) as \(h\to0\) (\(r(0):=0\)). If the support of \(K\) lies
  in a ball of radius \(R\), then \(|D|\le2R\), so
  \(T^{-\alpha}\big(C_{XX}(0)-C_{XX}(TD)\big)=c|D|^\alpha\big(1+r(TD)\big)\),
  which tends to \(c|D|^\alpha\) and, since \(r\) is bounded near \(0\),
  is bounded for small \(T\). By dominated convergence
  \(T^{-\alpha}a(\delta_{K_T})\,C_{XX}(0)\to c\,\mathbb{E}|D|^\alpha\),
  and \(\mathbb{E}|D|^\alpha>0\) because \(D\) has a density. For the
  exponential covariance, Lemma 8.6 gives explicitly \[
  \begin{aligned}\operatorname{Var}(X_{K_T})&=\frac{1}{T^2}\int_0^T\!\!\int_0^T\sigma^2e^{-|s-t|/\ell}\,ds\,dt=\frac{2\sigma^2}{T^2}\int_0^T(T-h)\,e^{-h/\ell}\,dh\\&=2\sigma^2\,\frac{x-1+e^{-x}}{x^2},\qquad x=\frac{T}{\ell}.\end{aligned}
  \] From \(e^{-x}=1-x+\tfrac{x^2}{2}-\tfrac{x^3}{6}+O(x^4)\),
  \(2(x-1+e^{-x})/x^2=1-x/3+O(x^2)\), so \(a(\delta_{K_T})=x/3+O(x^2)\).
  \(\blacksquare\)
\end{enumerate}

Proposition 8.16 shows that ``halving the zone scale quarters the
leading error term'' requires sufficiently smooth fields: it holds for
Lipschitz fields (Corollary 7.3) and for mean-square differentiable
stationary fields (part (ii)). If \(X\) is not mean-square
differentiable (\(\operatorname{tr}\Lambda_X=\infty\)) and \(K\) has
finite second moments, then \(a(\delta_{K_T})/T^2\to\infty\): in the
display of the proof of (ii), Fatou's lemma bounds the lower limit from
below by \(\int\omega^\top M_K\,\omega\,F_{XX}(d\omega)=\infty\), since
\(M_K\) is positive definite (\(K\) is a density). Under the condition
of (iii), halving the zone scale multiplies the leading term of \(a\) by
\(2^{-\alpha}\); the exponential covariance is the case \(\alpha=1\), in
which halving the zone scale only halves the leading term. The spectral
density of the exponential covariance (\(q=1\)) is proportional to
\(1/(1+\ell^2\omega^2)\), whose second moment is infinite, consistent
with the condition in (ii).

\textbf{Example (a first-order deviation).} In \(q=1\) let \(A\) and
\(B\) be independent zero-mean stationary fields, \(A\) with covariance
\(\sigma_A^2e^{-|h|/\ell}\) and \(B\) mean-square differentiable with
variance \(\sigma_B^2>0\), and put \(X=A\), \(Y=A+B\), so that
\(\rho=\sigma_A/\sqrt{\sigma_A^2+\sigma_B^2}\). Let \(K\) be the uniform
density on \([0,1]\) and \(g=|\hat K_T|^2\). Then \(F_{XX}=F_A\),
\(F_{YY}=F_A+F_B\) and \(\Gamma=F_A\), so
\(\rho(\delta_{K_T})=\big(1+F_B(g)/F_A(g)\big)^{-1/2}\); by (iii) and
(ii), \(F_A(g)=\sigma_A^2\big(1-T/(3\ell)+O(T^2)\big)\) and
\(F_B(g)=\sigma_B^2\big(1+O(T^2)\big)\). Since
\(\sigma_B^2/(\sigma_A^2+\sigma_B^2)=1-\rho^2\), expanding gives \[
\rho(\delta_{K_T})-\rho=-\frac{\rho(1-\rho^2)}{6\ell}\,T+O(T^2).
\] For rough fields, then, the deviation itself, and not only the bound,
can be of first order in the zone scale.

\textbf{Example (common signal plus independent noise).} Let \(X=S+N_X\)
and \(Y=S+N_Y\), where \(S\), \(N_X\), \(N_Y\) are independent
stationary fields with spectral densities \(f_S\), \(f_{N_X}\),
\(f_{N_Y}\) (bounded and continuous at the origin). Then
\(f_{XX}=f_S+f_{N_X}\), \(f_{YY}=f_S+f_{N_Y}\), \(c_{XY}=f_S\), and \[
\operatorname{coh}(0)=\frac{f_S(0)}{\sqrt{\big(f_S(0)+f_{N_X}(0)\big)\big(f_S(0)+f_{N_Y}(0)\big)}},\qquad
\rho=\frac{\sigma_S^2}{\sqrt{(\sigma_S^2+\sigma_{N_X}^2)(\sigma_S^2+\sigma_{N_Y}^2)}} .
\] If the common signal \(S\) is long-range and the noises are
short-range, then \(f_{N_X}(0)\) and \(f_{N_Y}(0)\) are small relative
to \(f_S(0)\) and \(\operatorname{coh}(0)\approx1\), while with large
noise variances the fine-scale correlation \(\rho\) can be small: this
is the ``aggregation increases correlation'' observed by Gehlke \& Biehl
(1934).

Figure 5 shows an example in which the correlation changes sign with
scale. The two fields share a long-range component \(S\), carry
mid-range components of opposite signs \(\pm M\), and have independent
short-range noise: \(X=S+M+N_1\), \(Y=S-M+N_2\). If the four components
are independent with spectral densities \(f_S,f_M,f_1,f_2\), then
\(f_{XX}=f_S+f_M+f_1\), \(f_{YY}=f_S+f_M+f_2\), \(c_{XY}=f_S-f_M\), and
\[
\operatorname{coh}(\omega)=\frac{f_S(\omega)-f_M(\omega)}{\sqrt{\big(f_S+f_M+f_1\big)(\omega)\,\big(f_S+f_M+f_2\big)(\omega)}} :
\] low frequencies are dominated by \(S\) and the coherence is positive;
intermediate frequencies are dominated by \(M\) and the coherence is
negative; at high frequencies there are only independent noises and the
coherence is close to zero. The fine-scale correlation is \(-0.40\); the
zonal correlation changes sign for block sides between 8 and 16 pixels;
as the blocks grow further it tends to \(\operatorname{coh}(0)=+0.91\).
The zoning here consists of regular blocks and does not depend on the
data; the change of correlation with scale comes entirely from the
different coherences of the two fields at different frequencies, and the
zonal correlations at all scales stay within the coherence bounds
\([-0.98,\,0.91]\) (Theorem 8.9).

\begin{figure}
\centering
\includegraphics[width=1\textwidth,height=\textheight]{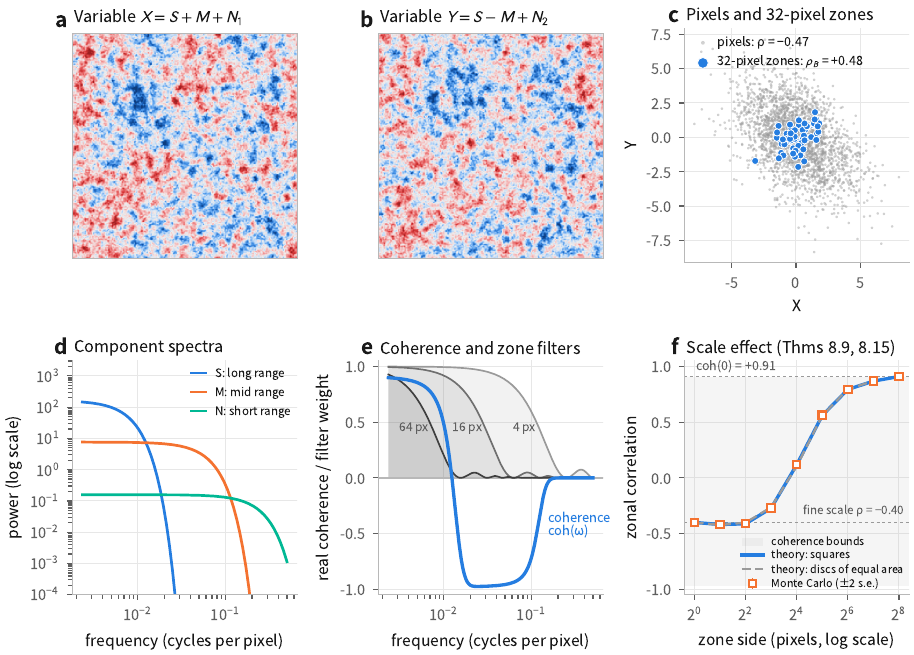}
\caption{Change of the zonal correlation with scale under
data-independent zonings. (a), (b) One realization of Gaussian fields on
a \(256\times256\) torus: \(X=S+M+N_1\), \(Y=S-M+N_2\), where \(S\) is a
shared long-range component, \(M\) a mid-range component entering with
opposite signs, and \(N_1\), \(N_2\) independent short-range noises. (c)
Pixel-level scatter (grey; sample correlation of this realization
\(-0.47\), theoretical value \(-0.40\)) and block means for blocks of
side 32 pixels (blue; \(\rho_B=+0.48\)). (d) Power spectra of the three
components. (e) Blue: real coherence \(\operatorname{coh}(\omega)\),
close to \(+0.91\) at the lowest frequencies, negative at intermediate
frequencies (minimum about \(-0.98\)), and zero at high frequencies
(where there are only independent noises); grey curves are the filters
\(|\hat K|^2\) of square blocks of side 64, 16 and 4 pixels: the larger
the block, the more only low frequencies are kept. (f) Zonal correlation
against block side: blue line, theoretical value for square blocks
(Theorem 8.9); grey dashed line, theoretical value for discs of equal
area; orange squares, mean of 1500 simulations (\(\pm2\) standard
errors); grey band, coherence bounds \([-0.98,\ 0.91]\). The curve
starts at the fine-scale correlation \(-0.40\), changes sign between 8
and 16 pixels, and tends to \(\operatorname{coh}(0)=+0.91\) (Theorem
8.15); the square and disc curves coincide at both ends and differ by at
most 0.013 in between.}
\end{figure}

\hypertarget{nested-hierarchies-the-martingale-scale-spectrum}{%
\section{9 Nested hierarchies: the martingale scale
spectrum}\label{nested-hierarchies-the-martingale-scale-spectrum}}

Aggregating level by level along a sequence of nested zonings decomposes
the fine-scale covariance orthogonally across scales (Theorem 9.2). Each
scale has its own level coherence \(\kappa_l\); the zonal correlation at
any level depends only on that level and the coarser ones, and is a
weighted average of these \(\kappa_l\) times a factor at most 1 (Theorem
9.3). This yields a sufficient condition for ``coarser aggregation,
higher correlation'' (Theorem 9.4) and a counterexample when it fails
(Example 9.5), as well as the fixed points of the coarsening map
(Proposition 9.6); for dyadic grid hierarchies of stationary fields all
these quantities can be computed from the spectrum (Proposition 9.7).
Sections 9.1--9.3 work on the support measure space
\((\Omega,\mathcal{F},\mu)\) of Section 3 and need no atomlessness
assumption; Section 9.4 works in the random-field setting of Section
8.1.

\hypertarget{orthogonal-decomposition}{%
\subsection{9.1 Orthogonal
decomposition}\label{orthogonal-decomposition}}

\textbf{Setting.} Let \(X,Y\in L^2(\mu)\) and let a finite nested
hierarchy \[
\mathcal{F}\supseteq\mathcal{G}_0\supseteq\mathcal{G}_1\supseteq\dots\supseteq\mathcal{G}_L=\{\varnothing,\Omega\}
\] be given, \(\mathcal{G}_0\) being the finest observed level
(\(\mathcal{G}_0=\mathcal{F}\) is the fine scale itself). The most
common case is \(\mathcal{G}_l=\sigma(\mathcal{Z}_l)\) for a sequence of
successively coarser zonings \(\mathcal{Z}_l\) (each zone of
\(\mathcal{Z}_{l+1}\) is a union of zones of \(\mathcal{Z}_l\)). Let
\(X^{(l)}=\mathbb{E}[X\mid\mathcal{G}_l]\), which by Lemma 3.1 is the
zone-mean function at level \(l\), and write the increments as
\(D_lX=X^{(l)}-X^{(l+1)}\) (\(l=0,\dots,L-1\)). Define the \textbf{scale
spectrum} and the \textbf{level coherences} \[
v_l(X)=\mathbb{E}\big[(D_lX)^2\big],\qquad c_l=\mathbb{E}\big[D_lX\cdot D_lY\big],\qquad
\kappa_l=\frac{c_l}{\sqrt{v_l(X)\,v_l(Y)}},
\] with \(\kappa_l=0\) when the denominator vanishes. Let
\(\rho_l=\operatorname{Corr}\big(X^{(l)},Y^{(l)}\big)\) (when defined):
it is the zonal correlation at level \(l\), weighted by zone weights;
when \(\mathcal{G}_0=\mathcal{F}\), \(\rho_0\) is the fine-scale
correlation.

\(v_l(X)\) is the variance component between levels \(l\) and \(l+1\):
it measures how much of the variation of \(X\) occurs ``within one zone
of level \(l+1\), between different zones of level \(l\)''. \(c_l\) is
the covariance component at the same level, and \(\kappa_l\) is the
correlation specific to that level. This parallels the variance
components of hierarchical (multilevel) models, except that the
decomposition below is an exact identity and does not depend on any
distributional assumption.

\textbf{Lemma 9.1 (orthogonality across scales).} For every \(l\),
\(D_lX\in L^2(\mu)\) and \(\mathbb{E}D_lX=0\); for \(l\neq m\),
\(\mathbb{E}[D_lX\cdot D_mY]=0\).

\emph{Proof.} By Lemma 3.5, \(X^{(l)}\in L^2\), so \(D_lX\in L^2\); by
the tower property (Lemma 3.2), \(\mathbb{E}X^{(l)}=\mathbb{E}X\) for
all \(l\), so \(\mathbb{E}D_lX=0\). Suppose \(l<m\).
\(D_mY=Y^{(m)}-Y^{(m+1)}\) is \(\mathcal{G}_m\)-measurable, and
\(m\ge l+1\) gives \(\mathcal{G}_m\subseteq\mathcal{G}_{l+1}\), so
\(D_mY\) is a \(\mathcal{G}_{l+1}\)-measurable \(L^2\) function. By the
tower property, \(\mathbb{E}[X^{(l)}\mid\mathcal{G}_{l+1}]=X^{(l+1)}\),
i.e.~\(D_lX=X^{(l)}-\mathbb{E}[X^{(l)}\mid\mathcal{G}_{l+1}]\). Lemma
3.5 applied to \(X^{(l)}\), the sub-\(\sigma\)-algebra
\(\mathcal{G}_{l+1}\) and \(W=D_mY\) gives
\(\mathbb{E}[D_lX\cdot D_mY]=0\). For \(l>m\) exchange the roles of
\(X\) and \(Y\). \(\blacksquare\)

\textbf{Theorem 9.2 (martingale scale-spectrum decomposition).} For
every \(0\le l'\le L-1\), \[
\begin{aligned}
&\operatorname{Cov}\big(X^{(l')},Y^{(l')}\big)=\sum_{l=l'}^{L-1}c_l,\qquad
\operatorname{Var}\big(X^{(l')}\big)=\sum_{l=l'}^{L-1}v_l(X),\\
&\rho_{l'}=\frac{\sum_{l\ge l'}c_l}{\sqrt{\sum_{l\ge l'}v_l(X)\,\sum_{l\ge l'}v_l(Y)}},
\end{aligned}
\] where \(\rho_{l'}\) is defined when the denominator is positive;
moreover \(|c_l|\le\sqrt{v_l(X)v_l(Y)}\). In particular, \(\rho_{l'}\)
depends only on level \(l'\) and the coarser levels: structure finer
than \(l'\) has no influence whatsoever on the zonal correlation at
level \(l'\).

\emph{Proof.}
\(X^{(L)}=\mathbb{E}[X\mid\{\varnothing,\Omega\}]=\mathbb{E}X\), so
telescoping gives \[
X^{(l')}-\mathbb{E}X=\sum_{l=l'}^{L-1}D_lX,\qquad Y^{(l')}-\mathbb{E}Y=\sum_{m=l'}^{L-1}D_mY .
\] Multiply the two and take expectations; by Lemma 9.1 the cross terms
with \(l\neq m\) have zero expectation, leaving
\(\sum_l\mathbb{E}[D_lX\,D_lY]=\sum_lc_l\). Taking \(Y=X\) gives the
variance formula. \(|c_l|\le\sqrt{v_l(X)v_l(Y)}\) is the Cauchy--Schwarz
inequality. \(\blacksquare\)

With \(L=2\), \(\mathcal{G}_0=\mathcal{F}\) and
\(\mathcal{G}_1=\sigma(\mathcal{Z})\), Theorem 9.2 is Lemma 5.1:
\(c_0=C_W\), \(c_1=C_B\), \(v_0(X)=W_X\), \(v_1(X)=B_X\).

\hypertarget{the-hierarchical-coherence-bound-and-the-direction-of-scale-effects}{%
\subsection{9.2 The hierarchical coherence bound and the direction of
scale
effects}\label{the-hierarchical-coherence-bound-and-the-direction-of-scale-effects}}

\textbf{Theorem 9.3 (hierarchical coherence bound and invariance).}

\begin{enumerate}
\def\labelenumi{(\roman{enumi})}
\item
  If \(\rho_{l'}\) is defined, then
  \(\rho_{l'}=\theta_{l'}\sum_{l\ge l'}\omega_l\kappa_l\), where \[
  \omega_l=\frac{\sqrt{v_l(X)v_l(Y)}}{\sum_{m\ge l'}\sqrt{v_m(X)v_m(Y)}},\qquad
  \theta_{l'}=\frac{\sum_{l\ge l'}\sqrt{v_l(X)v_l(Y)}}{\sqrt{\sum_{l\ge l'}v_l(X)\cdot\sum_{l\ge l'}v_l(Y)}}\in[0,1]
  \] (\(\omega_l\) are non-negative weights summing to 1; if
  \(\sum_{m\ge l'}\sqrt{v_m(X)v_m(Y)}=0\) then \(\rho_{l'}=0\)). Hence
  \[
  \min\Big(0,\min_{l\ge l'}\kappa_l\Big)\le\rho_{l'}\le\max\Big(0,\max_{l\ge l'}\kappa_l\Big),\qquad |\rho_{l'}|\le\max_{l\ge l'}|\kappa_l| .
  \]
\item
  If the scale spectra are proportional, i.e.~there is \(\lambda>0\)
  with \(v_l(Y)=\lambda v_l(X)\) for all \(l\), then \(\theta_{l'}=1\)
  and \[
  \rho_{l'}=\frac{\sum_{l\ge l'}v_l(X)\,\kappa_l}{\sum_{l\ge l'}v_l(X)}
  \] is the average of \(\{\kappa_l:l\ge l'\}\) weighted by \(v_l(X)\).
  In this case \(\rho_{l'}\) is the same at all levels where it is
  defined if and only if \(\kappa_l\) is the same at all levels with
  \(v_l(X)>0\).
\end{enumerate}

\emph{Proof.} (i) If \(v_l(X)v_l(Y)=0\), then
\(|c_l|\le\sqrt{v_l(X)v_l(Y)}\) gives \(c_l=0\), consistent with the
convention \(\kappa_l=0\); so \(c_l=\kappa_l\sqrt{v_l(X)v_l(Y)}\) for
all \(l\). Substituting into Theorem 9.2 gives
\(\rho_{l'}=\theta_{l'}\sum\omega_l\kappa_l\). By the Cauchy--Schwarz
inequality,
\(\sum_l\sqrt{v_l(X)}\sqrt{v_l(Y)}\le\big(\sum_lv_l(X)\big)^{1/2}\big(\sum_lv_l(Y)\big)^{1/2}\),
so \(\theta_{l'}\in[0,1]\). \(\sum\omega_l\kappa_l\) lies between
\(\min\kappa_l\) and \(\max\kappa_l\), and multiplying by
\(\theta_{l'}\in[0,1]\) only shrinks it towards \(0\).

\begin{enumerate}
\def\labelenumi{(\roman{enumi})}
\setcounter{enumi}{1}
\tightlist
\item
  Under proportionality, \(\sqrt{v_l(X)v_l(Y)}=\sqrt\lambda\,v_l(X)\)
  and \(\sum_lv_l(Y)=\lambda\sum_lv_l(X)\), so \(\theta_{l'}=1\) and
  \(\omega_l=v_l(X)/\sum_{m\ge l'}v_m(X)\). Write
  \(S_{l'}=\sum_{i\ge l'}v_i(X)\) and
  \(T_{l'}=\sum_{i\ge l'}v_i(X)\kappa_i\) (\(S_L=T_L=0\)). Then
  \(\rho_{l'}\) is defined if and only if \(S_{l'}>0\), in which case
  \(\rho_{l'}=T_{l'}/S_{l'}\), and
  \(T_{l'}=v_{l'}(X)\kappa_{l'}+T_{l'+1}\). ``\(\Leftarrow\)'': if
  \(\kappa_l=\kappa\) at all levels with \(v_l(X)>0\), then
  \(T_{l'}=\kappa S_{l'}\) and \(\rho_{l'}=\kappa\).
  ``\(\Rightarrow\)'': if \(\rho_{l'}=\rho_\ast\) at all levels with
  \(S_{l'}>0\), then \(T_{l'}=\rho_\ast S_{l'}\) for all \(l'\) (both
  sides vanish when \(S_{l'}=0\)). At levels with \(v_{l'}(X)>0\),
  \(v_{l'}(X)\kappa_{l'}=T_{l'}-T_{l'+1}=\rho_\ast(S_{l'}-S_{l'+1})=\rho_\ast v_{l'}(X)\),
  i.e.~\(\kappa_{l'}=\rho_\ast\). \(\blacksquare\)
\end{enumerate}

Theorem 9.3(ii) is the hierarchical counterpart of Theorem 8.11: with
proportional scale spectra, the zonal correlation is invariant across
levels if and only if the level coherences are everywhere the same.

\textbf{Theorem 9.4 (monotone scale effects).} Suppose the scale spectra
are proportional and \(v_l(X)>0\) (\(l=0,\dots,L-1\)). If \(\kappa_l\)
is non-decreasing in \(l\) (coarser scales are more coherent), then
\(\rho_{l'}\) is non-decreasing in \(l'\), i.e.~the coarser the
aggregation, the larger the zonal correlation; if \(\kappa\) is strictly
increasing, so is \(\rho\). If \(\kappa\) is non-increasing, the
conclusion is reversed.

\emph{Proof.} Use the notation of the proof of Theorem 9.3(ii). For
\(l'\le L-2\), \[
\rho_{l'}=\frac{v_{l'}(X)\,\kappa_{l'}+S_{l'+1}\,\rho_{l'+1}}{v_{l'}(X)+S_{l'+1}}
\] is a convex combination of \(\kappa_{l'}\) and \(\rho_{l'+1}\) with
positive weights. \(\rho_{l'+1}\) is a weighted average of
\(\{\kappa_l:l\ge l'+1\}\), so if \(\kappa\) is non-decreasing then
\(\rho_{l'+1}\ge\kappa_{l'+1}\ge\kappa_{l'}\), hence
\(\rho_{l'}\le\rho_{l'+1}\). If \(\kappa\) is strictly increasing,
\(\rho_{l'+1}\ge\kappa_{l'+1}>\kappa_{l'}\), and the weight of
\(\kappa_{l'}\) is positive, so \(\rho_{l'}<\rho_{l'+1}\). For
non-increasing \(\kappa\) reverse all inequalities. \(\blacksquare\)

\textbf{Remark (one aggregation step).} If the scale spectra are
proportional and \(v_l(X)>0\) for all \(l\) (no monotonicity of
\(\kappa\) is needed), the display in the proof of Theorem 9.4 gives,
for \(l'\le L-2\), \[
\rho_{l'}-\rho_{l'+1}=\frac{v_{l'}(X)}{S_{l'}}\big(\kappa_{l'}-\rho_{l'+1}\big),
\] so aggregating from level \(l'\) to level \(l'+1\) raises the zonal
correlation if and only if \(\kappa_{l'}<\rho_{l'+1}\), and leaves it
unchanged if and only if \(\kappa_{l'}=\rho_{l'+1}\). Theorem 9.4 is a
simple sufficient condition for the correlation to rise (or not to fall)
at every step.

The typical positive case is \(X=S+N_X\), \(Y=S+N_Y\), where the common
component \(S\) varies mainly at coarse scales and the independent
noises \(N_X\), \(N_Y\) only at fine scales. Then the fine-level
\(\kappa_l\) are close to \(0\), the coarse-level \(\kappa_l\) close to
\(1\), and aggregation increases the correlation (Figure 6). The
proportionality condition cannot be dropped:

\textbf{Example 9.5 (level coherence rises while the zonal correlation
falls).} Take \(L=2\) with scale spectra and level coherences \[
v_0(X)=1,\quad v_0(Y)=0.01,\quad v_1(X)=0.01,\quad v_1(Y)=1,\qquad \kappa_0=-0.9<\kappa_1=-0.5 .
\] Then \(c_0=-0.09\), \(c_1=-0.05\), and by Theorem 9.2 \[
\rho_1=\kappa_1=-0.5,\qquad \rho_0=\frac{-0.09-0.05}{\sqrt{1.01\times1.01}}=-\frac{0.14}{1.01}\approx-0.1386>\rho_1 .
\] The level coherence increases with scale, yet the zonal correlation
decreases with aggregation.

\emph{An explicit realization.} Let \(\Omega=\{(j,m):j,m\in\{1,2,3\}\}\)
with \(\mu\) uniform, \(\mathcal{G}_0=\mathcal{F}\) all subsets,
\(\mathcal{G}_1\) generated by the three top-level zones
\(A_j=\{(j,1),(j,2),(j,3)\}\), and \(\mathcal{G}_2\) trivial (the top
level needs at least 3 zones, since two weighted zone means always lie
on a line, so a two-zone correlation can only be \(\pm1\)). Take the
vectors \[
a=\sqrt{3/2}\,(1,-1,0),\qquad b=\sqrt{1/2}\,(1,1,-2)
\] in \(\mathbb{R}^3\); their components sum to zero,
\(\tfrac13\sum_ia_i^2=\tfrac13\sum_ib_i^2=1\) and
\(\tfrac13\sum_ia_ib_i=0\). For \(\kappa\in[-1,1]\) let
\(b_\kappa=\kappa a+\sqrt{1-\kappa^2}\,b\), so that
\(\tfrac13\sum_ib_{\kappa,i}^2=1\) and
\(\tfrac13\sum_ia_ib_{\kappa,i}=\kappa\). Define \[
X(j,m)=0.1\,a_j+a_m,\qquad Y(j,m)=b_{-0.5,j}+0.1\,b_{-0.9,m}.
\] Within the top-level zone \(A_j\) the second terms average to zero
over \(m\), so \(X^{(1)}(j,m)=0.1a_j\), \(Y^{(1)}(j,m)=b_{-0.5,j}\),
\(X^{(2)}=Y^{(2)}=0\), \(D_1X=0.1a_j\), \(D_0X=a_m\),
\(D_1Y=b_{-0.5,j}\) and \(D_0Y=0.1b_{-0.9,m}\). Direct computation gives
\(v_1(X)=0.01\), \(v_1(Y)=1\), \(c_1=0.1\times(-0.5)=-0.05\),
\(v_0(X)=1\), \(v_0(Y)=0.01\), \(c_0=0.1\times(-0.9)=-0.09\), which is
the scale spectrum above. The fine-scale correlation is
\(\rho_0\approx-0.139\), and the zonal correlation over the three
top-level zones is \(\rho_1=-0.5\).

\hypertarget{fixed-points-of-the-coarsening-map}{%
\subsection{9.3 Fixed points of the coarsening
map}\label{fixed-points-of-the-coarsening-map}}

\textbf{Proposition 9.6 (fixed points of the scale-spectrum flow are
geometric spectra).} Let \((v_l)_{l\ge0}\) be a sequence of positive
numbers with \(\sum_lv_l<\infty\) (an infinite hierarchy, for example
the dyadic aggregation of a stationary field in Proposition 9.7). Let \[
\hat s_{l'}(j)=\frac{v_{l'+j}}{\sum_{i\ge l'}v_i}\qquad(j\ge0)
\] be the normalized scale spectrum seen from level \(l'\). The
coarsening map \(\hat s_{l'}\mapsto\hat s_{l'+1}\) leaves \(\hat s\)
unchanged for every \(l'\) if and only if \(v_l=v_0r^l\) for some
\(r\in(0,1)\).

\emph{Proof.} ``\(\Leftarrow\)'':
\(\hat s_{l'}(j)=r^{l'+j}\big/\sum_{i\ge l'}r^i=r^j(1-r)\), independent
of \(l'\). ``\(\Rightarrow\)'': write \(S_l=\sum_{i\ge l}v_i\).
\(\hat s_{l'}(j)=\hat s_{l'+1}(j)\) for all \(j\) means
\(v_{l'+1+j}/v_{l'+j}=S_{l'+1}/S_{l'}=:r_{l'}\), independent of \(j\).
Taking \(j=1\) gives \(v_{l'+2}/v_{l'+1}=r_{l'}\); taking \(j=0\) for
\(l'+1\), the same ratio equals \(r_{l'+1}\). So \(r_{l'}=r_{l'+1}\) for
all \(l'\); call it \(r\). Then \(v_l=v_0r^l\), and \(\sum_lv_l<\infty\)
with \(v_l>0\) forces \(0<r<1\). \(\blacksquare\)

On the linear scale \(2^l\), a geometric spectrum is a power law
\(v_l\propto(2^l)^{\log_2r}\). So the fixed points of ``coarsen and
renormalize'' are exactly the scale spectra without a characteristic
scale; this is the same statement in the language of the renormalization
group.

\hypertarget{dyadic-hierarchies-of-stationary-fields}{%
\subsection{9.4 Dyadic hierarchies of stationary
fields}\label{dyadic-hierarchies-of-stationary-fields}}

\textbf{Proposition 9.7 (dyadic hierarchies of stationary fields).} Let
\((X,Y)\) satisfy the setting of Section 8.1, with spectral matrix
without an atom at the origin: \(F_{XX}(\{0\})=F_{YY}(\{0\})=0\). Fix
\(h_0>0\), let \(K_l\) be the uniform density on the cube
\([0,2^lh_0)^q\) (\(l=0,1,2,\dots\)), and let
\(\Delta_l=|\hat K_l|^2-|\hat K_{l+1}|^2\). Then:

\begin{enumerate}
\def\labelenumi{(\alph{enumi})}
\item
  \(\Delta_l\ge0\), and for every \(l'\) \[
  \operatorname{Cov}\big(X_{K_{l'}},Y_{K_{l'}}\big)=\sum_{l\ge l'}\int\Delta_l\,d\Gamma,\qquad
  \operatorname{Var}\big(X_{K_{l'}}\big)=\sum_{l\ge l'}\int\Delta_l\,dF_{XX},
  \] the series converging absolutely.
\item
  Let \[
  v_l(X)=\int\Delta_l\,dF_{XX},\qquad v_l(Y)=\int\Delta_l\,dF_{YY},\qquad c_l=\int\Delta_l\,d\Gamma .
  \] Then \(|c_l|\le\sqrt{v_l(X)v_l(Y)}\). If \(\tau\) is uniformly
  distributed on \(\{0,2^lh_0\}^q\) and independent of the fields, then
  \[
  \mathbb{E}\Big[\big(X_{K_l(\cdot-\tau)}-X_{K_{l+1}}\big)^2\Big]=v_l(X),\qquad
  \mathbb{E}\Big[\big(X_{K_l(\cdot-\tau)}-X_{K_{l+1}}\big)\big(Y_{K_l(\cdot-\tau)}-Y_{K_{l+1}}\big)\Big]=c_l ,
  \] i.e.~\(v_l\) and \(c_l\) are the design-level variance and
  covariance of the increments between adjacent levels.
\end{enumerate}

Hence \((v_l(X),v_l(Y),c_l)_{l\ge0}\) is an infinite-level scale
spectrum, and Theorem 9.3 holds for it with the minima and maxima over
\(l\ge l'\) replaced by infima and suprema, while Theorem 9.4 and
Proposition 9.6 hold for it provided \(v_l(X)>0\) for all \(l\) (their
proofs use only the formulas of Theorem 9.2 and
\(|c_l|\le\sqrt{v_l(X)v_l(Y)}\), and the infinite sums converge
absolutely); the level coherence \(\kappa_l\) is the coherence under the
band-pass filter \(\Delta_l\).

\emph{Proof.} Write \(a=2^lh_0\). A half-open cube of side \(2a\) is the
disjoint union of \(2^q\) half-open cubes of side \(a\),
\([0,2a)^q=\bigsqcup_{\tau\in\{0,a\}^q}\big(\tau+[0,a)^q\big)\), so \[
K_{l+1}=2^{-q}\sum_{\tau\in\{0,a\}^q}K_l(\cdot-\tau).\tag{9.1}
\] By Lemma 8.5(d), \[
\hat K_{l+1}(\omega)=\hat K_l(\omega)\prod_{i=1}^{q}\frac{1+e^{-ia\omega_i}}{2},\qquad |\hat K_{l+1}(\omega)|=|\hat K_l(\omega)|\prod_{i=1}^{q}\Big|\cos\frac{a\omega_i}{2}\Big|\le|\hat K_l(\omega)|,
\] so \(\Delta_l\ge0\). Telescoping gives
\(|\hat K_{l'}|^2=\sum_{l=l'}^{L-1}\Delta_l+|\hat K_L|^2\) (\(L>l'\)).
The uniform density on the cube \([0,T)^q\) (\(T=2^Lh_0\)) satisfies \[
\hat K_L(\omega)=\prod_{i=1}^q\frac{1-e^{-iT\omega_i}}{iT\omega_i},\qquad |\hat K_L(\omega)|=\prod_{i=1}^q\Big|\frac{\sin(T\omega_i/2)}{T\omega_i/2}\Big|
\] (factors with \(\omega_i=0\) equal to \(1\)). For \(\omega\neq0\)
some \(\omega_i\neq0\), and the corresponding factor is at most
\(2/(T|\omega_i|)\to0\) as \(L\to\infty\); at \(\omega=0\) the value is
\(1\). By \(|\hat K_L|\le1\), dominated convergence and
\(F_{XX}(\{0\})=0\), \(\int|\hat K_L|^2\,dF_{XX}\to0\); likewise for
\(F_{YY}\), and \(|\Gamma|\le\Lambda/2\) (Lemma 8.8) gives
\(|\Gamma|(\{0\})=0\), so \(\int|\hat K_L|^2\,d\Gamma\to0\). By Lemma
8.6 the left-hand sides equal \(\int|\hat K_{l'}|^2\,d\Gamma\) (or
\(dF_{XX}\)), and letting \(L\to\infty\) gives (a). The terms
\(\int\Delta_l\,dF_{XX}\ge0\) have bounded partial sums, so the series
converges; the series for \(\Gamma\) converges absolutely by
\(|c_l|\le(v_l(X)+v_l(Y))/2\) (see below).

\begin{enumerate}
\def\labelenumi{(\alph{enumi})}
\setcounter{enumi}{1}
\tightlist
\item
  For every Borel set \(A\), by Lemma 8.4(c) the matrix
  \(\begin{pmatrix}F_{XX}(A)&\Gamma(A)\\\Gamma(A)&F_{YY}(A)\end{pmatrix}\)
  is positive semi-definite, so for real \(\alpha,\beta\) the measure
  \(\alpha^2F_{XX}+2\alpha\beta\Gamma+\beta^2F_{YY}\) is positive, and
  its integral against the non-negative function \(\Delta_l\) is
  \(\alpha^2v_l(X)+2\alpha\beta c_l+\beta^2v_l(Y)\ge0\). This quadratic
  form in \((\alpha,\beta)\) is positive semi-definite, so
  \(c_l^2\le v_l(X)v_l(Y)\); in particular
  \(|c_l|\le(v_l(X)+v_l(Y))/2\).
\end{enumerate}

Write \(X_\tau=X_{K_l(\cdot-\tau)}\), and \(Y_\tau\) likewise. By (9.1)
and linearity of block means in the kernel,
\(2^{-q}\sum_\tau X_\tau=X_{K_{l+1}}\). By independence and Fubini's
theorem, averaging over \(\tau\), \[
\begin{aligned}
\mathbb{E}\big[(X_\tau-X_{K_{l+1}})(Y_\tau-Y_{K_{l+1}})\big]
={}&2^{-q}\sum_\tau\mathbb{E}[X_\tau Y_\tau]-\mathbb{E}\Big[\Big(2^{-q}\sum_\tau X_\tau\Big)Y_{K_{l+1}}\Big]\\
&-\mathbb{E}\Big[X_{K_{l+1}}\Big(2^{-q}\sum_\tau Y_\tau\Big)\Big]+\mathbb{E}\big[X_{K_{l+1}}Y_{K_{l+1}}\big].
\end{aligned}
\] By translation invariance in Lemma 8.6 the first term equals
\(\int|\hat K_l|^2\,d\Gamma\); the last three are
\(\pm\int|\hat K_{l+1}|^2\,d\Gamma\) with signs \(-,-,+\). So the
left-hand side is \(\int\Delta_l\,d\Gamma=c_l\). Taking \(Y=X\) gives
the first identity. \(\blacksquare\)

Cube kernels belong to \(L^2\), so by Theorem 8.15(ii): if the spectral
matrix has bounded densities continuous at the origin with
\(f_{XX}(0)f_{YY}(0)>0\), then
\(\rho(\delta_{K_{l'}})\to\operatorname{coh}(0)\) as \(l'\to\infty\).
The macro limit is determined by the coherence at the coarsest scales
alone.

\begin{figure}
\centering
\includegraphics[width=1\textwidth,height=\textheight]{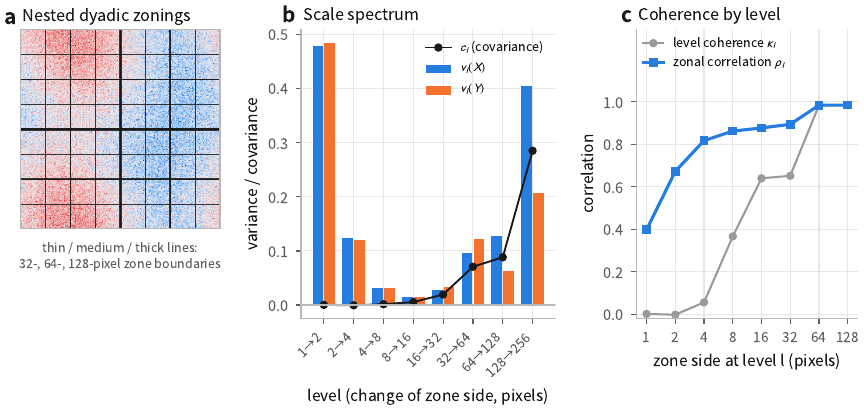}
\caption{Scale spectrum on a nested hierarchy. (a) The field \(X\) on a
\(256\times256\) grid with nested dyadic block zonings (levels doubling
from 1 to 256 pixels; thin, medium and thick lines are the boundaries of
32-, 64- and 128-pixel zones); \(X\) and \(Y\) share a large-scale
component and each carries independent pixel noise. (b) Scale spectrum:
variance components \(v_l(X)\), \(v_l(Y)\) (bars) and covariance
component \(c_l\) (black line) for each level (zone side doubling from
\(2^l\) to \(2^{l+1}\) pixels). The finest levels contain almost only
independent noise, \(c_l\approx0\); the common component is concentrated
at the coarsest levels. (c) Level coherence \(\kappa_l\) (grey) and the
zonal correlation \(\rho_l\) at the level with zone side \(2^l\) pixels
(blue). By Theorems 9.2 and 9.3, \(\rho_l\) equals a factor
\(\theta_l\le1\) times a weighted average of \(\kappa_m\) (\(m\ge l\));
since \(\kappa\) is close to 0 at fine scales and close to 1 at coarse
scales, \(\rho_l\) rises from 0.40 to 0.99 with aggregation. Theorem 9.4
gives a sufficient condition for such a monotone rise.}
\end{figure}

\hypertarget{identification-intervals-from-zonal-data-alone}{%
\section{10 Identification intervals from zonal data
alone}\label{identification-intervals-from-zonal-data-alone}}

This section considers the case in which only the zone weights and zone
means are observed and fine-scale data are unavailable. For variables
with values in \([0,1]\) (indicators or proportions), the fine-scale
covariance lies in an interval \([L,U]\) that can be computed from the
zonal data, and when the support measure is atomless every value in it
is attained (Theorem 10.1); the interval has zero width if and only if
in every zone at least one of the two variables has zone mean \(0\) or
\(1\) (Corollary 10.2). In the binary case these are the Duncan \& Davis
(1953) bounds, that is, the Fréchet (1951) bounds applied zone by zone.

\hypertarget{sharp-bounds}{%
\subsection{10.1 Sharp bounds}\label{sharp-bounds}}

\textbf{Theorem 10.1 (sharp identification bounds for the fine-scale
covariance).} Let \(\mathcal{Z}=\{A_1,\dots,A_k\}\) be a \(k\)-zone
zoning with zone weights \(p_j\) and known zone means
\(a_j,b_j\in[0,1]\). Let \(\bar a=\sum_jp_ja_j\),
\(\bar b=\sum_jp_jb_j\), \(C_B=\sum_jp_j(a_j-\bar a)(b_j-\bar b)\), and
\[
L=C_B+\sum_j p_j\big(\max(0,\,a_j+b_j-1)-a_jb_j\big),\qquad
U=C_B+\sum_j p_j\big(\min(a_j,b_j)-a_jb_j\big).
\]

\begin{enumerate}
\def\labelenumi{(\roman{enumi})}
\item
  For every support measure space, every zoning with zone weights
  \(p_j\), and all variables \(X,Y\) with values in \([0,1]\) and zone
  means \((a_j,b_j)\), \(L\le\operatorname{Cov}(X,Y)\le U\).
\item
  Conversely, if \(\mu\) is atomless, then for every zoning with zone
  weights \(p_j\) and every \(c\in[L,U]\) there are variables \(X,Y\)
  taking only the values \(0\) and \(1\), with zone means \((a_j,b_j)\)
  and \(\operatorname{Cov}(X,Y)=c\).
\end{enumerate}

Hence, if \(\mu\) is atomless, then over all fine-scale configurations
compatible with the zonal data \(\{(p_j,a_j,b_j)\}\), the set of values
of the fine-scale covariance is exactly \([L,U]\).

\emph{Proof.} (i) Write \(\tilde X=\mathbb{E}[X\mid\mathcal{G}]\) and
\(\tilde Y=\mathbb{E}[Y\mid\mathcal{G}]\) with
\(\mathcal{G}=\sigma(\mathcal{Z})\). By Lemma 5.1,
\(\operatorname{Cov}(X,Y)=C_B+C_W\), where
\(C_B=\operatorname{Cov}(\tilde X,\tilde Y)\) is the quantity computed
above from the zonal data and, by Lemma 3.1(c), \[
C_W=\mathbb{E}\big[(X-\tilde X)(Y-\tilde Y)\big]=\sum_jp_j\,\mathbb{E}\big[(X-a_j)(Y-b_j)\mid A_j\big]=\sum_jp_j\big(\mathbb{E}[XY\mid A_j]-a_jb_j\big).
\] On \([0,1]\), \(XY\le X\) and \(XY\le Y\), so
\(\mathbb{E}[XY\mid A_j]\le\min(a_j,b_j)\); \((1-X)(1-Y)\ge0\) gives
\(XY\ge X+Y-1\), which together with \(XY\ge0\) gives
\(\mathbb{E}[XY\mid A_j]\ge\max(0,\,a_j+b_j-1)\) (Fréchet 1951).
Substituting gives \(L\le\operatorname{Cov}(X,Y)\le U\).

\begin{enumerate}
\def\labelenumi{(\roman{enumi})}
\setcounter{enumi}{1}
\tightlist
\item
  Fix a zone \(A=A_j\) and write \(a=a_j\), \(b=b_j\), \(p=p_j\). By
  Lemma 4.5(b), the restriction of \(\mu\) to \(A\) is atomless; below
  we repeatedly use Theorem 4.3 (Sierpiński) to pick subsets of
  prescribed measure inside measurable subsets of \(A\). For
  \(t\in[0,1]\) choose \(P\subseteq A\) with \(\mu(P)=tp\) and let
  \(Q=A\setminus P\), so \(\mu(Q)=(1-t)p\).
\end{enumerate}

\emph{Comonotone construction on \(P\).} Choose \(F\subseteq P\) with
\(\mu(F)=\max(a,b)\,\mu(P)\) and then \(G\subseteq F\) with
\(\mu(G)=\min(a,b)\,\mu(P)\). If \(a\ge b\) put \(X=\mathbf{1}_F\),
\(Y=\mathbf{1}_G\) on \(P\); otherwise put \(X=\mathbf{1}_G\),
\(Y=\mathbf{1}_F\). Then \(\int_PX\,d\mu=a\,\mu(P)\),
\(\int_PY\,d\mu=b\,\mu(P)\) and
\(\int_PXY\,d\mu=\mu(G)=\min(a,b)\,\mu(P)\).

\emph{Countermonotone construction on \(Q\).} If \(a+b\le1\): choose
\(E\subseteq Q\) with \(\mu(E)=a\,\mu(Q)\) and then
\(E'\subseteq Q\setminus E\) with \(\mu(E')=b\,\mu(Q)\) (possible since
\(\mu(Q\setminus E)=(1-a)\mu(Q)\ge b\,\mu(Q)\)). If \(a+b>1\): choose
\(E\subseteq Q\) with \(\mu(E)=a\,\mu(Q)\), then \(G'\subseteq E\) with
\(\mu(G')=(a+b-1)\,\mu(Q)\) (note \(a+b-1\le a\)), and put
\(E'=(Q\setminus E)\cup G'\), so that
\(\mu(E')=(1-a+a+b-1)\mu(Q)=b\,\mu(Q)\). On \(Q\) put
\(X=\mathbf{1}_E\), \(Y=\mathbf{1}_{E'}\); then
\(\int_QX\,d\mu=a\,\mu(Q)\), \(\int_QY\,d\mu=b\,\mu(Q)\) and
\(\int_QXY\,d\mu=\mu(E\cap E')=\max(0,a+b-1)\,\mu(Q)\).

Together, \(X\) and \(Y\) take only the values \(0\) and \(1\) on \(A\),
have zone means \(a\) and \(b\), and \[
\mathbb{E}[XY\mid A]=t\min(a,b)+(1-t)\max(0,a+b-1).
\] Carry out this construction in every zone with the same \(t\). By the
formula in (i), \(\operatorname{Cov}(X,Y)=L+t\,(U-L)\). If \(U=L\) any
\(t\) works; otherwise \(t=(c-L)/(U-L)\in[0,1]\) gives
\(\operatorname{Cov}(X,Y)=c\). \(\blacksquare\)

\hypertarget{point-identification-and-the-width-of-the-interval}{%
\subsection{10.2 Point identification and the width of the
interval}\label{point-identification-and-the-width-of-the-interval}}

\textbf{Corollary 10.2 (point identification and interval width).} The
width of the interval is \[
U-L=\sum_j p_j\Big(\min(a_j,b_j)-\max(0,\,a_j+b_j-1)\Big)=\sum_jp_j\min\big(a_j,\ b_j,\ 1-a_j,\ 1-b_j\big),
\] and it vanishes if and only if every zone has \(a_j\in\{0,1\}\) or
\(b_j\in\{0,1\}\). If \(X,Y\) are indicators and \(0<\bar a,\bar b<1\),
the fine-scale variances \(\bar a(1-\bar a)\) and \(\bar b(1-\bar b)\)
are determined by the zonal data, and the fine-scale correlation lies in
the interval \[
\frac{[L,\,U]}{\sqrt{\bar a(1-\bar a)\,\bar b(1-\bar b)}} ;
\] if \(\mu\) is atomless, this interval is exactly the set of values of
the fine-scale correlation compatible with the zonal data.

\emph{Proof.} The first equality is the definition. The second is
checked term by term: if \(a+b\le1\), then \(a\le1-b\) and \(b\le1-a\),
and the term is \(\min(a,b)=\min(a,b,1-a,1-b)\); if \(a+b>1\), the term
is \(\min(a,b)-(a+b-1)=1-\max(a,b)=\min(1-a,1-b)\), and then \(1-b<a\)
and \(1-a<b\), so it also equals \(\min(a,b,1-a,1-b)\). By the proof of
Theorem 10.1(i) every term is non-negative, so the width vanishes if and
only if every term vanishes, i.e.~\(\min(a_j,b_j,1-a_j,1-b_j)=0\),
i.e.~\(a_j\in\{0,1\}\) or \(b_j\in\{0,1\}\). Indicators satisfy
\(X^2=X\), so
\(\operatorname{Var}X=\mathbb{E}X-(\mathbb{E}X)^2=\bar a(1-\bar a)\),
and likewise for \(Y\); by Theorem 10.1(i) the fine-scale correlation
lies in \([L,U]\) divided by this data-determined positive number, and
by Theorem 10.1(ii), when \(\mu\) is atomless, every value of this
interval occurs (the construction there uses indicators).
\(\blacksquare\)

Each zone contributes \(p_j\min(a_j,b_j,1-a_j,1-b_j)\) to the width: its
weight times the distance to \(\{0,1\}\) of whichever of its two zone
means is closer to \(0\) or \(1\). Zones whose means are close to \(0\)
or \(1\) nearly determine the fine-scale relation within them; zones
whose two means are both close to \(1/2\) say almost nothing. The width
\(U-L\) is computable from zonal data alone and measures how much of the
fine-scale relation the zonal data cannot determine; it can serve as an
index of the extent of the ecological fallacy. The ecological inference
method of King (1997) adds distributional assumptions within these
deterministic bounds; Theorem 10.1 shows that, for an atomless support
measure, without such assumptions no conclusion narrower than \([L,U]\)
is possible.

\textbf{Example.} Three zones of equal weight have zone means
\((a_j,b_j)=(0.2,0.3)\), \((0.5,0.5)\) and \((0.8,0.7)\). The three zone
means are collinear and the zonal correlation is \(\rho_B=+1\). One
computes \(\bar a=\bar b=0.5\), \(C_B=0.04\), \[
L=0.04+\tfrac13(-0.06-0.25-0.06)\approx-0.0833,\qquad U=0.04+\tfrac13(0.14+0.25+0.14)\approx0.2167,
\] and \(U-L=\tfrac13(0.2+0.5+0.2)=0.3\). For indicators, the fine-scale
correlation is only determined to lie in
\([L,U]/0.25\approx[-0.33,\ +0.87]\), and it may even be negative
(Figure 7(c)).

\begin{figure}
\centering
\includegraphics[width=1\textwidth,height=\textheight]{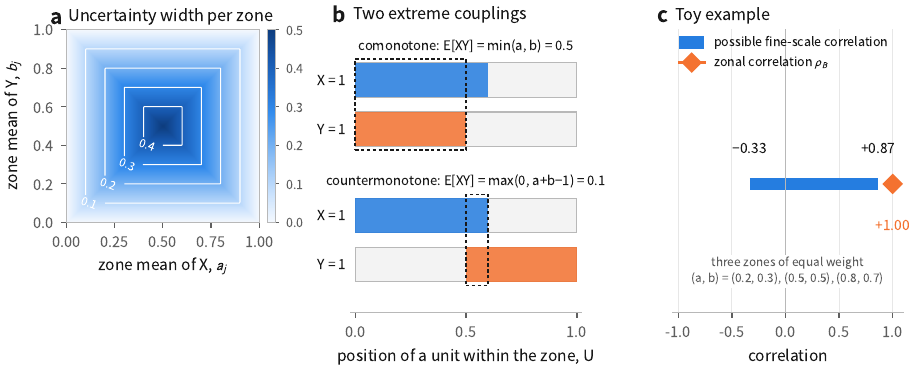}
\caption{Sharp identification bounds from zonal data alone. (a) The
contribution of a single zone, per unit zone weight, to the width
\(U-L\), \(\min(a,b)-\max(0,a+b-1)=\min(a,b,1-a,1-b)\), as a function of
the zone means \((a_j,b_j)\): zero on the boundary (either zone mean
equal to 0 or 1) and maximal, 0.5, at the centre (both zone means equal
to 1/2). (b) The two extreme couplings within one zone (\(a=0.6\),
\(b=0.5\)): with the fine-scale units of the zone arranged by position
\(U\), the comonotone coupling makes the units with \(X=1\) and those
with \(Y=1\) overlap as much as possible (\(\mathbb{E}[XY]=0.5\)), the
countermonotone coupling makes them overlap as little as possible
(\(\mathbb{E}[XY]=0.1\)); dashed boxes mark the overlap. Mixing the two
couplings in proportion gives every value in between (proof of Theorem
10.1). (c) Three zones of equal weight, \((a,b)=(0.2,0.3)\),
\((0.5,0.5)\), \((0.8,0.7)\): the zone means are collinear and the zonal
correlation is \(\rho_B=+1\), while the fine-scale correlation is only
determined to lie in \([-0.33,\ +0.87]\).}
\end{figure}

\hypertarget{interventions-on-zone-means-exact-linear-abstraction}{%
\section{11 Interventions on zone means: exact linear
abstraction}\label{interventions-on-zone-means-exact-linear-abstraction}}

When zone means are treated as macro variables, whether the macro
intervention ``set the zone mean \(\bar X\) to a given value'' has a
definite effect depends on which fine-scale interventions are permitted.
If arbitrary fine-scale interventions are permitted, the macro
intervention is well defined if and only if the within-zone effects are
homogeneous; if only uniform shifts within the zone are permitted, it is
always well defined, and its effect is the weighted average of the
within-zone effects (Theorem 11.2). For linear aggregation, this gives
an explicit condition of the kind studied under exact transformations
(Rubenstein et al.~2017) and abstractions (Beckers \& Halpern 2019) of
causal models.

\hypertarget{setting}{%
\subsection{11.1 Setting}\label{setting}}

Consider one zone containing fine-scale units \(i=1,\dots,n\) with
weights \(w_i>0\) (the mass of the support measure within the zone),
\(W=\sum_iw_i\). The fine-scale structural equations are \[
Y_i=\alpha_i+\beta_iX_i+U_i\qquad(i=1,\dots,n),
\] with no interference between units (the outcome of unit \(i\) is
affected only by \(X_i\)). The macro variables are the zone means
\(\bar X=\sum_iw_iX_i/W\) and \(\bar Y=\sum_iw_iY_i/W\). A fine-scale
intervention \(do(X=x)\) sets the vector \(X\) to \(x\in\mathbb{R}^n\);
the outcome is \(Y(x)\) and the corresponding zone means are
\(\bar X(x)\), \(\bar Y(x)\); \(\mathcal{I}\subseteq\mathbb{R}^n\) is
the class of permitted fine-scale interventions. The macro intervention
is \textbf{well defined} on \(\mathcal{I}\) if for all
\(x,x'\in\mathcal{I}\) and all \(U=(U_i)\), \(\bar X(x)=\bar X(x')\)
implies \(\bar Y(x)=\bar Y(x')\). For a zoning with several zones, the
results below apply zone by zone.

\textbf{Lemma 11.1 (proportional linear functionals).} Let \(S\) be a
real vector space and \(\varphi,\psi\) linear functionals on \(S\). Then
\(\varphi\) vanishes on \(\ker\psi\) if and only if \(\varphi=b\,\psi\)
for some \(b\in\mathbb{R}\). If \(\psi\neq0\) such \(b\) is unique; if
\(\psi=0\) the condition is \(\varphi=0\), and \(b\) is arbitrary.

\emph{Proof.} Sufficiency is clear. Necessity: if \(\psi=0\), then
\(\ker\psi=S\) and \(\varphi=0=b\psi\) for every \(b\). If
\(\psi\neq0\), choose \(s_0\in S\) with \(\psi(s_0)=1\). For any
\(s\in S\), \(s-\psi(s)s_0\in\ker\psi\), so
\(\varphi(s)=\psi(s)\varphi(s_0)\); take \(b=\varphi(s_0)\). If
\(\varphi=b\psi=b'\psi\), evaluating at \(s_0\) gives \(b=b'\).
\(\blacksquare\)

\hypertarget{the-condition-for-well-definedness}{%
\subsection{11.2 The condition for
well-definedness}\label{the-condition-for-well-definedness}}

\textbf{Theorem 11.2 (when a macro intervention is well defined).} Let
\(\mathcal{I}=x^0+S\) with \(x^0\in\mathbb{R}^n\) and
\(S\subseteq\mathbb{R}^n\) a linear subspace. The macro intervention is
well defined on \(\mathcal{I}\) if and only if there is
\(b\in\mathbb{R}\) with \[
\sum_i w_i\beta_i s_i=b\sum_i w_i s_i\qquad\text{for all }s\in S .\tag{11.1}
\] In that case \(\bar Y=c+b\bar X+\bar U\) on \(\mathcal{I}\), where
\(\bar U=\sum_iw_iU_i/W\) and \(c\) is a constant not depending on the
intervention; \(b\) is the macro causal effect (unique as soon as some
permitted intervention changes \(\bar X\)). In particular:

\begin{enumerate}
\def\labelenumi{(\alph{enumi})}
\item
  for \(S=\mathbb{R}^n\) (arbitrary fine-scale interventions), the macro
  intervention is well defined if and only if all \(\beta_i\) are equal;
\item
  for \(S=\operatorname{span}\{\mathbf{1}\}\) (uniform shifts within the
  zone), the macro intervention is always well defined, with
  \(b=\sum_iw_i\beta_i/W\), the weighted average of the within-zone
  effects.
\end{enumerate}

\emph{Proof.} For \(x,x'\in\mathcal{I}\) let \(s=x-x'\in S\). By the
structural equations, \[
\bar Y(x)-\bar Y(x')=\frac1W\sum_iw_i\beta_is_i=\frac{\varphi(s)}{W},\qquad \bar X(x)-\bar X(x')=\frac1W\sum_iw_is_i=\frac{\psi(s)}{W},
\] where \(\varphi(s)=\sum_iw_i\beta_is_i\) and
\(\psi(s)=\sum_iw_is_i\); \(U\) cancels in both differences. As \(x,x'\)
range over \(\mathcal{I}\), \(s\) ranges over \(S\), so well-definedness
is equivalent to: \(\varphi|_S\) vanishes on \(\ker(\psi|_S)\). By Lemma
11.1 (applied to the vector space \(S\)), this is equivalent to (11.1).
Taking \(x'=x^0\) gives
\(\bar Y(x)=\bar Y(x^0)+b\big(\bar X(x)-\bar X(x^0)\big)\); since
\(\bar Y(x^0)=\sum_iw_i(\alpha_i+\beta_ix^0_i)/W+\bar U\), we get
\(\bar Y=c+b\bar X+\bar U\) with
\(c=\sum_iw_i(\alpha_i+\beta_ix^0_i)/W-b\bar X(x^0)\).

\begin{enumerate}
\def\labelenumi{(\alph{enumi})}
\tightlist
\item
  For \(S=\mathbb{R}^n\), \(\psi\neq0\), and (11.1) with \(s=e_i\) reads
  \(w_i\beta_i=bw_i\) for all \(i\), i.e.~\(\beta_i=b\). (b) \(S\) is
  one-dimensional, \(\varphi(\mathbf{1})=\sum_iw_i\beta_i\) and
  \(\psi(\mathbf{1})=W\neq0\), so (11.1) holds with
  \(b=\sum_iw_i\beta_i/W\). \(\blacksquare\)
\end{enumerate}

\textbf{Example 11.3 (three implementations of the same macro
intervention).} A zone contains two units of equal weight with
\(\beta_1=2\) and \(\beta_2=0\). Three fine-scale interventions each
lower the zone mean \(\bar X\) by 1: lower both units by 1; lower only
unit 2 by 2; lower only unit 1 by 2. They lower \(\bar Y\) by 1, 0 and 2
respectively. So on \(S=\mathbb{R}^2\) ``the effect of lowering
\(\bar X\) by 1'' has no unique value (Theorem 11.2(a)); restricted to
uniform shifts within the zone it has the definite value \(1\), the
weighted average of \(\beta\) (Theorem 11.2(b)).

\hypertarget{zonal-effects-and-the-zoning}{%
\subsection{11.3 Zonal effects and the
zoning}\label{zonal-effects-and-the-zoning}}

\textbf{Corollary 11.4 (how zonal causal effects change with the
zoning).} Under the class of uniform within-zone shifts, regard
\(\beta\) as an integrable variable on the support measure space
\((\Omega,\mathcal{F},\mu)\). The zonal macro effect of any zoning
\(\{A_j\}\) is \(b_j=\mathbb{E}[\beta\mid A_j]\). The weighted average
over zones \(\sum_jp_jb_j=\mathbb{E}\beta\) does not depend on the
zoning; the distribution of zonal effects is dominated in the convex
order by that of the fine-scale effects, and its dispersion can change
with the zoning. If \(\mu\) is atomless, any statistic of the weighted
distribution of zonal effects that takes the same value for all
\(k\)-zone zonings must equal its value when every zonal effect equals
\(\mathbb{E}\beta\) (Theorem 5.5): only the average effect is a zoning
invariant.

\emph{Proof.} Theorem 11.2(b) extends verbatim from a zone of finitely
many units to a zone \(A\) with \(\mu(A)>0\) of a general support
measure space: with \(\beta\in L^1\), a uniform shift of \(X\) by \(t\)
on \(A\) changes \(\bar X=\mathbb{E}[X\mid A]\) by \(t\) and
\(\bar Y=\mathbb{E}[Y\mid A]\) by \(t\,\mathbb{E}[\beta\mid A]\). Hence
\(b_j=\mathbb{E}[\beta\mid A_j]\), the average of \(\beta\) over \(A_j\)
weighted by the support measure. By Lemma 3.1(d),
\(\sum_jp_jb_j=\mathbb{E}\beta\); by Corollary 3.6,
\(\mathbb{E}[\beta\mid\mathcal{G}]\le_{cx}\beta\). If \(\mu\) is
atomless and \(\beta\) is not almost surely constant, the dispersion
does change: under the 2-zone zoning of Corollary 4.9 both zonal effects
equal \(\mathbb{E}\beta\), whereas for \(t\) with \(0<\mu(\beta>t)<1\)
the zoning \(\{\beta>t\},\{\beta\le t\}\) has \(b_1>t\ge b_2\). The last
sentence is Theorem 5.5 applied to \(Z=\beta\) (\(d=1\)).
\(\blacksquare\)

The zonal regression slope of \(\bar Y_j\) on \(\bar X_j\) across zones
generally equals no macro causal effect: when \(\bar U_j\) is correlated
with \(\bar X_j\) across zones it is confounded, which is exactly the
bias discussed in Corollary 6.4 (data-dependent zonings can produce any
slope) and Theorem 8.14 (aggregation bias under data-independent
designs). Conditions for exact abstraction with spatial interference
between units or nonlinear fine-scale models are beyond the scope of
this paper (Section 13.3).

\hypertarget{numerical-checks}{%
\section{12 Numerical checks}\label{numerical-checks}}

All checks are carried out by the script
\texttt{scale\_theory\_verify.py}, provided as an ancillary file (Python
with NumPy and SciPy, fixed random seed 20261126; the results are
written to \texttt{verify\_results.json}). Identities should hold to
machine precision; existence results (Theorems 6.2 and 6.5) hold exactly
only for atomless support measures and can only be approximated on
finite samples or rasters, so what is checked is whether the error
shrinks as the sample grows or the grid is refined. Numerical checks do
not replace proofs; their role is to rule out errors in signs, constants
and conditions. Changing the random seed changes only the specific
numbers, not any of the conclusions.

\begin{longtable}[]{@{}
  >{\raggedright\arraybackslash}p{(\columnwidth - 6\tabcolsep) * \real{0.0488}}
  >{\raggedright\arraybackslash}p{(\columnwidth - 6\tabcolsep) * \real{0.1707}}
  >{\raggedright\arraybackslash}p{(\columnwidth - 6\tabcolsep) * \real{0.2683}}
  >{\raggedright\arraybackslash}p{(\columnwidth - 6\tabcolsep) * \real{0.5122}}@{}}
\caption{Summary of numerical checks}\tabularnewline
\toprule\noalign{}
\begin{minipage}[b]{\linewidth}\raggedright
ID
\end{minipage} & \begin{minipage}[b]{\linewidth}\raggedright
Result checked
\end{minipage} & \begin{minipage}[b]{\linewidth}\raggedright
Setting
\end{minipage} & \begin{minipage}[b]{\linewidth}\raggedright
Outcome
\end{minipage} \\
\midrule\noalign{}
\endfirsthead
\toprule\noalign{}
\begin{minipage}[b]{\linewidth}\raggedright
ID
\end{minipage} & \begin{minipage}[b]{\linewidth}\raggedright
Result checked
\end{minipage} & \begin{minipage}[b]{\linewidth}\raggedright
Setting
\end{minipage} & \begin{minipage}[b]{\linewidth}\raggedright
Outcome
\end{minipage} \\
\midrule\noalign{}
\endhead
\bottomrule\noalign{}
\endlastfoot
V0 & Theorem 5.2, Corollary 5.3 & 5000 fine-scale units, 37 random
zones, unequal weights & identity residual \(5.6\times10^{-17}\),
slope-decomposition residual 0 \\
V1 & Theorem 6.2, Corollary 6.3 (no connectivity, \(k=3\)) &
\(2\times10^6\) fine-scale units, fine-scale correlation 0.600;
randomized rounding of the fuzzy partition & targets \(-0.95\),
\(-0.5\), \(0\), \(0.5\), \(0.95\) achieved as \(-0.943\), \(-0.512\),
\(0.013\), \(0.501\), \(0.952\); for \(k=2\) always \(\pm1\) \\
V2 & Theorem 6.5 (connected zones, serpentine-lane construction) &
\(1440\times1440\) raster, \(120\times120\) fine cells, \(30\times30\)
coarse cells, fine-scale correlation 0.878 & targets \(-0.9\), \(-0.5\),
\(0\), \(0.5\), \(0.9\) achieved as \(-0.907\), \(-0.483\), \(0.009\),
\(0.512\), \(0.901\); every zone is a single 4-connected component;
refining the fine cells from 40 to 80 and 160 per side reduces the
maximal error from 0.075 to 0.022 and 0.011 \\
V3 & Theorem 7.1, Corollary 7.2, Proposition 7.4 & \(512\times512\)
raster; square blocks of side 2--128 pixels, random Voronoi zonings with
50--3000 zones & all deviations below the bound of Theorem 7.1 (16-pixel
blocks: deviation \(1.1\times10^{-3}\), bound \(8.7\times10^{-3}\)); the
bound of Corollary 7.2 is looser (0.35 at 16 pixels) and uninformative
for \(\tau\ge1\); the construction of Proposition 7.4 attains the bound
(both 0.47145) \\
V4 & Lemma 8.6, Theorem 8.9, Theorem 8.15 & Gaussian fields
\(X=S+M+N_1\), \(Y=S-M+N_2\) on a \(256\times256\) torus, 1500
simulations per scale & theory and simulation agree within 2 standard
errors; the whole curve lies within the coherence bounds
\([-0.98,\ 0.91]\); as the block side grows from 1 to 256 pixels the
zonal correlation moves from \(-0.40\) to
\(+0.91=\operatorname{coh}(0)\) (Figure 5(f)) \\
V5 & Theorem 8.11 (the intrinsic correlation model implies scale
invariance) & same spectrum, \(Y=0.6A+0.8B\) & all 13 designs (8 square
blocks, 4 discs and one elongated rectangle) equal 0.6000, error below
\(10^{-15}\) \\
V6 & Theorems 9.2--9.4 and Example 9.5 & \(256\times256\) raster, 8
dyadic block levels & cross-scale products below \(10^{-16}\),
decomposition residuals below \(2\times10^{-16}\), the hierarchical
coherence bound holds; Example 9.5 gives \(\rho_0=-0.1386\),
\(\rho_1=-0.5000\) \\
V7 & Theorem 10.1 (sharp bounds) & 200 random \([0,1]\)-valued
configurations with 25 zones each & no violation; the comonotone and
countermonotone constructions attain \(U\) and \(L\) exactly \\
V8 & Theorem 11.2 & 50 units with heterogeneous \(\beta_i\) & under
uniform shifts the macro effect equals the weighted mean of \(\beta\)
(both 0.98339); an arbitrary intervention preserving \(\bar X\) changes
\(\bar Y\) by 0.028 \\
V9 & Theorem 8.14 (aggregation bias iff
\(\operatorname{Re}F_{Xe}\neq0\)) & \(\operatorname{Cov}(X,e)=0\) at the
same location (numerically \(-7\times10^{-16}\)), \(\beta=1\) & as
blocks grow from 1 to 64 pixels the aggregated slope rises from 1.000 to
1.342; theory and simulation differ by less than one standard error \\
\end{longtable}

In V4 the two theoretical curves (square blocks and discs of equal area)
coincide at both ends and differ by at most 0.013 at intermediate
scales: the micro and macro limits do not depend on the zone shape,
which affects only intermediate scales (Theorem 8.15). The whole curve
lies within the coherence bounds (Theorem 8.9), and the sign change
occurs between 8 and 16 pixels. In V3 the bound of Theorem 7.1 lies
above the actual deviation at every scale and grows, like the deviation,
roughly as the square of the zone scale (Figure 4(e), Corollary 7.3);
the Lipschitz bound of Corollary 7.2 is larger by more than an order of
magnitude, as explained in Section 7.2. In V2 the error of the connected
construction decreases roughly as \(1/N\) as the fine cells are refined
(Figure 8(b)).

\begin{figure}
\centering
\includegraphics[width=1\textwidth,height=\textheight]{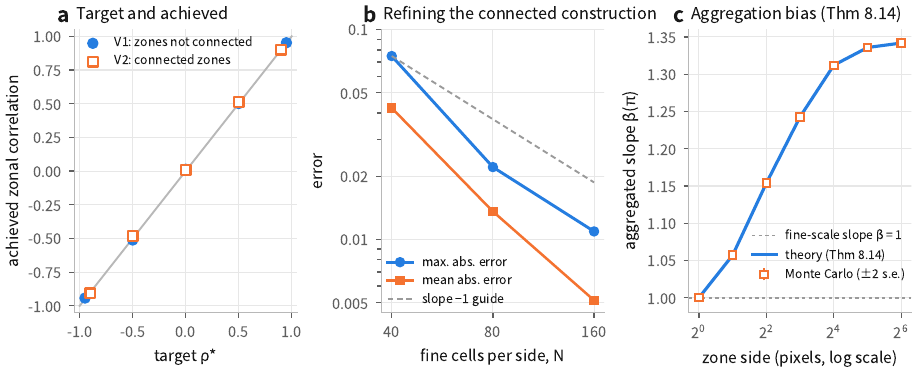}
\caption{Numerical checks of the existence theorems and of aggregation
bias. (a) Target versus achieved zonal correlations in checks V1 (no
connectivity, blue dots) and V2 (every zone connected, open orange
squares); grey line, the diagonal. (b) Grid refinement of the connected
construction in V2: as the number \(N\) of fine cells per side is
refined from 40 to 80 and 160, the maximal (blue) and mean (orange)
absolute errors over the five targets decrease roughly as \(1/N\); grey
dashed line, reference slope \(-1\) (log--log axes). (c) Check V9:
\(\operatorname{Cov}(X,e)=0\) at the same location and the fine-scale
slope is \(\beta=1\) (grey dashed line), yet the aggregated regression
slope rises to 1.34 as the block side grows; blue line, the theoretical
value of Theorem 8.14; orange squares, simulation results (\(\pm2\)
standard errors).}
\end{figure}

\hypertarget{discussion}{%
\section{13 Discussion}\label{discussion}}

\hypertarget{main-contributions}{%
\subsection{13.1 Main contributions}\label{main-contributions}}

This paper places the scale problem of zonal statistics in a single
framework: a zoning of a support measure space is a
sub-\(\sigma\)-algebra, and aggregation is conditional expectation.
Within this framework, the zoning effect, the support-scale effect and
multiscale association each receive a precise answer. In the language of
spatial statistics, the core conclusions come down to three statements.

\textbf{First, the zoning effect: data-dependent zonings are not
constrained by fine-scale relations.} Under an atomless support measure,
as long as there is no exact linear relation between the variables, for
a fixed number of zones \(k\ge3\) and any prescribed zone weights,
data-dependent zonings make the zonal correlation take every value in
\([-1,1]\) exactly; the zonal correlation matrix ranges over exactly the
correlation matrices of rank at most \(k-1\) (Theorem 6.2, Corollary
6.3), and the zonal regression slope can take any real value (Corollary
6.4). Requiring every zone to be connected in the unit square still lets
the zonal correlation approach any value in \([-1,1]\) arbitrarily
closely (Theorem 6.5, Corollary 6.6). The connected zones constructed
are long thin lanes running across the domain, so what really constrains
the zonal correlation is compactness, not connectivity. This gives the
limiting form of the phenomenon observed in simulation studies such as
Openshaw \& Taylor (1979): when the zoning may be drawn by looking at
the data, zonal association by itself is no evidence about fine-scale
relations.

\textbf{Second, support scale and aggregation bias: how well the error
can be controlled.} ``Smaller zones, smaller error'' is only a
direction. The contribution here is an error analysis: the deviation of
the zonal correlation from the fine-scale correlation is bounded
explicitly in terms of the fine-scale correlation and the two
within-zone variance shares \(a\) and \(b\) (Theorem 7.1). The bound can
be computed from fine-scale data alone, without any model assumption,
and it is attained whenever it does not exceed \(1-|\rho|\) (Proposition
7.4). For Lipschitz fields the within-zone variance shares are
controlled by the square of the ratio of zone diameter to variation
length, and the leading term of the error bound is of second order,
\(\kappa\delta^2\), in the largest zone diameter \(\delta\): halving the
largest diameter reduces the leading term to one quarter (Corollary
7.3), and the second order cannot be improved (Example 7.5). The rate
depends on the smoothness of the fields: for mean-square differentiable
stationary fields the design-level within-zone variance shares, and
hence the bound, again vanish at second order, with an exact constant
given by the spectral second moments, whereas for rough fields with
\(C(0)-C(h)\sim c|h|^\alpha\), \(\alpha<2\), they vanish only at order
\(\alpha\) (first order for the exponential covariance), and the
deviation itself can then be of first order (Proposition 8.16 and the
example after it). In applications one should compute \(a\) and \(b\)
from the data and use Theorem 7.1; the role of the Lipschitz form is to
explain why the error is of second order.

\textbf{Third, multiscale association: a necessary and sufficient
condition for scale invariance.} For data-independent zoning designs,
``the correlation may change with scale'' is common knowledge; what this
paper answers is under which condition the correlation stays the same
across scales and zone shapes. The design-level zonal correlation takes
the same non-zero value for all connected zone shapes and their finite
mixtures if and only if the power spectra of the two fields are
proportional and the co-spectrum (the real part of the cross-spectrum)
is a non-zero multiple of the common power spectrum, which is the
intrinsic correlation model of multivariate geostatistics (Theorem
8.11(a)); it vanishes for all of them if and only if the co-spectrum
vanishes (Theorem 8.11(b)). When this condition fails and the fine-scale
correlation is non-zero, some design made of at most two connected zone
shapes must have a zonal correlation different from the fine-scale
correlation (Corollary 8.12). Conversely, observing similar or even
identical zonal correlations at a few chosen scales or zone shapes does
not imply the proportionality (Proposition 8.13): testing scale
invariance requires testing whether the spectra are proportional, not
comparing correlations at finitely many scales.

Around the third statement the paper also shows: the zonal correlation
is a weighted average of the real coherence times a factor at most 1,
and therefore always lies in
\([\min(0,\operatorname{ess\,inf}\operatorname{coh}),\)
\(\max(0,\operatorname{ess\,sup}\operatorname{coh})]\) (Theorem 8.9), in
sharp contrast with the first statement; as zones grow without bound,
the zonal correlation converges to the zero-frequency coherence,
whatever the zone shape, for square-integrable kernels and spectral
densities that are bounded, continuous and positive at the origin
(Theorem 8.15); and the aggregated regression slope is unbiased for all
connected designs if and only if the symmetrized cross-covariance
between the explanatory variable and the error vanishes at all lags, so
that being uncorrelated at the same location is not enough (Theorem
8.14).

\textbf{Other results.} Under an atomless support measure, a zonal
statistic that takes the same value for all \(k\)-zone zonings must take
the value it has when every zone mean equals the overall mean (Theorem
5.5). On nested hierarchies the covariance decomposes orthogonally
across scales, and the zonal correlation at a level depends only on that
level and the coarser ones; ``coarser aggregation, higher correlation''
holds when the scale spectra are proportional and the level coherence is
non-decreasing from finer to coarser levels, and can fail when the
spectra are not proportional (Theorems 9.2--9.4, Example 9.5); the fixed
points of the coarsening map are exactly the geometric scale spectra
(Proposition 9.6). From zonal data alone, the fine-scale covariance of
\([0,1]\)-valued variables lies in a computable interval that is sharp
under an atomless support measure (Theorem 10.1, Corollary 10.2). A
macro intervention setting zone means to given values is well defined if
and only if the within-zone effects are homogeneous along the permitted
intervention directions (Theorem 11.2).

\textbf{Relation to known results.} The results of Section 3 are
standard; Theorem 4.3 (Sierpiński 1922), Theorem 4.6 (Lyapunov 1940),
Theorem 4.8 and its fixed-weight form Corollary 4.10 (Dvoretzky, Wald \&
Wolfowitz 1951), Lemma 5.1 and Theorem 5.2 (an equivalent form of
Robinson 1950), Theorems 8.1--8.3 (the uniqueness theorem and the
theorems of Bochner and Cramér 1940), Lemma 8.6 (the block-covariance
formula; Journel \& Huijbregts 1978; Cressie 1993; Gotway \& Young
2002), the sufficiency direction of Theorem 8.11 (Wackernagel 2003),
Lemma 9.1 and Theorem 9.2 (orthogonality of martingale increments) and
the binary case of Theorem 10.1 (Duncan \& Davis 1953; Fréchet 1951) are
known; Theorem 8.15(ii) echoes Tiao (1972) and Vargas-Guzmán, Warrick \&
Myers (1999). The other main results, namely Theorem 5.5, Theorem 6.2,
Corollaries 6.3, 6.4 and 6.6, Theorem 6.5, Theorem 7.1, Corollaries 7.2
and 7.3, Proposition 7.4, Theorem 8.9, the necessity direction of
Theorem 8.11, Corollary 8.12, Proposition 8.13, Theorem 8.14,
Proposition 8.16, Theorems 9.3 and 9.4, Propositions 9.6 and 9.7, the
general case of Theorem 10.1, Corollary 10.2, Theorem 11.2 and Corollary
11.4, are established in this paper.

\hypertarget{implications-for-spatial-analysis}{%
\subsection{13.2 Implications for spatial
analysis}\label{implications-for-spatial-analysis}}

Faced with the question ``will this zonal conclusion change with the
zoning, and by how much?'', one can proceed in the following order; each
step corresponds to a proved result.

\begin{enumerate}
\def\labelenumi{\arabic{enumi}.}
\tightlist
\item
  \textbf{Is the statistic a function of the overall mean only?} If so,
  it does not depend on the zoning (Lemma 3.1). Otherwise, under an
  atomless support measure it can be invariant only by always taking the
  value it has when every zone mean equals the overall mean (Theorem
  5.5); for variables without an exact linear relation, zonal variances,
  correlations and slopes are not invariant (Section 5.3, Corollaries
  6.3 and 6.4).
\item
  \textbf{Was the zoning drawn with reference to the data, or selected
  from many candidate zonings?} If so, the zonal correlation and the
  zonal regression slope can take any value (Theorem 6.2, Corollaries
  6.3 and 6.4), and requiring connected zones does not prevent this
  (Theorem 6.5). Fine-scale data should be used instead, or the
  distribution of results over many admissible zonings should be
  reported (for example by sampling zonings under connectivity
  constraints; DeFord, Duchin \& Solomon 2021).
\item
  \textbf{Are fine-scale data available?} If so, compute the two
  within-zone variance shares \(a\) and \(b\) and substitute them into
  Theorem 7.1 to obtain how far the zonal correlation can deviate from
  the fine-scale correlation. If both shares exceed one half, the bound
  is uninformative.
\item
  \textbf{Is the zoning independent of the data, and can the two
  variables be regarded as jointly stationary fields?} If so, examine
  the cross-spectrum: the zonal correlation lies within the coherence
  bounds of Theorem 8.9; it takes the same non-zero value for all
  connected designs if and only if the intrinsic correlation model holds
  (Theorem 8.11); if it does not hold and the fine-scale correlation is
  non-zero, some design of at most two connected zone shapes deviates
  from it, and if the fine-scale correlation is zero but the co-spectrum
  does not vanish, some such design has a non-zero zonal correlation
  (Corollary 8.12); comparisons at finitely many scales cannot replace a
  test of the spectra (Proposition 8.13); very large zones estimate the
  zero-frequency coherence \(\operatorname{coh}(0)\) (Theorem 8.15); the
  aggregated regression slope is unbiased if and only if the symmetrized
  lagged cross-covariance vanishes (Theorem 8.14).
\item
  \textbf{Is there a nested hierarchy?} If so, decompose the covariance
  across scales by Theorem 9.2 and report \(v_l\), \(c_l\) and
  \(\kappa_l\) at each level; ``the more aggregated, the stronger the
  correlation'' is guaranteed when the scale spectra are proportional
  and \(\kappa_l\) is non-decreasing from finer to coarser levels
  (Theorem 9.4); without proportionality it can fail even when
  \(\kappa_l\) increases (Example 9.5).
\item
  \textbf{Are only zonal data available, while the target is the
  fine-scale relation?} First report the bounds \([L,U]\), which are
  sharp under an atomless support measure, and the width \(U-L\)
  (Theorem 10.1, Corollary 10.2), and then decide whether to introduce
  distributional assumptions (King 1997).
\item
  \textbf{Is a zonal result to be interpreted as an intervention
  effect?} First state the permitted class of interventions: under
  uniform within-zone shifts the macro effect is well defined and equals
  the weighted average of the within-zone effects (Theorem 11.2);
  extending it to arbitrary fine-scale interventions requires
  homogeneous within-zone effects.
\end{enumerate}

These steps apply to any variable summarized over zones. In
public-health monitoring, both the choice of area-based socioeconomic
measure and the geographic level affect the measured health inequalities
(Krieger et al.~2002), but the results here are not limited to
population: what is aggregated within zones may be population counts and
health outcomes, but also land-cover proportions, pollutant
concentrations, incomes or remotely sensed retrievals; the support
measure is correspondingly population, area or sampling weights (Section
3.1). For rate-type variables the support measure is usually the count
of the denominator, for example population; switching to area weights or
equal weights per zone changes the meaning of the ``fine-scale
correlation'' itself.

\hypertarget{limitations}{%
\subsection{13.3 Limitations}\label{limitations}}

\textbf{The atomless idealization.} The existence results of Sections
4--6 and the sharpness statement of Theorem 10.1(ii) hold exactly under
atomless support measures. Actual census units or rasters are finite,
and each unit can only be assigned to a zone as a whole, so the results
hold only approximately; checks V1 and V2 show that the approximation
error shrinks as the units are refined, but quantitative results for the
discrete case have not been proved.

\textbf{Stationarity.} Section 8 and Proposition 9.7 assume that the two
fields are jointly weakly stationary and ignore boundary effects of the
study region. For non-stationary fields (for example with large-scale
trends or local heterogeneity), the spectral matrix must be replaced by
local spectra or non-stationary covariances, and the corresponding
necessary and sufficient conditions need to be reformulated.

\textbf{Second moments and linear models.} The results concern
correlation coefficients, covariances and linear regression slopes.
Apart from the convex-order monotonicity of Proposition 3.7 and the
impossibility theorem 5.5, results corresponding to Sections 6--9 for
nonlinear zonal statistics such as the Gini coefficient, segregation
indices or Moran's index have not been established; Section 11 treats
only linear structural equations without interference.

\textbf{Design level versus sample level.} The zonal correlation
\(\rho(\pi)\) of Section 8 is a design-level quantity, the expectation
of the zonal second moments under the law of the fields. The sample
zonal correlation from one realization on a bounded study region
fluctuates around it; we give no quantitative bound on this fluctuation.

\hypertarget{use-of-generative-ai}{%
\section{Use of generative AI}\label{use-of-generative-ai}}

The author directed the research, guided the development of the
manuscript and its proofs, and reviewed the completed work. Under the
author's direction, Claude (Anthropic), a large language model, provided
substantial assistance in drafting the text and proofs, producing the
figures, writing the numerical verification script
\texttt{scale\_theory\_verify.py}, and checking the references. The
author takes full responsibility for the final manuscript and all claims
it contains.

\hypertarget{references}{%
\section{References}\label{references}}

\begin{itemize}
\tightlist
\item
  Amrhein, C. G. (1995). Searching for the elusive aggregation effect:
  evidence from statistical simulations. \emph{Environment and Planning
  A}, 27(1), 105--119. \url{https://doi.org/10.1068/a270105}
\item
  Arbia, G. (1989). \emph{Spatial Data Configuration in Statistical
  Analysis of Regional Economic and Related Problems}. Kluwer,
  Dordrecht.
\item
  Arcaute, E., Hatna, E., Ferguson, P., Youn, H., Johansson, A., \&
  Batty, M. (2015). Constructing cities, deconstructing scaling laws.
  \emph{Journal of the Royal Society Interface}, 12(102), 20140745.
  \url{https://doi.org/10.1098/rsif.2014.0745}
\item
  Banerjee, S., Carlin, B. P., \& Gelfand, A. E. (2014).
  \emph{Hierarchical Modeling and Analysis for Spatial Data} (2nd ed.).
  CRC Press, Boca Raton, FL.
\item
  Beckers, S., \& Halpern, J. Y. (2019). Abstracting causal models.
  \emph{Proceedings of the AAAI Conference on Artificial Intelligence},
  33(01), 2678--2685. \url{https://doi.org/10.1609/aaai.v33i01.33012678}
\item
  Bradley, J. R., Wikle, C. K., \& Holan, S. H. (2017). Regionalization
  of multiscale spatial processes by using a criterion for spatial
  aggregation error. \emph{Journal of the Royal Statistical Society:
  Series B}, 79(3), 815--832. \url{https://doi.org/10.1111/rssb.12179}
\item
  Chilès, J.-P., \& Delfiner, P. (2012). \emph{Geostatistics: Modeling
  Spatial Uncertainty} (2nd ed.). Wiley, Hoboken, NJ.
\item
  Cramér, H. (1940). On the theory of stationary random processes.
  \emph{Annals of Mathematics}, 41(1), 215--230.
  \url{https://doi.org/10.2307/1968827}
\item
  Cressie, N. (1993). \emph{Statistics for Spatial Data} (rev. ed.).
  Wiley, New York.
\item
  Cressie, N. (1996). Change of support and the modifiable areal unit
  problem. \emph{Geographical Systems}, 3(2--3), 159--180.
\item
  Daw, R., Bradley, J. R., Wikle, C. K., \& Holan, S. H. (2025). A
  criterion for aggregation error for multivariate spatial data.
  arXiv:2312.12287 (v3). \url{https://arxiv.org/abs/2312.12287}
\item
  DeFord, D., Duchin, M., \& Solomon, J. (2021). Recombination: a family
  of Markov chains for redistricting. \emph{Harvard Data Science
  Review}, 3(1). \url{https://doi.org/10.1162/99608f92.eb30390f}
\item
  Duncan, O. D., \& Davis, B. (1953). An alternative to ecological
  correlation. \emph{American Sociological Review}, 18(6), 665--666.
  \url{https://doi.org/10.2307/2088122}
\item
  Duque, J. C., Anselin, L., \& Rey, S. J. (2012). The max-p-regions
  problem. \emph{Journal of Regional Science}, 52(3), 397--419.
  \url{https://doi.org/10.1111/j.1467-9787.2011.00743.x}
\item
  Durrett, R. (2019). \emph{Probability: Theory and Examples} (5th ed.).
  Cambridge University Press, Cambridge.
\item
  Dvoretzky, A., Wald, A., \& Wolfowitz, J. (1951). Relations among
  certain ranges of vector measures. \emph{Pacific Journal of
  Mathematics}, 1(1), 59--74.
  \url{https://doi.org/10.2140/pjm.1951.1.59}
\item
  Elton, J., \& Hill, T. P. (1992). Fusions of a probability
  distribution. \emph{Annals of Probability}, 20(1), 421--454.
  \url{https://doi.org/10.1214/aop/1176989936}
\item
  Folland, G. B. (1999). \emph{Real Analysis: Modern Techniques and
  Their Applications} (2nd ed.). Wiley.
\item
  Fotheringham, A. S., \& Wong, D. W. S. (1991). The modifiable areal
  unit problem in multivariate statistical analysis. \emph{Environment
  and Planning A}, 23(7), 1025--1044.
  \url{https://doi.org/10.1068/a231025}
\item
  Fréchet, M. (1951). Sur les tableaux de corrélation dont les marges
  sont données. \emph{Annales de l'Université de Lyon, Section A, Série
  3}, 14, 53--77.
\item
  Gehlke, C. E., \& Biehl, K. (1934). Certain effects of grouping upon
  the size of the correlation coefficient in census tract material.
  \emph{Journal of the American Statistical Association}, 29(185A),
  169--170. \url{https://doi.org/10.1080/01621459.1934.10506247}
\item
  Gelfand, A. E., Zhu, L., \& Carlin, B. P. (2001). On the change of
  support problem for spatio-temporal data. \emph{Biostatistics}, 2(1),
  31--45. \url{https://doi.org/10.1093/biostatistics/2.1.31}
\item
  Gelman, A., Shor, B., Bafumi, J., \& Park, D. (2007). Rich state, poor
  state, red state, blue state: What's the matter with Connecticut?
  \emph{Quarterly Journal of Political Science}, 2(4), 345--367.
  \url{https://doi.org/10.1561/100.00006026}
\item
  Goodman, L. A. (1953). Ecological regressions and behavior of
  individuals. \emph{American Sociological Review}, 18(6), 663--664.
  \url{https://doi.org/10.2307/2088121}
\item
  Gotway, C. A., \& Young, L. J. (2002). Combining incompatible spatial
  data. \emph{Journal of the American Statistical Association}, 97(458),
  632--648. \url{https://doi.org/10.1198/016214502760047140}
\item
  Greenland, S., \& Robins, J. (1994). Invited commentary: ecologic
  studies --- biases, misconceptions, and counterexamples.
  \emph{American Journal of Epidemiology}, 139(8), 747--760.
  \url{https://doi.org/10.1093/oxfordjournals.aje.a117069}
\item
  Guan, Y., Page, G. L., Reich, B. J., Ventrucci, M., \& Yang, S.
  (2023). Spectral adjustment for spatial confounding.
  \emph{Biometrika}, 110(3), 699--719.
  \url{https://doi.org/10.1093/biomet/asac069}
\item
  Hodges, J. S., \& Reich, B. J. (2010). Adding spatially-correlated
  errors can mess up the fixed effect you love. \emph{The American
  Statistician}, 64(4), 325--334.
  \url{https://doi.org/10.1198/tast.2010.10052}
\item
  Holt, D., Steel, D. G., \& Tranmer, M. (1996). Area homogeneity and
  the modifiable areal unit problem. \emph{Geographical Systems},
  3(2--3), 181--200.
\item
  Journel, A. G., \& Huijbregts, C. J. (1978). \emph{Mining
  Geostatistics}. Academic Press, London.
\item
  King, G. (1997). \emph{A Solution to the Ecological Inference Problem:
  Reconstructing Individual Behavior from Aggregate Data}. Princeton
  University Press.
\item
  Krieger, N., Chen, J. T., Waterman, P. D., Soobader, M.-J.,
  Subramanian, S. V., \& Carson, R. (2002). Geocoding and monitoring of
  US socioeconomic inequalities in mortality and cancer incidence: does
  the choice of area-based measure and geographic level matter? The
  Public Health Disparities Geocoding Project. \emph{American Journal of
  Epidemiology}, 156(5), 471--482.
  \url{https://doi.org/10.1093/aje/kwf068}
\item
  Lindenstrauss, J. (1966). A short proof of Liapounoff's convexity
  theorem. \emph{Journal of Mathematics and Mechanics}, 15, 971--972.
\item
  Lyapunov, A. A. (1940). Sur les fonctions-vecteurs complètement
  additives. \emph{Bulletin de l'Académie des Sciences de l'URSS, Série
  Mathématique} (Izvestiya Akademii Nauk SSSR, Seriya Matematicheskaya),
  4(6), 465--478 (in Russian).
\item
  Openshaw, S. (1977). A geographical solution to scale and aggregation
  problems in region-building, partitioning and spatial modelling.
  \emph{Transactions of the Institute of British Geographers}, 2(4),
  459--472. \url{https://doi.org/10.2307/622300}
\item
  Openshaw, S. (1984). \emph{The Modifiable Areal Unit Problem}. CATMOG
  38. Geo Books, Norwich.
\item
  Openshaw, S., \& Taylor, P. J. (1979). A million or so correlation
  coefficients: three experiments on the modifiable areal unit problem.
  In N. Wrigley (Ed.), \emph{Statistical Applications in the Spatial
  Sciences} (pp.~127--144). Pion, London.
\item
  Paciorek, C. J. (2010). The importance of scale for
  spatial-confounding bias and precision of spatial regression
  estimators. \emph{Statistical Science}, 25(1), 107--125.
  \url{https://doi.org/10.1214/10-STS326}
\item
  Popoviciu, T. (1935). Sur les équations algébriques ayant toutes leurs
  racines réelles. \emph{Mathematica (Cluj)}, 9, 129--145.
\item
  Robinson, W. S. (1950). Ecological correlations and the behavior of
  individuals. \emph{American Sociological Review}, 15(3), 351--357.
  \url{https://doi.org/10.2307/2087176}
\item
  Rockafellar, R. T. (1970). \emph{Convex Analysis}. Princeton
  University Press, Princeton, NJ.
\item
  Rubenstein, P. K., Weichwald, S., Bongers, S., Mooij, J. M., Janzing,
  D., Grosse-Wentrup, M., \& Schölkopf, B. (2017). Causal consistency of
  structural equation models. \emph{Proceedings of the 33rd Conference
  on Uncertainty in Artificial Intelligence (UAI 2017)}.
  \url{https://auai.org/uai2017/proceedings/papers/11.pdf}
\item
  Rudin, W. (1987). \emph{Real and Complex Analysis} (3rd ed.).
  McGraw-Hill.
\item
  Rudin, W. (1991). \emph{Functional Analysis} (2nd ed.). McGraw-Hill,
  New York.
\item
  Shaked, M., \& Shanthikumar, J. G. (2007). \emph{Stochastic Orders}.
  Springer, New York. \url{https://doi.org/10.1007/978-0-387-34675-5}
\item
  Sierpiński, W. (1922). Sur les fonctions d'ensemble additives et
  continues. \emph{Fundamenta Mathematicae}, 3, 240--246.
  \url{https://doi.org/10.4064/fm-3-1-240-246}
\item
  Simpson, E. H. (1951). The interpretation of interaction in
  contingency tables. \emph{Journal of the Royal Statistical Society:
  Series B}, 13(2), 238--241.
  \url{https://doi.org/10.1111/j.2517-6161.1951.tb00088.x}
\item
  Stein, M. L. (1999). \emph{Interpolation of Spatial Data: Some Theory
  for Kriging}. Springer, New York.
\item
  Theil, H. (1954). \emph{Linear Aggregation of Economic Relations}.
  North-Holland, Amsterdam.
\item
  Tiao, G. C. (1972). Asymptotic behaviour of temporal aggregates of
  time series. \emph{Biometrika}, 59(3), 525--531.
  \url{https://doi.org/10.1093/biomet/59.3.525}
\item
  Tobler, W. R. (1989). Frame independent spatial analysis. In M. F.
  Goodchild \& S. Gopal (Eds.), \emph{The Accuracy of Spatial Databases}
  (pp.~115--122). Taylor \& Francis, London.
\item
  Vargas-Guzmán, J. A., Warrick, A. W., \& Myers, D. E. (1999).
  Multivariate correlation in the framework of support and spatial
  scales of variability. \emph{Mathematical Geology}, 31(1), 85--103.
  \url{https://doi.org/10.1023/A:1007545431115}
\item
  Wackernagel, H. (2003). \emph{Multivariate Geostatistics: An
  Introduction with Applications} (3rd ed.). Springer.
\item
  Wakefield, J. (2008). Ecologic studies revisited. \emph{Annual Review
  of Public Health}, 29, 75--90.
  \url{https://doi.org/10.1146/annurev.publhealth.29.020907.090821}
\item
  Yaglom, A. M. (1987). \emph{Correlation Theory of Stationary and
  Related Random Functions, Vol. I: Basic Results}. Springer, New York.
\end{itemize}

\end{document}